\documentclass[11pt]{article}
\usepackage{amsfonts,amsmath,amssymb}
\usepackage{graphicx}
\usepackage{hyperref}
\usepackage{color}
\usepackage{subcaption}

\usepackage[normalem]{ulem}
\usepackage{tikz}
\usetikzlibrary{calc,arrows,positioning,fit,petri,intersections}

\newcommand{\plaat}[2]{{#2}}

\newcommand{\detail}[1]{\par\noi{\bf[Proof detail\ }{#1}
\hfill{\bf ]}\par\noi\hspace{-4pt}}
\renewcommand{\detail}[1]{}

\newcommand{\dis}{\displaystyle}

\newcommand{\noi}{\noindent}

\newcommand{\halmos}{\rule{1ex}{1.4ex}}
\def \qed {\nopagebreak{\hspace*{\fill}$\halmos$\medskip}}
\newcommand{\med}{\medskip}

\newtheorem{theorem}{Theorem}[section]
\newtheorem{proposition}[theorem]{Proposition}
\newtheorem{corollary}[theorem]{Corollary}
\newtheorem{conjecture}[theorem]{Conjecture}
\newtheorem{exercise}[theorem]{Exercise}
\newtheorem{lemma}[theorem]{Lemma}

\newtheorem{remark}[theorem]{Remark}
\newtheorem{defi}[theorem]{Definition}

\newcommand{\bd}{\begin{defi}}
\newcommand{\ed}{\end{defi}}
\newcommand{\bex}{\begin{exercise}}
\newcommand{\eex}{\end{exercise}}
\newcommand{\bt}{\begin{theorem}}
\newcommand{\et}{\end{theorem}}
\newcommand{\bl}{\begin{lemma}}
\newcommand{\el}{\end{lemma}}
\newcommand{\bp}{\begin{proposition}}
\newcommand{\ep}{\end{proposition}}
\newcommand{\bcor}{\begin{corollary}}
\newcommand{\ecor}{\end{corollary}}
\newcommand{\br}{\begin{remark}\rm}
\newcommand{\er}{\end{remark}}
\newcommand{\bcon}{\begin{conjecture}}
\newcommand{\econ}{\end{conjecture}}

\renewcommand{\theequation}{\thesection .\arabic{equation}}

\newcommand{\be}{\begin{equation}}
\newcommand{\ee}{\end{equation}}

\newcommand{\ben}{\begin{equation*}}
\newcommand{\een}{\end{equation*}}
\newcommand{\bea}{\begin{eqnarray}}
\newcommand{\eea}{\end{eqnarray}}
\newcommand{\bean}{\begin{eqnarray*}}
\newcommand{\eean}{\end{eqnarray*}}

\newcommand{\ba}{\begin{array}}
\newcommand{\ea}{\end{array}}
\newcommand{\bc}{\be\begin{array}{r@{\,}c@{\,}l}}
\newcommand{\ec}{\end{array}\ee}

\newcommand{\bet}{\beta}

\newcommand{\de}{\delta}
\newcommand{\De}{\Delta}
\newcommand{\eps}{\varepsilon}

\newcommand{\sig}{\sigma}

\newcommand{\si}{\ensuremath{\sigma}}

\newcommand{\Bi}{{\cal B}}
\newcommand{\Ci}{{\cal C}}
\newcommand{\Di}{{\cal D}}
\newcommand{\Ei}{{\cal E}}
\newcommand{\Fi}{{\cal F}}

\newcommand{\Hi}{{\cal H}}

\newcommand{\Ki}{{\cal K}}
\newcommand{\Li}{{\cal L}}
\newcommand{\cL}{{\cal L}}
\newcommand{\Mi}{{\cal M}}
\newcommand{\Ni}{{\cal N}}

\newcommand{\Ti}{{\cal T}}

\newcommand{\Wi}{{\cal W}}

\newcommand{\R}{{\mathbb R}}
\newcommand{\N}{{\mathbb N}}
\newcommand{\Z}{{\mathbb Z}}

\newcommand{\Q}{{\mathbb Q}}

\renewcommand{\P}{{\mathbb P}}
\newcommand{\E}{{\mathbb E}}

\newcommand{\li}{\langle}
\newcommand{\re}{\rangle}

\newcommand{\up}{\uparrow}
\newcommand{\down}{\downarrow}
\newcommand{\sub}{\subset}

\newcommand{\Asto}[1]{\underset{{#1}\to\infty}{\Longrightarrow}}
\newcommand{\aston}[1]{\underset{{#1}\to 0}{\longrightarrow}}
\newcommand{\Aston}[1]{\underset{{#1}\to 0}{\Longrightarrow}}

\newcommand{\ov}{\overline}

\newcommand{\ffrac}[2]{{\textstyle\frac{{#1}}{{#2}}}}

\newcommand{\di}{\mathrm{d}}

\newcommand{\start}{\sigma}
\newcommand{\final}{\hat\sigma}

\begin{document}

\makeatletter\@addtoreset{equation}{section}
\makeatother\def\theequation{\thesection.\arabic{equation}}

\renewcommand{\labelenumi}{{(\roman{enumi})}}

\newcommand{\Rc}{R^2_{\rm c}}
\newcommand{\Wl}{\Wi^{\rm l}}
\newcommand{\Wr}{\Wi^{\rm r}}
\newcommand{\El}{\Ei^{\rm l}}
\newcommand{\Er}{\Ei^{\rm r}}
\newcommand{\lp}{l}
\newcommand{\rp}{r}
\newcommand{\pp}{p}
\newcommand{\ul}{{\rm l}}
\newcommand{\ur}{{\rm r}}

\title{\vspace*{-1cm}The Brownian web, the Brownian net, and their universality}
\author{Emmanuel Schertzer$^{\,1}$ \and Rongfeng~Sun$^{\,2}$ \and Jan~M.~Swart$^{\,3}$}


\maketitle

\footnotetext[1]{UPMC University Paris 6, Laboratoire de Probabilit\'es et Mod\`eles Al\'eatoires, CNRS UMR 7599, Paris, France. Email: emmanuel.schertzer@upmc.fr}

\footnotetext[2]{Department of Mathematics, National University of Singapore, 10 Lower Kent Ridge Road, 119076 Singapore. Email:
matsr@nus.edu.sg}

\footnotetext[3]{Institute of Information Theory and Automation of the ASCR (\' UTIA),  Pod vod\'arenskou v\v e\v z\' i 4,
18208 Praha 8, Czech Republic.  Email: swart@utia.cas.cz}

\begin{abstract}\noi

The Brownian web is a collection of one-dimensional coalescing Brownian motions starting from everywhere
in space and time, and the Brownian net is a generalization that also allows branching. They appear in the
diffusive scaling limits of many one-dimensional interacting particle systems with branching and coalescence.
This article gives an introduction to the Brownian web and net, and how they arise in the scaling limits of various
one-dimensional models, focusing mainly on coalescing random walks and random walks in i.i.d.\ space-time random
environments. We will also briefly survey models and results connected to the Brownian web and net, including
alternative topologies, population genetic models, true self-repelling motion, planar aggregation, drainage networks,
oriented percolation, black noise and critical percolation. Some open questions are discussed at the end.

\end{abstract}

\vspace{.3cm}

\noi
\vspace{12pt}

\setcounter{tocdepth}{2}
{\setlength{\parskip}{-2pt}\tableofcontents}

\section{Introduction}

The Brownian web originated from the work of Arratia's Ph.D.\ thesis~\cite{A79}, where he studied diffusive scaling limits of coalescing random walk paths starting from everywhere on $\Z$, which can be seen as the spatial genealogies of the population in the dual voter model on $\Z$. Arratia showed that the collection of coalescing random walks converge to a collection of coalescing Brownian motions on $\R$, starting from every point on $\R$ at time $0$. Subsequently, Arratia~\cite{A81} attempted to generalize his result by
constructing a system of coalescing Brownian motions starting from everywhere in the space-time plane $\R^2$, which would be the
scaling limit of coalescing random walk paths starting from everywhere on $\Z$ at every time $t\in\R$. However, the manuscript~\cite{A81} was never completed, even though fundamental ideas have been laid down. This topic
remained dormant until T\'oth and Werner~\cite{TW98} discovered a surprising
connection between the one-dimensional space-time coalescing Brownian motions
Arratia tried to construct, and an unusual process called the {\em true self-repelling motion},
which is repelled by its own local time profile. Building on ideas from~\cite{A81},
T\'oth and Werner~\cite{TW98} gave a construction of the system of space-time
coalescing Brownian motions, and then used it to construct the true
self-repelling motion.

On the other hand, Fontes, Isopi, Newman and Stein~\cite{FINS01} discovered that this system of space-time coalescing Brownian
motions also arises in the study of aging and scaling limits of one-dimensional
spin systems. To establish weak convergence of discrete models to the system of
coalescing Brownian motions, Fontes et al~\cite{FINR02, FINR04} introduced
a topology such that the system of coalescing Brownian motions starting
from every space-time point can be realized as a random variable taking values in
a Polish space, and they named this random variable {\em the Brownian web}.
An extension to the Brownian web was later introduced by the authors in~\cite{SS08},
and independently by Newman, Ravishankar and Schertzer in~\cite{NRS10}. This object was
named {\em the Brownian net} in~\cite{SS08}, where the coalescing paths in the Brownian
web are also allowed to branch. To counter the effect of instantaneous coalescence,
the branching occurs at an effectively  ``infinite'' rate.

The Brownian web and net have very interesting properties. Their construction
is non-trivial due to the uncountable number of starting points in space-time. Coalescence
allows one to reduce the system to a countable number of starting points. In fact, the
collection of coalescing paths starting from every point on $\R$ at time $0$ immediately
becomes locally finite when time becomes positive, similar to the phenomenon of {\em coming
down from infinity} in Kingman's coalescent (see e.g.~\cite{B09}). In fact, the Brownian web
can be regarded as the spatial analogue of Kingman's coalescent, with the former arising
as the limit of genealogies of the voter model on $\Z$, and the latter arising as the limit of
genealogies of the voter model on the complete graph. The key tool in the analysis of the Brownian
web, as well as the Brownian net, is its self-duality, similar to the self-duality
of critical bond percolation on $\Z^2$. Duality allows one to show that there exist
random space-time points where multiple paths originate, and one can give a complete
classification of these points. The Brownian web and net also admit a coupling,
where the web can be constructed by sampling paths in the net, and conversely,
the net can be constructed from the web by Poisson marking a set of ``pivotal''
points in the web and turning these into points where paths can branch. The latter
construction is similar to the construction of scaling limits of near-critical planar
percolation from that of critical percolation~\cite{CFN06, GPS13a, GPS13b}.

The Brownian web and net give rise to a new universality class. In particular, they are expected to arise as the universal
scaling limits of one-dimensional interacting particle systems with coalescence, resp.\ branching-coalescence.
One such class of models are population genetic models with resampling and selection, whose spatial genealogies undergo branching
and coalescence. Establishing weak convergence to the Brownian web or net can also help in the study of the discrete particle systems themselves. Related models which have been shown to converge to the Brownian web include coalescing
random walks~\cite{NRS05}, succession lines in Poisson trees~\cite{FFW05, CV14, FVV14} and
drainage network type models~\cite{CDF09, CV11, RSS13}. Interesting connections with the Brownian web and net
have also emerged from many unexpected sources, including supercritical oriented percolation on $\Z^{1+1}$~\cite{AS11},
planar aggregation models~\cite{NT12, NT15}, true self-avoiding random walks on $\Z$~\cite{T95, TW98}, random matrix
theory~\cite{TZ11, TYZ12}, and also one-dimensional random walks in i.i.d.\ space-time random environments~\cite{SSS14}. There are also close parallels between the Brownian web and the scaling limit of critical planar percolation, which are the only known examples of two-dimensional {\em black noise}~\cite{T04a, T04b, SS11, EF12}.

The goal of this article is to give an introduction to the Brownian web and net, their basic
properties, and how they arise in the scaling limits of one-dimensional interacting particle
systems with branching and coalescence. We will focus on the key ideas, while referring many
details to the literature. Our emphasis is naturally biased toward our own research.
However, we will also briefly survey related work, including the many interesting connections
mentioned above. We have left out many other closely related studies, including
diffusion-limited reactions~\cite{DbA88, bABD90} where a dynamic phase transition is observed for branching-coalescing random walks,
the propagation of cracks in a sheet~\cite{D15}, rill erosion~\cite{DW09} and directed Abelian Sandpile Model~\cite{D06}, quantum spin chains~\cite{KPWH95}, etc,  which all lie within the general framework of non-equilibrium critical phenomena discussed in the physics surveys~\cite{P97, H00}.

The rest of this article is organized as follows. In
Section~\ref{S:web}, we will construct and give a characterization of the Brownian web and study
its properties. In Section~\ref{S:net}, we do the same for the Brownian net. In Section~\ref{S:couple},
we introduce a coupling between the Brownian web and net and show how one can be constructed from the other.
In Section~\ref{S:rwre}, we will explain how the Brownian web and net can be used to construct the scaling limits of one-dimensional random walks in
i.i.d.\ space-time random environments. In Section~\ref{S:conv}, we formulate
convergence criteria for the Brownian web, which are then applied to coalescing
random walks. We will also discuss strategies for proving convergence to the Brownian net.
In Section~\ref{S:disc}, we survey other interesting models and results connected to the Brownian web and net.
Lastly, in Section~\ref{S:open}, we conclude with some interesting open questions.

\section{The Brownian web}\label{S:web}

\begin{figure}[tp] 
\begin{center}
\includegraphics[width=15cm]{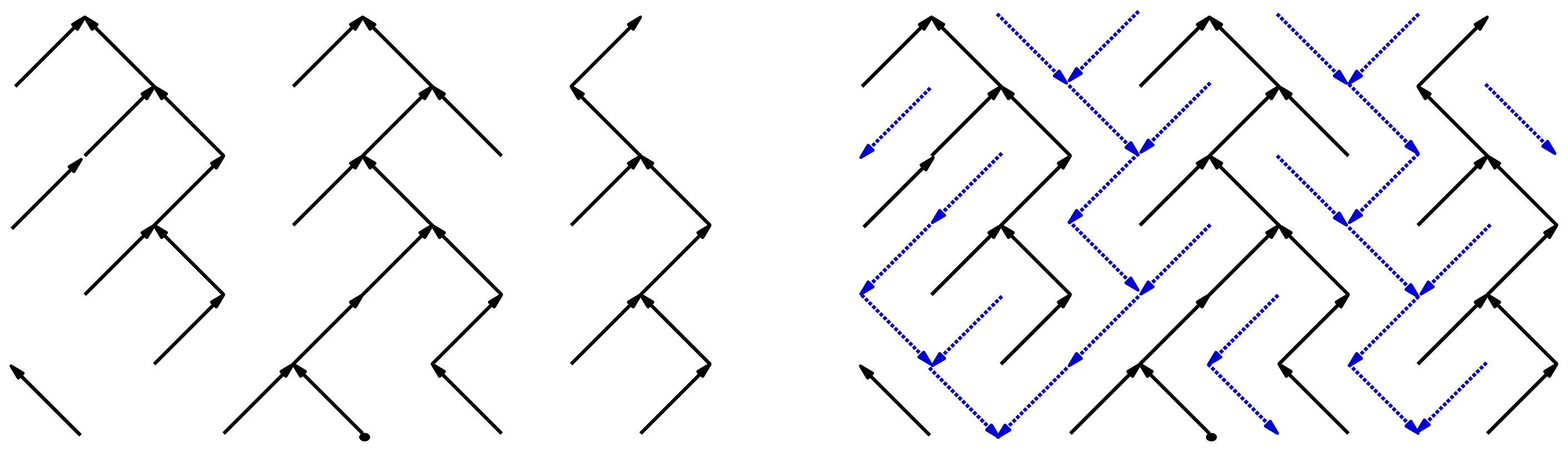}
\caption{Discrete space-time coalescing random walks on $\Z^2_{\rm even}$, and its dual on $\Z^2_{\rm odd}$.}
\label{fig:webdual}
\end{center}
\end{figure}

The Brownian web is best motivated by its discrete analogue, the collection of discrete time coalescing
simple symmetric random walks on $\Z$, with one walker starting from every site in the space-time lattice
$\Z^2_{\rm even}:=\{(x,n)\in \Z^2: x+n \mbox{ is even}\}$. The restriction to the sublattice $\Z^2_{\rm even}$
is necessary due to parity. Figure~\ref{fig:webdual} illustrates a graphical construction, where from each
$(x,n)\in \Z^2_{\rm even}$, an independent arrow is drawn from $(x,n)$ to either $(x-1, n+1)$ or $(x+1, n+1)$
with probability $1/2$ each, determining whether the walk starting at $x$ at time $n$ should move to $x-1$ or
$x+1$ at time $n+1$. The objects of interest for us are the collection of upward random walk paths (obtained by
following the arrows) starting from every space-time lattice point. The question is:
\begin{itemize}
\item[\bf Q.1] What is the diffusive scaling
limit of this collection of coalescing random walk paths, if space and time are scaled by $1/\sqrt{n}$ and $1/n$
respectively?
\end{itemize}
Intuitively, it is not difficult to see that the limit should be a collection of coalescing Brownian motions,
starting from everywhere in the space-time plane $\R^2$. This is what we will call the {\em Brownian web}. However,
a conceptual difficulty arises, namely that we need to construct the joint realization of {\em uncountably} many Brownian motions.
Fortunately it turns out that coalescence allows us to reduce the construction to only a countable collection of coalescing Brownian motions.

Note that in Figure~\ref{fig:webdual}, we have also drawn a collection of downward arrows connecting points in the odd space-time
lattice $\Z^2_{\rm odd}:= \{(x,n)\in \Z^2: x+n \mbox{ is odd}\}$, which are dual to the upward arrows by the constraint
that the upward and backward arrows do not cross each other. This is the same duality as that for planar bond percolation,
and the collection of upward arrows uniquely determine the downward arrows, and vice versa. The collection of downward
arrows determine a collection of coalescing random walk paths running backward in time, with one walker starting from each
site in $\Z^2_{\rm odd}$. We may thus strengthen {\bf Q.1} to the following:
\begin{itemize}
\item[\bf Q.2] What is the diffusive scaling limit of the joint realization of the collection of forward and
backward coalescing random walk paths?
\end{itemize}
Observe that the collection of backward coalescing random walk paths has the same distribution as the forward collection,
except for a rotation in space-time by $180^o$ and a lattice shift. Therefore, the natural answer to {\bf Q.2} is that the limit consists of
two collections of coalescing Brownian motions starting from everywhere in space-time, one running forward in time and the
other backward, and the two collections are equally distributed except for a time-reversal. This is what we will call the
(forward) {\em Brownian web} and the {\em dual (backward) Brownian web}.

In the discrete system, we observe that the collection of
forward and the collection of backward coalescing random walk paths uniquely determine each other by the constraint that forward and backward
paths cannot cross. It is natural to expect the same for their continuum limits, namely that the Brownian web and
the dual Brownian web almost surely uniquely determine each other by the constraint that their paths cannot cross.

The heuristic considerations above, based on discrete approximations, outline the key properties that we expect the Brownian web to
satisfy and provide a guide for our analysis.

Before proceeding to a proper construction of the Brownian web and
establishing its basic properties, we first define a suitable
Polish space in which the Brownian web takes its value. This will be essential to prove weak convergence to the Brownian web.

\subsection{The space of compact sets of paths} \label{SS:hausdorff}

Following Fontes et al~\cite{FINR04}, we regard the collection of colaescing Brownian motions as a {\em set} of space-time paths, which can be shown to be almost surely relatively compact if space and time are suitably compactified. It is known that given a Polish space $E$ (the space of paths in our case), the space of compact subsets of $E$, equipped with the induced Hausdorff topology, is a Polish space itself. Therefore a natural space for the Brownian web is the space of {\em compact sets of paths} (after compactifying space and time), with the Brownian web taken to be the almost sure closure of the set of colaescing Brownian motions. This {\em paths topology} was inspired by a similar topology proposed by Aizenman~\cite{A98} to study two-dimensional percolation configurations as a closed sets of curves, called the {\em percolation web}, which was then studied rigorously by Aizenman and Burchard in~\cite{AB99}. We now give the details.

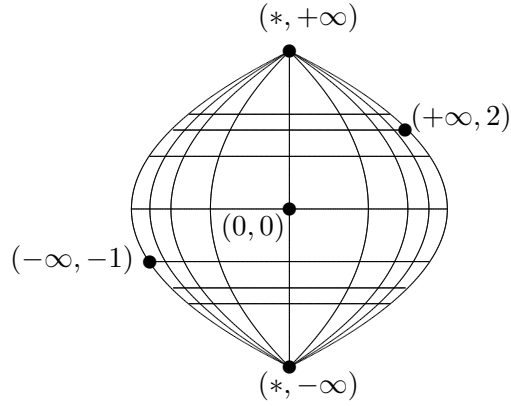
\begin{figure}
\begin{center}
\setlength{\unitlength}{.7cm}
\begin{picture}(10,8)(-5,-4)
\linethickness{.4pt}
\qbezier(0,-3)(0,0)(0,3)
\qbezier(-3,0)(0,0)(3,0)
\linethickness{.4pt}
\qbezier(0,-3)(-6,0)(0,3)
\qbezier(0,-3)(-5.3,0)(0,3)
\qbezier(0,-3)(-4.5,0)(0,3)
\qbezier(0,-3)(-3,0)(0,3)
\qbezier(0,-3)(3,0)(0,3)
\qbezier(0,-3)(4.5,0)(0,3)
\qbezier(0,-3)(5.3,0)(0,3)
\qbezier(0,-3)(6,0)(0,3)
\qbezier(-1.9,-1.8)(0,-1.8)(1.9,-1.8)
\qbezier(-2.2,-1.5)(0,-1.5)(2.2,-1.5)
\qbezier(-2.65,-1)(0,-1)(2.65,-1)
\qbezier(-2.65,1)(0,1)(2.65,1)
\qbezier(-2.2,1.5)(0,1.5)(2.2,1.5)
\qbezier(-1.9,1.8)(0,1.8)(1.9,1.8)

\put(0,-3){\circle*{.25}}
\put(0,3){\circle*{.25}}
\put(0,0){\circle*{.25}}
\put(2.2,1.5){\circle*{.25}}
\put(-2.65,-1){\circle*{.25}}

\put(-.6,-3.5){$(\ast,-\infty)$}
\put(-.6,3.5){$(\ast,+\infty)$}
\put(-1.3,-.5){$(0,0)$}
\put(2.3,1.6){$(+\infty,2)$}
\put(-5.3,-1.1){$(-\infty,-1)$}

\end{picture}
\caption[The compactification $\Rc$ of $\R^2$.]{The compactification
$\Rc$ of $\R^2$.}\label{fig:comp}
\end{center}
\end{figure}

We first compactify $\R^2$. Let $\Rc$ denote the completion of the space-time plane $\R^2$ w.r.t.\ the metric
\be\label{rho}
\rho\big((x_1, t_1), (x_2,t_2)\big) = \left|\tanh(t_1)-\tanh(t_2)\right|
\ \vee\ \left|\frac{\tanh(x_1)}{1+|t_1|}-\frac{\tanh(x_2)}{1+|t_2|}\right|.
\ee
Note that $\Rc$ can be identified with the continuous image of
$[-\infty, \infty]^2$ under a map that identifies the line
$[-\infty,\infty]\times\{\infty\}$ with a single point $(*, \infty)$, and the
line $[-\infty,\infty]\times\{-\infty\}$ with the point $(*,-\infty)$, see
Figure~\ref{fig:comp}.

A path $\pi$ in $\Rc$, whose starting time we denote by
$\sigma_\pi\in [-\infty,\infty]$, is a mapping $\pi :
[\sigma_\pi,\infty] \to [-\infty, \infty] \cup\{*\}$ such that
$\pi(\infty)=*$, $\pi(\sigma_\pi)=*$ if $\sigma_\pi=-\infty$, and $t
\to (\pi(t), t)$ is a continuous map from $[\sigma_\pi,\infty]$ to
$(\Rc, \rho)$. We then define $\Pi$ to be the space of all
paths in $\Rc$ with all possible starting times in $[-\infty,\infty]$.
Endowed with the metric
\be\label{PId}
\!
d(\pi_1, \pi_2) = \Big|\!\tanh(\sigma_{\pi_1})-\tanh(\sigma_{\pi_2})\Big|
\ \vee \sup_{t\geq \sigma_{\pi_1} \wedge \sigma_{\pi_2}}\!
\left|\frac{\tanh(\pi_1(t\vee \sigma_{\pi_1}))}{1+|t|}
-\frac{\tanh(\pi_2(t\vee \sigma_{\pi_2}))}{1+|t|}\right|, \!\!\!
\ee
$(\Pi, d)$ is a complete separable metric space. Note that
convergence in the metric $d$ can be desrcibed as locally uniform convergence
of paths plus convergence of starting times. (The metric
$d$ differs slightly from the original choice in \cite{FINR04}, which
is somewhat less natural as explained in the appendix of \cite{SS08}.)

Let $\Hi$ denote the {\em space of compact subsets of $(\Pi, d)$},
equipped with the Hausdorff metric
\be\label{dH}
d_{\Hi}(K_1, K_2) = \sup_{\pi_1\in K_1} \inf_{\pi_2\in K_2}\!\!
d(\pi_1, \pi_2)\ \vee  \sup_{\pi_2\in K_2}
\inf_{\pi_1\in K_1} d(\pi_1, \pi_2),
\ee
and let $\Bi_\Hi$ be the Borel $\sigma$-algebra associated with $d_\Hi$.
\begin{exercise}
Show that $(\Hi, d_\Hi)$ is a complete separable metric space.
\end{exercise}

\bex\label{E:comsub}
Let $K\subset \Pi$ be compact. Show that $\Ki:=\{ \bar A: A\subset K\}$ is
a compact subset of $\Hi$.
\eex
For further properties of $(\Hi, d_\Hi)$, such as textcolor{red}{a} criterion for the convergence of a
sequence of elements in $\Hi$, or necessary and sufficient conditions for
the precompactness of a subset of $\Hi$, see e.g.~\cite[Appendix B]{SSS14}.

We will construct the Brownian web $\Wi$ as an $(\Hi, \Bi_\Hi)$-valued random variable.
The following notational convention will be adopted in the rest of this article:
\begin{itemize}
\item For $K\in\Hi$ and $A\subset\Rc$, $K(A)$ will denote the set of paths in $K$ with starting points in $A$.

\item When $A=\{z\}$ for $z\in\Rc$, we also write $K(z)$ instead of $K(\{z\})$.
\end{itemize}

\subsection{Construction and characterization of the Brownian web}

The basic ideas in constructing the Brownian web are the following. First we can construct coalescing Brownian motions starting from a deterministic countable dense subset $\Di$ of the space-time plane $\R^2$. It is easily seen that coalescence forces paths started at typical points outside $\Di$ to be squeezed between coalescing paths started from $\Di$. Therefore to construct paths starting from outside $\Di$, we only need to take the closure of the set of paths starting from $\Di$. Lastly one shows that the law of the random set of paths constructed does not depend on the choice of $\Di$. This construction procedure is effectively contained in the following result from~\cite[Theorem 2.1]{FINR04}, which gives a
characterization of the Brownian web $\Wi$ as an $(\Hi, \Bi_\Hi)$-valued random variable, i.e., a random compact set of paths.

\bt[Characterization of the Brownian web]\label{T:webchar}
There exists an $(\Hi, \Bi_\Hi)$-valued random variable $\Wi$, called
the standard Brownian web, whose distribution is uniquely determined
by the following properties:
\begin{itemize}
\item[{\rm (a)}] For each deterministic $z\in\R^2$, almost surely
there is a unique path $\pi_z\in \Wi(z)$.

\item[{\rm (b)}] For any finite deterministic set of points $z_1, \ldots, z_k
  \in\R^2$, the collection $(\pi_{z_1}, \ldots, \pi_{z_k})$
  is distributed as coalescing Brownian motions.

\item[{\rm (c)}] For any deterministic countable dense subset
$\Di \subset\R^2$, almost surely, $\Wi$ is the closure of
$\{\pi_z : z\in\Di\}$ in $(\Pi, d)$.
\end{itemize}
\et
{\bf Proof Sketch.} We will sketch the main ideas and ingredients and refer the details to~\cite{FINR03, FINR04}. The main steps are:
\begin{itemize}
\item[\bf (1)] Let $\Di=\{(x,t): x,t \in \Q\}$ and construct the collection of coalescing Brownian motions
$\Wi(\Di):=\{\pi_z\}_{z\in\Di}$, where $\pi_z$ is the Brownian motion starting at $z$.

\item[\bf (2)] Show that $\Wi(\Di)$ is almost surely a pre-compact set in the space of paths $(\Pi, d)$, and hence $\Wi:=\overline{\Wi(\Di)}$ defines a random compact set, i.e., an $(\Hi, \Bi_\Hi)$-valued random variable.

\item[\bf (3)] Show that properties (a) and (b) hold for $\Wi$, which can be easily seen to imply that property (c) also holds for $\Wi$.
\end{itemize}
The above steps construct a random variable $\Wi$ satisfying properties (a)--(c). Its law is uniquely determined, since if $\widetilde\Wi$ is another random variable satisfying the same properties, then both $\Wi(\Di)$ and $\widetilde \Wi(\Di)$ are coalescing Brownian motions starting from $\Di$, and hence can be coupled to equal almost surely. Property (c) then implies that $\Wi=\widetilde \Wi$ almost surely under this coupling.

{\bf Step (1).} Fix an order for points in $\Di$, so that $\Di=\{z_k\}_{k\in\N}$. Coalescing Brownian motions $(\pi_k)_{k\in\N}$ starting respectively from $(z_k)_{k\in\N}$ can be constructed inductively from independent Brownian motions $(\tilde \pi_k)_{k\in\N}$ starting from $(z_k)_{k\in\N}$. First let $\pi_1 :=\tilde\pi_1$. Assuming that $\pi_1,\ldots, \pi_k$ have already been constructed from $\tilde \pi_1, \ldots, \tilde \pi_k$, then we define the path $\pi_{k+1}$ to coincide with the independent Brownian motion $\tilde \pi_{k+1}$ until the first
time $\tau$ when it meets one of the already constructed coalescing paths, say $\pi_j$, for some $1\leq j\leq k$. From time $\tau$ onward, we just set $\pi_{k+1}$ to coincide with $\pi_j$. It is not difficult to see that for any $k\in\N$, $(\pi_i)_{1\leq i\leq k}$ is a collection of coalescing Brownian motions characterized by the property that, different paths evolve as independent Brownian motions when they are apart, and evolve as the same Brownian motion from the time when they first meet. Furthermore, any subset of a collection of coalescing Brownian motions is also a collection of coalescing Brownian motions.

{\bf Step (2).} The main idea is the following. The compactification of space-time as shown in Figure~\ref{fig:comp} allows us to approximate $\R^2$ by a large space-time box $\Lambda_{L,T}:=[-L, L]\times [-T, T]$, and proving precompactness of $\Wi(\Di)$ can be reduced to proving the equi-continuity of paths in $\Wi(\Di)$ restricted to $\Lambda_{L,T}$ (for further details, see~\cite[Appendix B]{FINR04}). More precisely, it suffices
to show that for any $\eps>0$, almost surely we can choose $\delta>0$ such that the modulus of continuity
\begin{equation}\label{psiwi}
\psi_{\Wi(\Di), L, T}(\delta) := \sup\{ \, |\pi_z(t)-\pi_z(s)|\, :\,  z\in \Di,  (\pi_z(s), s)\in \Lambda_{L,T},t\in [s, s+\delta]\}\ \leq\ \eps.
\end{equation}
Assuming w.l.o.g.\ that $\eps, \delta \in \Q$, we will control $\psi_{\Wi(\Di), L, T}(\delta)>\eps$ in terms of the modulus of continuity of coalescing Brownian motions starting from the grid
$$
G_{\eps, \delta}:=\{(m\eps/4, n\delta): m,n\in\N\}\cap [-L-\eps, L+\eps]\times [-T-\delta, T] \ \subset \Di.
$$
Indeed, $\psi_{\Wi(\Di), L, T}(\delta)>\eps$ means that $|\pi_z(t)-\pi_z(s)|>\eps$ for some $z\in \Di$ with $(\pi_z(s), s)\in \Lambda_{L,T}$ and
$t\in [s, s+\delta]$. Then there exists a point in the grid $\tilde z=(\tilde x, \tilde t)\in G_{\eps, \delta}$ with $s\in [\tilde t, \tilde t+\delta)$ and $\tilde x\in (\pi(s)\wedge \pi(t)+ \eps/4, \pi(s)\wedge \pi(t)-\eps/4)$. Since $\pi_z$ and $\pi_{\tilde z}$ are coalescing Brownian
motions, either $\pi_{\tilde z}$ coalesces with $\pi_z$ before time $t$, or $\pi_{\tilde z}$ avoids $\pi_z$ up to time $t$. Either way, we must
have
$$
\sup_{h\in [0, 2\delta]} |\pi_{\tilde z}(\tilde t+h) - \pi_{\tilde z}(\tilde t)|\geq \eps/4.
$$
Denote this event by $E^{\eps, \delta}_z$. Then
$$
\begin{aligned}
\P(\psi_{\Wi(\Di), L, T}(\delta) >\eps) \leq \P\Big(\bigcup_{z\in G_{\eps, \delta}} E^{\eps, \delta}_z \Big) & \leq \sum_{z\in G_{\eps, \delta}} \P(E^{\eps, \delta}_z)  = |G_{\eps,\delta}| \P(\sup_{h\in [0,2\delta]} |B_h|\geq \eps/4) \\
& \leq C_{L,T} \eps^{-1}\delta^{-1} e^{-c\eps^2/\delta},
\end{aligned}
$$
where $c, C_{L,T}>0$, $B$ is a standard Brownian motion, and we have used the reflection principle to bound the tail probability for $\sup |B|$.
Since $\P(\psi_{\Wi(\Di), L, T}(\delta) >\eps)\to 0$ as $\delta\downarrow 0$, this implies \eqref{psiwi}. Therefore $\Wi(\Di)$ is a.s.\ precompact, and $\Wi:=\overline{\Wi(\Di)}$ defines an $(\Hi, \Bi_\Hi)$-valued random variable.

{\bf Step (3).} We first show that for each $z=(x,t)\in\R^2$, almost surely $\Wi(z)$, the paths in $\Wi$ starting at $z$, contains a unique path.
Let $z_n^-=(x-\eps_n, t-\delta_n)\in \Di$, $z_n^+=(x+\eps_n, t-\delta_n)\in\Di$, with $\eps_n, \delta_n\downarrow 0$, and let $\tau_n$ be the time when $\pi_{z_n^-}$ and $\pi_{z_n^+}$ coalesce. Note that on the event
$$
E_n:= \Big\{\pi_{z_n^-}(t) < x < \pi_{z_n^+}(t), \tau_n \leq t+\frac{1}{n} \Big\},
$$
every path in $\Wi(z)$ must be enclosed between $\pi_{z_n^-}$ and $\pi_{z_n^+}$, and hence is uniquely determined from time
$\tau_n\leq t+\frac{1}{n}$ onward. It is easy to see that we can choose $\eps_n\downarrow 0$ sufficiently fast, and $\delta_n\downarrow 0$
much faster than $\eps_n^2$, such that $\P(E_n)\to 1$ as $n\to\infty$. In particular, almost surely, $E_n$ occurs infinitely often,
which implies that the paths in $\Wi(z)$ all coincide on $(t,\infty)$ and hence $\Wi(z)$ contains a unique path.

To show that $\Wi$ satisfies property (b), let us fix $z_1, \ldots, z_k\in\R^2$. For each $1\leq i\leq k$, let $z_{n,i}\in \Di$ with $z_{n,i}\to z_i$ as $n\to\infty$. By the a.s.\ compactness of $\Wi$, and the fact that $\Wi(z_i)$ a.s.\ contains a unique path $\pi_{z_i}$ by property (a) that we just verified, we must have $\pi_{z_{n,i}}\to \pi_{z_i}$ in $(\Pi, d)$ for each $1\leq i\leq k$. In particular, as a sequence of $\Pi^k$-valued random variables, $(\pi_{z_{n,i}})_{1\leq i\leq k}$ converges in distribution to $(\pi_{z_{i}})_{1\leq i\leq k}$. On the other hand, as a subset of $\Wi(\Di)$, $(\pi_{z_{n,i}})_{1\leq i\leq k}$ is a collection of coalescing Brownian motions, and it is easy to show that as their starting points converge, they converge in distribution to a collection of coalescing Brownian motions starting from $(z_i)_{1\leq i\leq k}$. Therefore $(\pi_{z_{i}})_{1\leq i\leq k}$ is distributed as a collection of coalescing Brownian motions.

Lastly to show that $\Wi$ satisfies property (c), let $\Di'$ be another countable dense subset of $\R^2$. Clearly $\overline{\Wi(D')}\subset \Wi$. To show the converse, $\Wi\subset \overline{\Wi(D')}$, it suffices to show that for each $z\in \Di$, $\pi_z\in \overline{\Wi(D')}$. This can
be seen by taking a sequence $z'_n\in \Di'$ with $z'_n\to z$, for which we must have $\pi_{z_n'}\to \pi_z\in \overline{\Wi(D')}$ by the compactness of $\overline{\Wi(D')}\subset \Wi$ and the fact that $\Wi(z)=\{\pi_z\}$.
\qed

\subsection{The Brownian web and its dual}\label{SS:webdual}

As discussed at the beginning of Section~\ref{S:web}, similar to the duality between forward and backward coalescing random walks shown in Figure~\ref{fig:webdual}, we expect the Brownian web $\Wi$ also to have a dual $\widehat \Wi$. Such a duality provides a powerful tool for analyzing properties of the Brownian web. Since the dual Brownian web $\widehat \Wi$ should be a collection of coalescing paths running backward in time, we first define the space in which $\widehat \Wi$ takes its values.

Given $z=(x,t)\in R_c^2$, which is identified with $[-\infty, \infty]^2$ where $[-\infty, \infty]\times \{\pm\infty\}$ is contracted to a single point $(*, \pm\infty)$, let $-z$ denote $(-x, -t)$. Given a set $A\subset R_c^2$, let $-A$ denote $\{-z: z\in A\}$. Identifying
each path $\pi \in \Pi$ with its graph as a subset of $R_c^2$, $\hat \pi :=-\pi$ defines a path running backward in time, with starting time $\hat \sigma_{\hat \pi}=-\sigma_\pi$. Let $\widehat \Pi:=-\Pi$ denote the set of all such backward paths, equipped with a metric $\hat d$ that is inherited from $(\Pi, d)$ under the mapping $-$. Let $\widehat \Hi$ be the space of compact subsets of $(\widehat \Pi, \hat d)$, equipped with the Hausdorff metric $d_{\widehat \Hi}$ and Borel $\sigma$-algebra $\Bi_{\widehat \Hi}$. For any $K\in \Hi$, we will let $-K$ denote the set $\{-\pi: \pi \in K\}\in \widehat \Hi$.

The following result characterizes the joint law of the Brownian web $\Wi$ and its dual $\widehat \Wi$ as a random variable taking values in
$\Hi\times \widehat\Hi$, equipped with the product $\sigma$-algebra.

\bt[Characterization of the double Brownian web]\label{T:dwebchar}
There exists an $\Hi\times \widehat \Hi$-valued random variable $(\Wi, \widehat \Wi)$, called the
double Brownian web $($with $\widehat \Wi$ called the dual Brownian web$)$, whose distribution is uniquely determined by the following properties:
\begin{itemize}
\item[{\rm (a)}] $\Wi$ and $-\widehat \Wi$ are both distributed as the standard Brownian web.

\item[{\rm (b)}] Almost surely, no path $\pi_z\in \Wi$ crosses any path $\hat \pi_{\hat z}\in \widehat \Wi$ in the sense that, $z=(x,t)$
and $\hat z=(\hat x, \hat t)$ with $t<\hat t$, and $(\pi_z(s_1)-\hat \pi_{\hat z}(s_1)) (\pi_z(s_2)-\hat \pi_{\hat z}(s_2))<0$ for some
$t<s_1<s_2<\hat t$.
\end{itemize}
Furthermore, for each $z\in\R^2$, $\widehat \Wi(z)$ a.s.\ consists of a single path $\hat \pi_z$ which is the unique path in $\widehat \Pi$ that does not cross any path in $\Wi$, and thus $\widehat \Wi$ is a.s.\ determined by $\Wi$ and vice versa.
\et
{\bf Proof Sketch.} The existence of a double Brownian web $(\Wi, \widehat \Wi)$ satisfying properties (a)--(b) is most easily derived as scaling limits of forward and backward coalescing random walks. We defer this to Section~\ref{S:conv}, after we introduce general criteria for convergence to the Brownian web.

Let us first prove that if $(\Wi, \widehat \Wi)$ satisfies properties (a)-(b), then almost surely $\Wi$ uniquely determines $\widehat\Wi$. Indeed,
fix a deterministic $z=(x,t)\in\R^2$. By the characterization of the Brownian web $\Wi$, $\Wi(\Q^2)$ is a collection of colaescing Brownian motions, with a.s.\ one Brownian motion starting from each point in $\Q^2$. Since Brownian motion has zero probability of hitting a deterministic space-time point, there is zero probability that $z$ lies on $\pi_{z'}$ for any $z'\in\Q^2$. Therefore for any $s\in \Q$
with $s<t$, property (b) implies that for any path $\hat\pi \in \widehat\Wi(z)$, we must have
$$
\hat \pi(s) = \sup \{ y\in \Q: \pi_{(y,s)}(t)<x\} = \inf\{ y\in \Q: \pi_{(y,s)}(t)>x\}.
$$
In other words, $\hat\pi$ is uniquely determined at rational times and hence at all times, and
$\widehat\Wi(z)$ contains a unique path. It follows that $\widehat \Wi(\Q^2)$ is a.s.\ uniquely
determined by $\Wi$, and hence so is $\widehat\Wi = \overline{\widehat\Wi(\Q^2)}$.

Lastly we show that the distribution of $(\Wi, \widehat\Wi)$ is uniquely determined by properties (a) and (b). Indeed, if $(\Wi', \widehat\Wi')$ is another double Brownian web, then $\Wi$ and $\Wi'$ can be coupled so that they equal a.s. As we have just shown, $\Wi$ a.s.\ uniquely determines
$\widehat\Wi$, and $\Wi'$ determines $\widehat\Wi'$. Therefore $\widehat\Wi=\widehat\Wi'$ a.s., and $(\Wi', \widehat\Wi')$ has the same distribution as $(\Wi, \widehat \Wi)$.
\qed

\br\label{R:skoref} One can characterize the joint law of paths in $(\Wi, \widehat\Wi)$ starting from a finite deterministic set of points. Similar to the construction of coalescing Brownian motions, we can construct one path at a time. To add a new forward path to an existing collection, we follow an independent Brownian motion until it either meets an existing forward Brownian motion, in which case they coalesce, or it meets an existing dual Brownian motion, in which case it is {\em Skorohod reflected} by the dual Brownian motion. For further details, see~\cite{STW00}. Extending this pathwise construction to a countable dense set of starting points $\Di$ and then taking closure, this gives a direct construction of $(\Wi, \widehat\Wi)$, which is formulated in~\cite[Theorem 3.7]{FINR06}.
\er

\br\label{R:webdualchar}
In light of Theorem~\ref{T:dwebchar}, one may wonder whether $\Wi$ a.s.\ consists of {\em all} paths in $\Pi$ which do not cross any path in $\widehat\Wi$, and vice versa. The answer is no, and $\Wi$ is actually the minimal compact set of paths that do not cross any path in $\widehat\Wi$ while still containing paths starting from every point in $\R^2$. More non-crossing paths can be added to $\Wi$ by extending paths in $\Wi$ backward in time, following paths in $\widehat\Wi$ (see~\cite{FN06}). Such paths can be excluded if we impose the further restriction that no path can enter from outside any open region enclosed by a pair of paths in $\widehat\Wi$. This is called the {\em wedge characterization} of the Brownian web, to be discussed in more detail in Remark~\ref{R:dualchar}.
\er

\subsection{The coalescing point set}\label{S:coptset}
The coupling between the Brownian web and its dual given in Theorem~\ref{T:dwebchar} allows one to deduce interesting properties for the Brownian web. The first result is on the density of paths in the Brownian web $\Wi$ started at time $0$.

Given the Brownian web $\Wi$, and a closed set $A\subset \R$, define the {\em coalescing point set} by
\begin{equation}\label{coalptset}
\xi^A_t := \{ y\in \R: y=\pi(t) \mbox{ for some } \pi \in \Wi(A\times \{0\}) \}, \qquad t\geq 0.
\end{equation}
In words, $\xi^A_t$ is the set of points in $\R$ that lie on some path in $\Wi$ that start from $A$ at time $0$. Note that this process is monotone in the sense that if $A\subset B$, then $\xi^A_t\subset \xi^B_t$ a.s.\ for all $t\geq 0$.

It turns out that even if started from the whole line, $\xi^\R_t$ becomes a.s.\ locally finite as soon as $t>0$, as the following density
result shows. Such a {\em coming down from infinity} phenomena also appears in Kingman's coalescent, see e.g.~\cite{B09}.

\bp[Density of the coalescing point set]\label{P:webdensity} Let $\xi^\R_\cdot$ be the coalescing point set defined from the Brownian web
$\Wi$ as in \eqref{coalptset}. Then for all $t>0$ and $a<b$,
\begin{equation}\label{webdensity}
\E[ |\xi^\R_t \cap [a,b]| ] = \frac{b-a}{\sqrt{\pi t}}.
\end{equation}
\ep
{\bf Proof.} Let $\widehat \Wi$ be the dual Brownian web determined a.s.\ by $\Wi$, as in Theorem~\ref{T:dwebchar}. Observe that by the non-crossing property between paths in $\Wi$ and $\widehat \Wi$, $\xi^\R_t \cap (a,b)\neq\emptyset$ implies that the paths $\hat \pi_a, \hat \pi_b\in \widehat\Wi$, starting respectively at $(a,t)$ and $(b,t)$, do not coalesce in the time interval $(0,t)$ (i.e., $\hat\tau\leq 0$ if
$\hat\tau$ denotes the time when $\hat\pi_a$ and $\hat \pi_b$ coalesce). Conversely, if $\hat \tau\leq 0$, then any path in $\Wi$ started
from $[\hat \pi_a(0), \hat\pi_b(0)]$ at time $0$ will hit $[a,b]$ at time $t$, i.e., $\xi^\R_t\cap [a,b]\neq \emptyset$. Thus,
\begin{equation}\label{wd1}
\P(\xi^\R_t \cap (a,b)\neq \emptyset) \leq \P(\hat\tau \leq 0) \leq \P(\xi^\R_t \cap [a,b]\neq \emptyset),
\end{equation}
where we observe that
$$
\begin{aligned}
\P(\hat\tau\leq 0) = \P\Big(\sup_{s\in (0,t)} (B_2(s)-B_1(s)) \leq b-a\Big) & = \P\Big(\sup_{s\in (0,t)} B(s) \leq \frac{b-a}{\sqrt 2}\Big) \\
& = \frac{1}{\sqrt{2\pi t}}\int_{-\frac{b-a}{\sqrt 2}}^{\frac{b-a}{\sqrt 2}} e^{-\frac{x^2}{2t}} {\rm d}x
\sim \frac{b-a}{\sqrt{\pi t}} \quad \mbox{as } b-a\downarrow 0.
\end{aligned}
$$
By \eqref{wd1}, this implies that
$$
\P(x\in \xi^\R_t) =0 \quad \mbox{ for all } x\in \R,
$$
and the inequalities in \eqref{wd1} are in fact all equalities.

We can then apply monotone convergence theorem to obtain
$$
\begin{aligned}
\E[|\xi^\R_t \cap [a,b]| ] & = \lim_{n\to\infty} \E\Big[\,\Big| \Big\{ 1\leq i< (b-a)2^n: \xi^\R_t \cap \big(a+\frac{i-1}{2^n}, a+\frac{i}{2^n}\big)\neq\emptyset\Big\}\Big|\,\Big] \\
& = \lim_{n\to\infty} (b-a) 2^n \P(\xi^\R_t \cap (0, 2^{-n})\neq\emptyset) = \lim_{n\to\infty} (b-a) 2^n \frac{2^{-n}}{\sqrt{\pi t}}
= \frac{b-a}{\sqrt{\pi t}}.
\end{aligned}
$$
This concludes the proof of the proposition.
\qed

As a corollary of Proposition~\ref{P:webdensity}, we show that when paths in the Brownian web converge, they converge in a strong sense (see e.g.~\cite[Lemma 3.4]{SS08}).
\bcor[Strong convergence of paths in $\Wi$]\label{C:pathconv}
Let $\Wi$ be the standard Brownian web. Almost surely, for any sequence $\pi_n\in \Wi$ with $\pi_n\to \pi\in\Wi$, the time of coalescence $\tau_n$ between $\pi_n$ and $\pi$ must tend to $\sigma_\pi$ as $n\to\infty$.
\ecor

\begin{exercise}
Deduce Corollary~\ref{C:pathconv} from Proposition~\ref{P:webdensity}.
\end{exercise}
%

\br\label{R:pfaffian}
Apart from its density, we actually know quite a bit more about the coalescing point set $\xi^\R_t$. It has been shown by Tribe et al \cite{TZ11, TYZ12} that $\xi^\R_t$ is in fact a {\em Pfaffian point process}, whose kernel also appears in the real Ginibre random matrix ensemble. Furthermore, $\xi^\R_t$ (and more generally $\xi^A_t$ for any $A\subset \R$) can be shown (see e.g.~\cite[Appendix C]{GSW15} and \cite{MRTZ06}) to be negatively associated in the sense that for any $n\in\N$ and any disjoint open intervals $O_1,\cdots, O_n$, we have
\begin{equation}
\P(\cap_{i=1}^n \{\xi^\R_t\cap O_i\neq\emptyset\}) \leq \prod_{i=1}^n \P(\xi^\R_t\cap O_i\neq\emptyset).
\end{equation}
For any $B\subset \R$ with positive Lebesgue measure, we also have
\begin{equation}
\P(|\xi^\R_t \cap B|\geq m+n) \leq \P(|\xi^\R_t \cap B|\geq m) \P(|\xi^\R_t \cap B|\geq n) \qquad \forall\, m,n\in\N.
\end{equation}
On a side note, we remark that when $A\subset \R$ is a finite set, determinantal formulas have also been derived for the distribution of $\xi^A_t$ in~\cite[Prop.~9]{W07}.
\er

\subsection{Special points of the Brownian web}\label{S:webpts}

We have seen in Theorem~\ref{T:webchar} that for each deterministic $z\in \R^2$, almost surely the Brownian web $\Wi$ contains a unique path starting from $z$. However, it is easily seen that there must exist random points $z\in\R^2$ where $\Wi(z)$ contains multiple paths. Indeed,
consider paths in $\Wi$ starting from $\R$ at time $0$. Proposition~\ref{P:webdensity} shows that these paths coalesce into a locally
finite set of points $\xi^\R_t$ at any time $t>0$. Each point $x_i\in \xi^\R_t$ (with $x_i<x_{i+1}$ for all $i\in\Z$) can be traced back to an interval $(u_i, u_{i+1})$ at time $0$, where all paths starting there pass through the space-time point $(x_i, t)$. At the boundary $u_i$
between two such intervals, we note however that $\Wi((u_i, 0))$ must contain at least two paths, which are limits of paths in $\Wi$ starting from $(u_{i-1}, u_i)$, resp.\ $(u_i, u_{i+1})$, at time $0$. Are there random space-time points where more than two paths originate? It turns out that we can give a complete classification of the type of multiplicity we see almost surely in a Brownian web. The main tool to accomplish this is the self-duality of the Brownian web discussed in Section~\ref{SS:webdual}.

First we give a classification scheme for $z\in\R^2$ according to the multiplicity of paths in $\Wi$ entering and leaving $z$. We say a path $\pi$
enters $z=(x,t)$ if $\sigma_\pi<t$ and $\pi(t)=x$, and $\pi$ leaves $z$ if $\sigma_\pi\leq t$ and $\pi(t)=x$. Two paths $\pi$ and $\pi'$ leaving
$z$ are defined to be equivalent, denoted by $\pi\sim^z_{\rm out}\pi'$, if $\pi=\pi'$ on $[t,\infty)$. Two paths $\pi$ and $\pi'$ entering $z$
are defined to be equivalent, denoted by $\pi\sim^z_{\rm in}\pi'$, if $\pi=\pi'$ on $[t-\eps,\infty)$ for some $\eps>0$. Note that $\sim^z_{\rm in}$ and $\sim^z_{\rm out}$ are equivalence relations.

Let $m_{\rm in}(z)$, resp.\ $m_{\rm out}(z)$, denote the number of equivalence
classes of paths in $\Wi$ entering, resp.\ leaving, $z$, and let $\hat m_{\rm in}(z)$ and $\hat m_{\rm out}(z)$ be defined similarly for the
dual Brownian web $\hat\Wi$. Given a realization of the Brownian web $\Wi$, points $z\in\R^2$ are classified according to the value of $(m_{\rm in}(z),m_{\rm out}(z))$. We divide points of type (1,2) further into types $(1,2)_{\rm l}$ and $(1,2)_{\rm r}$, where the subscript
$\rm l$ (resp.\ $\rm r$) indicates that the left (resp.\ right) of the two outgoing paths is the continuation of the (up to equivalence)
unique incoming path. Points in the dual Brownian web $\widehat\Wi$ are labeled according to their type in the Brownian web obtained by
rotating the graph of $\widehat\Wi$ in $\R^2$ by $180^{\rm o}$.

\begin{figure}[tp] 
\begin{center}
\includegraphics[width=16cm]{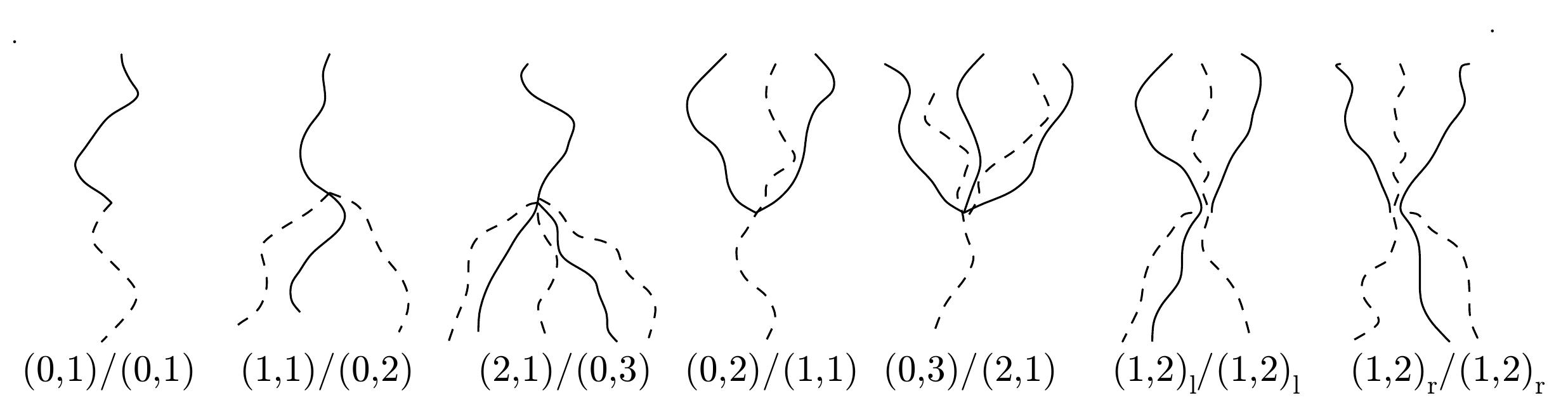}
\caption{Special points of the Brownian web.}
\label{fig:webspecialpts}
\end{center}
\end{figure}

We are now ready to state the following classification result (see also
\cite[Proposition~2.4]{TW98} and \cite[Theorems~3.11--3.14]{FINR06}), illustrated in Figure~\ref{fig:webspecialpts}.

\bt[Special points of the Brownian web]\label{T:classweb}
Let $(\Wi, \widehat \Wi)$ be the standard Brownian web and its dual. Then almost surely, each
$z\in\R^2$ satisfies
\begin{equation}\label{minout}
m_{\rm out}(z)= \hat m_{\rm in}(z)+1 \qquad \mbox{and}  \qquad \hat m_{\rm out}(z)=m_{\rm in}(z)+1,
\end{equation}
and $z$ is of one of the following seven types according to $(m_{\rm in}(z), m_{\rm out}(z))/(\hat m_{\rm in}(z), \hat m_{\rm out}(z))$:
$$
(0,1)/(0,1), \ (1,1)/(0,2),\ (0,2)/(1,1),\  (2,1)/(0,3),\ (0,3)/(2,1),\  (1,2)_{\rm l}/(1,2)_{\rm l},\  (1,2)_{\rm r}/(1,2)_{\rm r}.
$$
Almost surely,
\begin{itemize}
\item[\rm(i)] the set of points of type $(0,1)/(0,1)$ has full Lebesgue measure in $\R^2$;

\item[\rm(ii)] points of type $(1,1)/(0,2)$ are points in the set $\bigcup_{\pi \in \Wi} \{(\pi(t), t): t>\sigma_\pi\}$, excluding points of type $(1,2)/(1,2)$ and $(2,1)/(0,3)$;

\item[\rm(iii)] the set of points of type $(2,1)/(0,3)$ consists of points at which two paths in $\Wi$ coalesce and is countable;

\item[\rm(iv)] points of type $(1,2)_{\rm l}/(1,2)_{\rm l}$ are points of intersection between some $\pi\in\Wi$ and $\hat\pi\in\widehat\Wi$, with $\sigma_\pi<t<\hat\sigma_{\hat\pi}$, $\pi(s)\leq \hat\pi(s)$ for all $s\in [\sigma_\pi, \hat\sigma_{\hat\pi}]$, and $\pi$ intersects
    $\hat\pi$ at $(\pi(t),t)=(\hat\pi(t), t)$.
\end{itemize}
Similar statements hold for the remaining three types by symmetry.
\et
{\bf Proof.} We first prove relation \eqref{minout}. Let $z=(x,t)$, and assume that $\hat m_{\rm in}(z)=k$ for some $k\in \N_0:=\{0\}\cup\N$. Then there exist $\eps>0$ and $k$ ordered paths $\hat\pi_1, \ldots, \hat\pi_k\in\widehat\Wi$ starting at time $t+\eps$, such that these paths are disjoint on $(t,t+\eps]$ and coalesce together at time $t$ at position $x$. Note that the ordered paths $(\hat\pi_i)_{1\leq i\leq k}$ divide the space-time strip $\R\times (t, t+\eps)$ into $k+1$ regions $(I_i)_{1\leq i\leq k+1}$, where $I_1$ is the region to the left of $\hat\pi_1$,
$I_i$ is the region between $\hat\pi_{i-1}$ and $\hat\pi_{i}$ for each $2\leq i\leq k$, and $I_{k+1}$ is the region to the right of $\hat\pi_{k}$.
From the interior of each region $I_i$, we can pick a sequence of starting points $(z^n_i)_{n\in\N}$ with $z^n_i\to z$. Since paths in $\Wi$ and
$\widehat\Wi$ do not cross, as formulated in Theorem~\ref{T:dwebchar}, each path $\pi^n_i \in \Wi(z^n_i)$ must stay confined in $\overline{I_i}$ in the time interval $[t, t+\eps]$, and so must any subsequential limit of $(\pi^n_i)_{n\in\N}$. Such subsequential limits must exist by the almost sure compactness of $\Wi$, and each subsequential limit is a path $\pi_i\in \Wi(z)$. Therefore $\Wi(z)$ must contain at least $k+1$ distinct paths, one contained in $\overline{I_i}$ for each $1\leq i\leq k+1$. Furthermore, each $\overline{I_i}$ cannot contain more than one path in $\Wi(z)$.
Indeed, if $\overline{I_i}$ contains two distinct path $\pi, \pi'\in\Wi(z)$, then any path $\hat\pi\in \widehat\Wi$ started strictly between $\pi$ and $\pi'$ on the time interval $(t, t+\eps)$ must enter $z$, and $\hat\pi$ is distinct from $(\hat\pi_i)_{1\leq i\leq k}$, which contradicts the assumption that $\hat m_{\rm in}(z)=k$. Therefore we must have $m_{\rm out}(z)=k+1=\hat m_{\rm in}(z)+1$.

Similar considerations as above show that if $z$ is of type $(1,2)_{\rm l}$ in $\Wi$, then it must be of the same type in $\widehat\Wi$. The same holds for type $(1,2)_{\rm r}$.

We now show that the seven types of points listed are all there is. Note that it suffices to show that almost surely $m_{\rm in}(z)+m_{\rm out}(z)\leq 3$ for each $z\in\R^2$. There are four possible cases of $m_{\rm in}(z)+m_{\rm out}(z)>3$, which we rule out one by one:
\begin{itemize}
\item[\rm (a)] For some $z\in\R^2$, $m_{\rm in}(z)\geq 3$ and $m_{\rm out}(z)\geq 1$. Note that Corollary~\ref{C:pathconv} implies that every path $\pi\in\Wi$ coincides with some path in $\Wi(\Q^2)$ on $[\sigma_\pi+\eps, \infty)$, for any given $\eps>0$. Therefore the event that $m_{\rm in}(z)\geq 3$ for some $z\in\R^2$ is contained in the event that three distinct Brownian motions among $\Wi(\Q^2)$ coalesce at the same time. Such an event has probability zero, because there are countably many ways of choosing three Brownian motions form $\Wi(\Q^2)$, and conditioned on two Brownian motions coalescing at a given space-time point, there is zero probability that a third independent Brownian motion (which evolves independently before coalescing) would visit the same space-time point.

\item[\rm (b)] For some $z\in\R^2$, $m_{\rm in}(z)\geq 2$ and $m_{\rm out}(z)\geq 2$. In this case, $\hat m_{\rm in}(z)\geq 1$, and again by Corollary~\ref{C:pathconv}, the event we consider is contained in the event that there exist two paths $\pi_1, \pi_2\in \Wi(\Q^2)$ and a path $\hat\pi\in \widehat\Wi(\Q^2)$, such that $\hat\pi$ passes through the point of coalescence between $\pi_1$ and $\pi_2$. Such an event has probability 0, since conditioned on $\pi_1$ and $\pi_2$ up to the time of their coalescence, $\hat\pi$ is an independent Brownian motion with zero probability of hitting a given space-time point -- the point of coalescence between $\pi_1$ and $\pi_2$.

\item[\rm (c)] For some $z\in\R^2$, $m_{\rm in}(z)\geq 1$ and $m_{\rm out}(z)\geq 3$. In this case, $\hat m_{\rm out}(z)\geq 2$ and $\hat m_{\rm in}(z)\geq 2$, which is equivalent to Case (b) by the symmetry between $\Wi$ and $\widehat\Wi$.

\item[\rm (d)] For some $z\in\R^2$, $m_{\rm out}(z)\geq 4$. In this case, $\hat m_{\rm in}(z)\geq 3$, which is equivalent to Case (a).
\end{itemize}

We leave the verification of statements (i)--(iv) as an exercise.
%
%
%
%
\qed

\begin{exercise}
Verify statements (i)--(iv) in Theorem~\ref{T:classweb}.
\end{exercise}

\section{The Brownian net}\label{S:net}

\begin{figure}[tp] 
\begin{center}
\includegraphics[width=15cm]{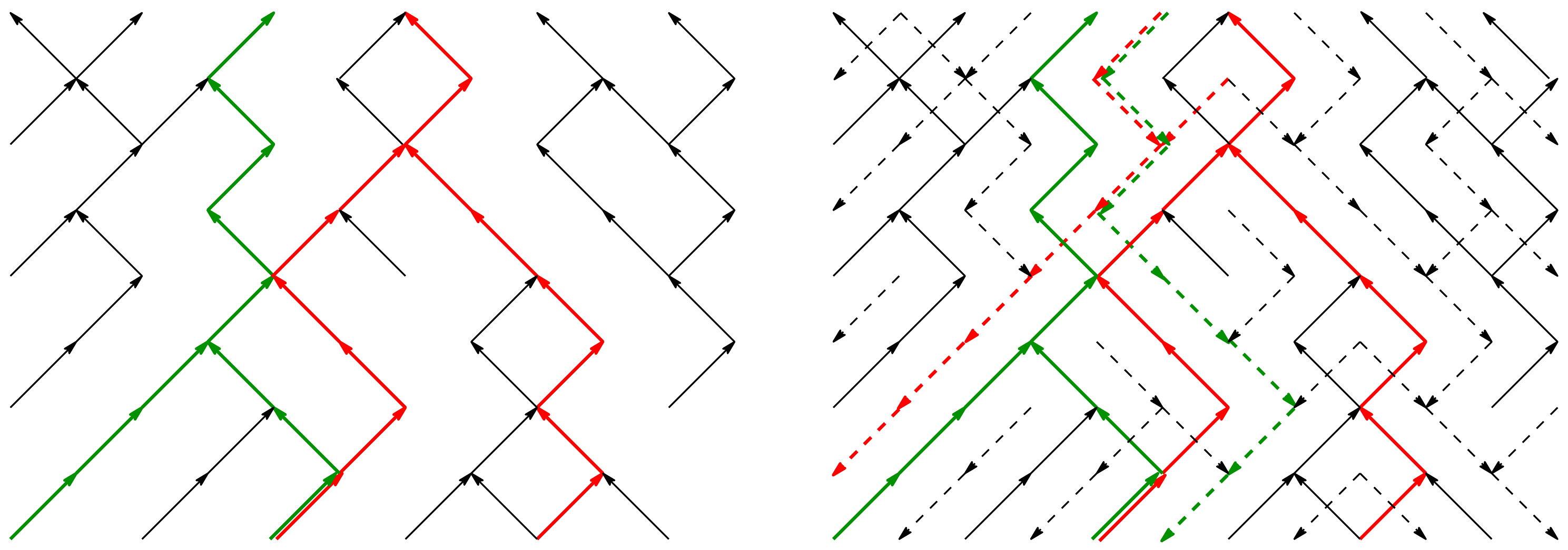}
\caption{Discrete space-time branching-coalescing random walks on $\Z^2_{\rm even}$, and its dual on $\Z^2_{\rm odd}$.}
\label{fig:netdual}
\end{center}
\end{figure}

The Brownian net generalizes the Brownian web by allowing paths to branch. The existence of
such an object is again motivated by its discrete analogue, the collection of discrete time branching-coalescing
simple symmetric random walks on $\Z$. Figure~\ref{fig:netdual} gives an illustration: from each
$(x,n)\in \Z^2_{\rm even}$, an arrow is drawn from $(x,n)$ to either $(x-1, n+1)$
or $(x+1, n+1)$ with probability $(1-\eps)/2$ each,
representing whether the walk starting at $x$ at time $n$ should move to $x-1$ or $x+1$ at time $n+1$;
and with probability $\eps$, arrows are drawn from $(x,n)$ to both $(x-1, n+1)$ and $(x+1, n+1)$, so that the
walk starting at $(x,n)$ branches into two walks, with one moving to $x-1$ and the other to $x+1$ at time $n+1$.

If we consider the collection of all upward random walk paths obtained by following the arrows, then the natural
question is: when space-time is scaled diffusively, could this random collection of paths have a non-trivial
limit? To have an affirmative answer to this question, it is necessary to choose the branching
probability $\eps$ to depend suitably on the diffusive scaling parameter. More precisely:

\begin{itemize}
\item[\bf Q.1] If space-time is rescaled by $S_\eps (x,t):=(\eps x, \eps^2 t)$, and the branching probability is chosen to
be $b\eps$ for some $b>0$, then what is the scaling limit of the collection of branching-coalescing random walk paths as $\eps\downarrow 0$?
\end{itemize}

The limit will be what we call the {\em Brownian net $\Ni_b$ with branching parameter $b$}. For simplicity, we will focus on the case $b=1$, with
the limit called the {\em standard Brownian net} $\Ni$. The fact that a non-trivial scaling limit exists with the above choice of the branching probability is hinted by the following observation. Instead of considering the collection of all random walk paths, let us first restrict our attention to
two special subsets: the set of {\em leftmost}, resp.\ {\em rightmost}, random walk paths where the random walk always follows the arrow
to the left, resp.\ right, whenever it encounters a branching point (see Figure~\ref{fig:netdual}). Note that the collection
of leftmost paths is a collection of coalescing random walks with drift $-\eps$, which ensures that each path under the diffusive scaling
$S_\eps (x,t)=(\eps x, \eps^2 t)$ converges to a Brownian motion with drift $-1$. Therefore we expect the collection of
leftmost paths to converge to a variant of the Brownian web, $\Wi^{\rm l}$, which consists of coalescing Brownian motions
with drift $-1$. Similarly, we expect the collection of rightmost paths to converge to a limit $\Wi^{\rm r}$, which
consists of coalescing Brownian motions with drift $+1$. Of course, $\Wi^{\rm l}$ and $\Wi^{\rm r}$ are coupled in
a non-trivial way.

The above observation explains the choice of the branching probability, and we see that any limit of the branching-coalescing
random walk paths must contain the two coupled Brownian webs $(\Wi^{\rm l}, \Wi^{\rm r})$. The questions that remain are:
\begin{itemize}
\item[\rm (A)] How to characterize the joint law of $(\Wi^{\rm l}, \Wi^{\rm r})$?

\item[\rm (B)] Can we construct the scaling limit of branching-coalescing random walks from $(\Wi^{\rm l}, \Wi^{\rm r})$?
\end{itemize}

To answer (A), it suffices to characterize the joint distribution $l_{z_1}\in \Wi^{\rm l}(z_1), \ldots, l_{z_k}\in \Wi^{\rm l}(z_k)$ and $r_{z'_1}\in \Wi^{\rm r}(z'_1), \ldots, r_{z'_{k'}}\in\Wi^{\rm r}(z'_{k'})$ for a finite collection of $(z_i)_{1\leq i\leq k}$ and $(z'_i)_{1\leq i\leq k'}$ in $\R^2$. An examination of their discrete analogue suggests that:
\begin{itemize}
\item the paths $(l_{z_1}, \ldots, l_{z_k}, r_{z'_1},\ldots, r_{z'_{k'}})$ evolve independently when they are apart;

\item the {\em leftmost paths} $(l_{z_1}, \ldots, l_{z_k})$ coalesce when they meet, and the same is true for the {\em rightmost paths} $(r_{z'_1},\ldots, r_{z'_{k'}})$;

\item a pair of leftmost and rightmost paths $(l_{z_i}, r_{z'_j})$ solves the following pair of SDEs:
\bc\label{coupsde}
\dis\di L_t&=&\dis 1_{\{L_t\neq R_t\}}\di B^{\rm l}_t
+1_{\{L_t=R_t\}}\di B^{\rm s}_t-\di t,\\[5pt]
\dis\di R_t&=&\dis 1_{\{L_t\neq R_t\}}\di B^{\rm r}_t
+1_{\{L_t=R_t\}}\di B^{\rm s}_t+\di t,
\ec
where the leftmost path $L$ and the rightmost path $R$ are driven by independent Brownian motions
$B^{\rm l}$ and $B^{\rm r}$ when they are apart, and driven by the same Brownian motion $B^{\rm s}$ (independent
of $B^{\rm l}$ and $B^{\rm r}$) when they coincide; furthermore, $L$ and $R$ are subject to the constraint that $L_t\leq R_t$ for all
$t\geq T:=\inf\{u\geq \sigma_L \vee \sigma_R: L_u\leq R_u\}$, with $\sigma_L$ and $\sigma_R$ being the starting times of
$L$ and $R$.
\end{itemize}
It turns out that the SDE \eqref{coupsde} has a unique weak solution, and the above properties uniquely determine the joint law of
$(l_{z_1}, \ldots, l_{z_k}, r_{z'_1},\ldots, r_{z'_{k'}})$, which we will call {\em left-right coalescing Brownian motions}. Extending
the starting points to a countable dense set in $\R^2$, and then taking closure of the resulting set of leftmost, resp.\ rightmost
paths a.s.\ determines $(\Wi^{\rm l}, \Wi^{\rm r})$, which we will call the {\em left-right Brownian web}.

To answer (B), observe that in the discrete case, all random walk paths can be constructed by hopping back and forth between
leftmost and rightmost random walk paths. This suggests a similar approach to construct the scaling limit of
the set of all branching-coalescing random walk paths, which we will call {\em the Brownian net} $\Ni$. More precisely, to construct
$\Ni$ from $(\Wi^{\rm l}, \Wi^{\rm r})$, we simply consider the set of all paths that can be obtained by {\em hopping} a finite number of times between paths in $\Wi^{\rm l}$ and $\Wi^{\rm r}$, and then take its closure.

The above considerations led to the original construction of the Brownian net $\Ni$ in~\cite{SS08}, called the {\em hopping construction}.

From Figure~\ref{fig:netdual}, it is easily seen that the branching-coalescing random walks on $\Z^2_{\rm even}$ a.s.\ uniquely determine a dual collection of branching-coalescing random walks on $\Z^2_{\rm odd}$, running backward in time. Furthermore, the two systems are equally distributed apart from a rotation in space-time by $180^o$ and a lattice shift. Therefore in the scaling limit, we expect the left-right Brownian web
$(\Wi^{\rm l}, \Wi^{\rm r})$ to have a dual $(\widehat\Wi^{\rm l}, \widehat\Wi^{\rm r})$, which determines a dual Brownain net $\widehat \Ni$.
As for the Brownian web, such a duality provides a powerful tool. In particular, it leads to a second construction of the Brownian net, called the {\em wedge construction} in~\cite{SS08}.

Besides the hopping and wedge constructions of the Brownian net, there are two more constructions, called the {\em mesh construction}, also developed in~\cite{SS08}, and the {\em marking construction} developed by Newman, Ravishankar and Schertzer
in~\cite{NRS10}, where the Brownian net was conceived independently from~\cite{SS08}. The mesh construction is based on the observation that,
given the left-right Brownian web $(\Wi^{\rm l}, \Wi^{\rm r})$, there exist space-time regions (called {\em meshes}) with their left boundaries being rightmost paths, and their right boundaries being leftmost paths. Such unusual configurations makes these meshes forbidden regions, where no paths can enter. The mesh construction asserts that the Brownian net consists of all paths which do not enter meshes.

In contrast to the hopping construction, which is an {\em outside-in} approach where the Brownian net is constructed from its outermost paths -- the leftmost and rightmost paths, the {\em marking construction} developed in~\cite{NRS10} is an {\em inside-out} approach, where one
starts from a Brownian web and then constructs the Brownian net by adding branching points. In the discrete setting, this amounts to turning coalescing random walks into branching-coalescing random walks by changing each lattice point independently into a branching point with probability $\eps$. In the continuum setting, this turns out to require {\em Poisson marking} the set of $(1,2)$ points of the Brownian web (cf.~Theorem~\ref{T:classweb}) and turning them into branching points, so that the incoming Brownian web path can continue along either of the two outgoing Brownian web paths. We will introduce the marking construction in detail in Section~\ref{S:couple}, where we will study couplings between the Brownian web and net.

In the rest of the this section, we will define the left-right Brownian web $(\Wl, \Wr)$, give the hopping, wedge and mesh constructions of the Brownian net, and study various properties of the Brownian net, including the branching-coalescing point set, the backbone of the Brownian net, and special points of the Brownian net.

\subsection{The left-right Brownian web and its dual}\label{S:lrweb}

The discussions above show that the key object in the construction of the Brownian net $\Ni$ is the left-right Brownian web $(\Wl, \Wr)$, which should be the diffusive scaling limit of the collections of leftmost and rightmost branching-coalescing random walk paths with branching probability $\eps$. In turn, the key ingredient in the construction of $(\Wl, \Wr)$ is the pair of left-right SDEs in \eqref{coupsde}, which can be shown to be well-posed.

\bp{\bf(The left-right SDE)}\label{P:sdeuni} For each initial state $(L_0,R_0)\in\R^2$, there exists a unique weak
solution to the SDE \eqref{coupsde} subject to the constraint that $L_t\leq R_t$ for all $t\geq T:=\inf\{s\geq 0:L_s=R_s\}$. Furthermore,
almost surely, if $I:=\{t\geq T: L_t=R_t\}\neq\emptyset$, then $I$ is nowhere dense perfect set with positive Lebesgue measure.
\ep
{\bf Proof Sketch.} We sketch the basic idea and refer to~\cite[Prop.~2.1 \& 3.1]{SS08} for details. Assume w.l.o.g.\ that $L_0=R_0=0$. Define
$$
T_t:= \int_0^t 1_{\{ L_s<R_s\}}{\rm d}s, \qquad S_t:= \int_0^t 1_{\{ L_s= R_s\}}{\rm d}s,
$$
and
$$
\tilde B^{\rm l}_{T_t}:= \int_0^t 1_{\{L_s<R_s\}} {\rm d}B^{\rm l}_s, \quad \tilde B^{\rm r}_{T_t}:= \int_0^t 1_{\{L_s<R_s\}} {\rm d}B^{\rm r}_s, \quad \tilde B^{\rm s}_{S_t}:= \int_0^t 1_{\{L_s=R_s\}} {\rm d}B^{\rm s}_s.
$$
Then $(L_t, R_t)$ solves
\begin{equation}\label{LRt}
\begin{aligned}
L_t  = \tilde B^{\rm l}_{T_t} + \tilde B^{\rm s}_{S_t} -t, \\
R_t  = \tilde B^{\rm r}_{T_t} + \tilde B^{\rm s}_{S_t} +t, \\
\end{aligned}
\end{equation}
and the difference $D_t:=R_t-L_t$ satisfies
\begin{equation}\label{Dt}
D_t  = (\tilde B^{\rm r}-\tilde B^{\rm r})_{T_t} + 2t = \sqrt{2} \tilde B_{T_t} + 2T_t +2S_t
\end{equation}
with the constraint that $D$ is non-negative, where $\tilde B:= \frac{1}{\sqrt 2} (\tilde B^{\rm r}-\tilde B^{\rm l})$
is also a standard Brownian motion.

Since $T_t = \int_0^t 1_{\{ L_s<R_s\}}{\rm d}s = \int_0^t 1_{\{ D_s>0\}}{\rm d}s$, it is easily seen that $T_t$ must be continuous and
strictly in creasing in $t$. Therefore $\tau:=T_t$ admits an inverse $T^{-1}\tau =t$. Rewriting the equation \eqref{Dt} for $D$ with
respect to the variable $\tau$, the time $D$ spent at the origin, we obtain
\begin{equation}
\tilde D_\tau :=\frac{1}{\sqrt 2}D_{T^{-1}\tau} = \tilde B_\tau + \sqrt{2} \tau + \sqrt{2}S_{T^{-1}\tau},
\end{equation}
where $\tilde D_\tau$ can be regarded as a transformation of $\tilde B_\tau + \sqrt{2} \tau$ by adding an increasing function $\tilde S_\tau:=\sqrt{2}S_{T^{-1}\tau}$, which increases only when $\tilde D_\tau=0$ such that $\tilde D_\tau$ stays non-negative. Such an equation is known as a {\rm Skorohod equation}, with $\tilde D_t$ being the Skorohod reflection of $\tilde B_\tau +\sqrt{2} \tau$ at the origin. Such a Skorohod equation admits a pathwise unique solution~\cite[Section 3.6.C]{KS91}, with
$$
\tilde S_\tau:=\sqrt{2}S_{T^{-1}\tau} = -\inf_{0\leq s\leq \tau}(\tilde B_\tau + \sqrt{2} \tau),
$$
which is in fact also the local time at origin for the drifted Brownian motion $\tilde B_\tau + \sqrt{2} \tau$ reflected at the origin.

Having determined $S_{T^{-1}\tau}$ and $D_{T^{-1}\tau}$ almost surely from $\tilde B=\frac{1}{\sqrt 2} (\tilde B^{\rm r}-\tilde B^{\rm l})$, to recover $D_t$, we only need to make a time change from $\tau$ back to $t:= T_t+ S_t=\tau+S_{T^{-1}\tau}$. Note that this time change has no effect when $D$ is away from $0$, but adds positive Lebesgue time when $D$ is at $0$. Therefore in contrast to $D_{T^{-1}\tau}$, which is a drifted Brownian motion reflected instantaneously at the origin, $D_t$ is the same Brownian motion {\em sticky reflected} at the origin (see e.g.~\cite{W02} and the references therein for further details on sticky reflected Brownian motions). Similarly, from $T_t=\tau$ and $S_t=S_{T^{-1}\tau}$, we can construct $(L_t, R_t)$ in \eqref{LRt}. From the same arguments, it is also easily seen that $I:=\{t\geq 0: L_t=R_t\}$ is almost surely a nowhere dense perfect set with positive Lebesgue measure.
\qed

Having characterized the interaction of a single pair of leftmost and rightmost paths, we can now construct a collection of {\em left-right coalescing Brownian motions} $(l_{z_1}, \ldots, l_{z_k}, r_{z'_1},\ldots, r_{z'_{k'}})$ with the properties that: (1) the paths evolve independently when they do not coincide; (2) $(l_{z_1}, \ldots, l_{z_k})$, resp.\ $(r_{z'_1},\ldots, r_{z'_{k'}})$, is distributed as a collection of coalescing Brownian motions with drift $-1$, resp.\ $+1$; (3) every pair $(l_{z_i}, r_{z_j})$ is a weak solution to the left-right SDE \eqref{coupsde}. The construction can be carried out inductively. Assume w.l.o.g.\ that the paths all start at the same time. Then
\begin{itemize}
\item Let the paths evolve independently until the first time two paths meet.

\item If this pair of paths are of the same type, then let them coalesce and iterate the construction with one path less than before.

\item If this pair of paths are of different types, then let them evolve as a left-right pair solving the SDE \eqref{coupsde}, and let all
other paths evolve independently, until the first time two paths (other than paths in the same left-right pair) meet.

\item If the two meeting paths are of the same type, then let them coalesce and iterate the construction.

\item If the two meeting paths are of different types, then let them form a left-right pair (breaking whatever pair relations they were in),
and iterate the construction.
\end{itemize}
It is easily seen that this iterative construction terminates after a finite number of steps, when either a single path, or a single pair of leftmost and rightmost paths remains.
\medskip

We are now ready to characterize the {\em left-right Brownian web} $(\Wl, \Wr)$.
\bt{\bf(The left-right Brownian web and its dual)}\label{T:lrweb}
There exists an $\Hi^2$-valued random variable $(\Wl,\Wr)$, called the (standard) left-right Brownian web, whose distribution is
uniquely determined by the following properties:
\begin{itemize}
\item[\rm (i)] The left Brownian web $\Wl$ $($resp.\ right Brownian web $\Wr$$)$ is distributed as a Brownian web $\Wi$ tilted with drift $-1$ $($resp.\ $+1$$)$, i.e., $\Wl$ $($resp.\ $\Wr$$)$ has the same distribution as the image of $\Wi$ under the space-time transformation $(x,t)\to (x- t, t)$ $($resp.\ $(x,t)\to (x+ t, t)$$)$.

\item[\rm (ii)] For any finite deterministic set of points $z_1,\ldots,z_k,z'_1,\ldots,z'_{k'}\in\R^2$, the collection of paths
$(\lp_{z_1}, \cdots, \lp_{z_k}; \rp_{z'_1}, \cdots, \rp_{z'_{k'}})$ is distributed as a family of left-right coalescing Brownian
motions.
\end{itemize}
Furthermore, almost surely there exists a dual left-right Brownian web $(\widehat \Wi^{\rm l}, \widehat \Wi^{\rm r})\in \widehat\Hi^2$, such
that $(\Wl, \widehat\Wi^{\rm l})$ $($resp.\ $(\Wr, \widehat\Wi^{\rm r})$$)$ is distributed as $(\Wi, \widehat\Wi)$  tilted with drift $-1$ $($resp.\ $+1)$, and $(-\widehat \Wi^{\rm l}, -\widehat \Wi^{\rm r})$ has the same distribution as $(\Wl, \Wr)$.
\et
{\bf Proof Sketch.} The existence and uniqueness of a left-right Brownian web satisfying properties (i)--(ii) follow the same argument as for the Brownian web in Theorem~\ref{T:webchar}. The almost sure existence of a dual left-right Brownian web follows from the duality of the Brownian web. The fact that $(\widehat \Wi^{\rm l}, \widehat \Wi^{\rm r})$ has the same distribution as $(\Wl, \Wr)$, except for rotation around the origin in space-time by $180^o$, can be derived by taking the diffusive scaling limits of their discrete counterparts, where such an equality in distribution is trivial. See~\cite[Theorem~5.3]{SS08} for further details.
\qed

\begin{exercise}\label{E:lrwebcross}
Show that a.s., no path in $\Wl$ can cross any path in $\Wr$ or $\widehat\Wi^{\rm r}$ from left to right, where $\pi_1$
is said to cross $\pi_2$ from left to right if there exist $s<t$ such that $\pi_1(s)<\pi_2(s)$ and $\pi_1(t)>\pi_2(t)$. Similarly, paths in $\Wr$ cannot cross paths in $\Wl$ or $\widehat\Wi^{\rm l}$ from right to left.
\end{exercise}

\subsection{The hopping construction of the Brownian net}

We are now ready to construct the Brownian net $\Ni$ by allowing paths to hop back and forth between paths in the left-right Brownian web $(\Wl, \Wr)$, where a path $\pi$ obtained by hopping from $\pi_1$ to $\pi_2$ at time $t$, with $\pi_1(t)=\pi_2(t)$, is defined by $\pi:=\pi_1$ on $[\sigma_{\pi_1}, t]$ and $\pi:=\pi_2$ on $[t,\infty)$.

Given $\pi_1$ and $\pi_2$ with $\pi_1(t)=\pi_2(t)$, which are in the scaling limit $\Ni$ of branching-coalescing random walk paths, it is not guaranteed that $\pi$ obtained by hopping from $\pi_1$ to $\pi_2$ at time $t$ is also in $\Ni$. Indeed, even though $\pi_1(t)=\pi_2(t)$,
$\pi_1$ and $\pi_2$ may still arise as limits of random walk paths which do not meet, for which hopping is not possible. Therefore $\pi$ has
no approximating analogue among the branching-coalescing random walk paths, and hence $\pi$ may not be in $\Ni$. One remedy is to allow hopping
from $\pi_1$ to $\pi_2$ at time $t$ only if the two paths {\em cross at time $t$}, i.e., $\sigma_{\pi_1}, \sigma_{\pi_2}<t$, there exist times
$t^-<t^+$ with $(\pi_1(t^-)-\pi_2(t^-))(\pi_1(t^+)-\pi_2(t^+))<0$, and $t=\inf\{s\in (t^-, t^+): (\pi_1(t^-)-\pi_2(t^-))(\pi_1(s)-\pi_2(s))<0\}$.
We call $t$ the {\em crossing time} between $\pi_1$ and $\pi_2$. If $\pi_1, \pi_2\in\Ni$ and they cross at time $t$, then it is easily seen by discrete approximations that the path $\pi$ obtained by hopping from $\pi_1$ to $\pi_2$ at time $t$ must also be in $\Ni$.

Given a set of paths $K$, let $\Hi_{\rm cross}(K)$ denote the set of paths obtained by hopping a finite number of times among paths in $K$ at crossing times. The Brownian net $\Ni$ can then be constructed by setting $\Ni:= \overline{\Hi_{\rm cross}(\Wl\cup\Wr)}$. Here is the hopping characterization of the Brownian net from~\cite[Theorem~1.3]{SS08}.

\bt[Hopping Characterization of the Brownian net]\label{T:netchar}
There exists an $(\Hi, \Bi_\Hi)$-valued random variable $\Ni$, called the standard Brownian net, whose distribution is uniquely determined
by the following properties:
\begin{itemize}
\item[\rm (i)] For each $z\in\R^2$, $\Ni(z)$ a.s.\ contains a
unique left-most path $l_z$ and right-most path $r_z$.\med

\item[\rm (ii)] For any finite deterministic set of points
$z_1,\ldots,z_k,z'_1,\ldots,z'_{k'}\in\R^2$, the collection of paths
$(l_{z_1},\ldots,l_{z_k},r_{z'_1},\ldots,r_{z'_{k'}})$ is distributed
as a family of left-right coalescing Brownian motions.\med

\item[\rm (iii)] For any deterministic countable dense sets
$\Di^{\rm l},\Di^{\rm r}\sub\R^2$,
\be\label{hopchar}
\Ni=\ov{\Hi_{\rm cross}\big(\{l_z:z\in\Di^{\rm l}\}
\cup\{r_z:z\in\Di^{\rm r}\}\big)}\quad{\rm a.s.}
\ee
\end{itemize}
\et
{\bf Proof.} The uniqueness in law of a random variable $\Ni$ satisfying the above properties is easily verified by the same argument as that for the Brownian web in Theorem~\ref{T:webchar}. For existence, we can just define $\Ni:= \overline{\Hi_{\rm cross}(\Wl\cup\Wr)}$, where
$(\Wl, \Wr)$ is the standard left-right Brownian web.

To show that $\Ni$ satisfies properties (i)--(iii), first note that $\Hi_{\rm cross}(\Wl\cup\Wr)$ satisfies properties (i)--(ii), where for each deterministic $z$, $\Hi_{\rm cross}(\Wl\cup\Wr)$ contains a leftmost element
and a rightmost element, which are just $l_z\in \Wl(z)$ and $r_z\in\Wr(z)$. We claim that taking closure of $\Hi_{\rm cross}(\Wl\cup\Wr)$
does not change the leftmost and rightmost element starting from any given $z=(x,t)$.

Indeed, if $\Ni(z)$ contains any path $\pi$ with $\pi(t')<l_z(t')-\eps$ for some $t'>t$ and $\eps>0$, then there exists a sequence $\pi_n\in \Hi_{\rm cross}(\Wl\cup\Wr)$ starting
from $z_n$, such that $\pi_n\to \pi$ and $\pi_n(t')\leq l_z(t')-\eps$ for all $n$ large. Since $l_{z_n}\in \Wl(z_n)$ is the leftmost path
in $\Hi_{\rm cross}(\Wl\cup\Wr)$ starting from $z_n$, we have $l_{z_n}\leq \pi_n$, and hence also $l_{z_n}(t')\leq l_z(t')-\eps$ for all
$n$ large. However, this is impossible because $l_{z_n}\to l_z$, and hence the time of coalescence between $l_{z_n}$ and $l_z$ tends to $t$
as $n\to\infty$ by Corollary~\ref{C:pathconv}.

We can therefore conclude that $\Ni:= \overline{\Hi_{\rm cross}(\Wl\cup\Wr)}$ also satisfies properties (i)--(ii). Since any path in $\Wl\cup\Wr$
can be approximated by paths in $\Wl(\Di^{\rm l})\cup\Wr(\Di^{\rm r})$ in the strong sense as in Corollary~\ref{C:pathconv}, it is easily seen that $\Hi_{\rm cross}(\Wl\cup\Wr) \subset \overline{\Hi_{\rm cross}(\Wl(\Di^{\rm l})\cup\Wr(\Di^{\rm r}))}$, and hence $\Ni$ also satisfies property (iii).

What we have left out in the proof is the a.s.\ pre-compactness of $\Hi_{\rm cross}(\Wl\cup\Wr)$, which is needed for $\Ni$ to qualify as
an $(\Hi, \Bi_\Hi)$-valued random variable. We leave this as an exercise.
\qed

\begin{exercise}
Show that almost surely, $\Hi_{\rm cross}(\Wl\cup\Wr)$ is pre-compact.
\end{exercise}

\begin{exercise}\label{E:netcross}
Show that a.s., no path in $\Ni$ can cross any path in $\Wr$ or $\widehat\Wi^{\rm r}$ from left to right, or cross any path in $\Wl$ or $\widehat\Wi^{\rm l}$ from right to left, where the definition of crossing is as in Exercise~\ref{E:lrwebcross}.
\end{exercise}

\br
Of course, we need to justify that the Brownian net characterized in Theorem~\ref{T:netchar} is indeed the scaling limit of branching-coalescing random walks, as motivated at the start of Section~\ref{S:net}. Such a convergence result was established in \cite[Section 5.3]{SS08}.
\er

\subsection{The wedge construction of the Brownian net}

The wedge and mesh constructions of the Brownian net $\Ni$ are both based on the observation that there are certain forbidden regions in space-time where Brownian net paths cannot enter. It turns out that the Brownian net can also be characterized as the set of paths that do not enter these forbidden regions. In the wedge construction, these forbidden regions, called {\em wedges}, are defined from the dual left-right Brownian web $(\widehat \Wi^{\rm l}, \widehat \Wi^{\rm r})$, while in the mesh construction, these forbidden regions, called {\em meshes}, are defined from the left-right Brownian web $(\Wl, \Wr)$.

\begin{figure}[tp] 
\centering
\plaat{5cm}{\includegraphics[width=3.5cm]{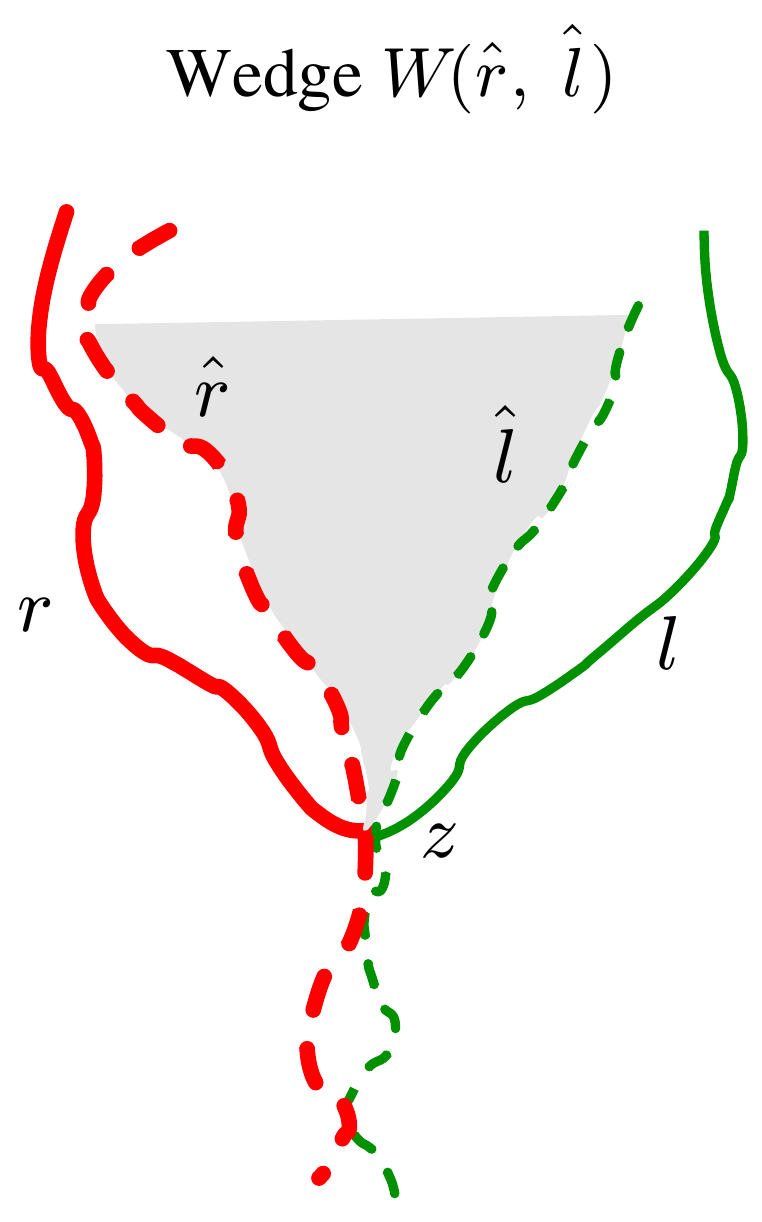}}
\caption{A wedge $W(\hat\rp,\hat\lp)$ with bottom point $z$.}
\label{fig:wedge}
\end{figure}

\bd[Wedges]
Let $(\Wl, \Wr, \widehat\Wi^{\rm l}, \widehat\Wi^{\rm r})$ be the standard left-right Brownian web and its dual.
For any $\hat r\in\widehat\Wi^{\rm r}$ and $\hat l\in\widehat\Wi^{\rm l}$ that are ordered with $\hat r(s)<\hat l(s)$ at the
time $s:=\final_{\hat r}\wedge\final_{\hat l}$, let $T:=\sup\{t<s:\hat r(t)=\hat l(t)\}$ $($possibly equals $-\infty)$ be the first
hitting time of $\hat r$ and $\hat l$. We call the open set $($see Figure~\ref{fig:wedge}$)$
\be\label{wedge}
W=W(\hat r,\hat l):=\big\{(x,u)\in\R^2:T<u<s,\ \hat r(u)<x<\hat l(u)\big\}
\ee
a {\em wedge} of $(\widehat\Wi^{\rm l}, \widehat\Wi^{\rm r})$ with left and right boundary $\hat r$ and $\hat l$ and bottom point $z:=(\hat r(T), T)=(\hat l(T), T)$. A path $\pi$ is said to {\em enter $W$ from outside} if there exist $\sigma_\pi\leq s<t$ such that $(\pi(s), s)\not\in \overline{W}$ and $(\pi(t), t)\in W$.
\ed

\bt[Wedge characterization of the Brownian net]\label{T:dualchar}
Let $(\Wl,\Wr,\widehat\Wi^{\rm l},\widehat\Wi^{\rm r})$ be the standard left-right Brownian web and its dual.
Then almost surely,
\be\label{dualchar}
\Ni=\big\{\pi\in\Pi:\pi\mbox{ does not enter any wedge of $(\widehat\Wi^{\rm l},\widehat\Wi^{\rm r})$
{f}rom outside}\big\}
\ee
is the standard Brownian net associated with $(\Wl,\Wr)$, i.e., $\Ni= \overline{\Hi_{\rm cross}(\Wl\cup\Wr)}$.
\et
\br\label{R:dualchar} The wedge characterization can also be applied to the Brownian web $\Wi$ with both $\widehat\Wi^{\rm l}$ and $\widehat\Wi^{\rm r}$
replaced by the dual Brownian web $\widehat\Wi$. Indeed, $\Wi$ can be seen as a degenerate Brownian net $\Ni_b$ with branching parameter
$b=0$, where $\Ni_b$ can be constructed in the same way as the standard Brownian net $\Ni$ with $b=1$, except that the left-right coalescing
Brownian motions in $(\Wl, \Wr)$ now have drift $\mp b$ respectively. For $\Wi$, the wedge characterization is stronger than requiring paths not to cross any path in $\widehat \Wi$ (cf.~Remark~\ref{R:webdualchar}), because it also prevents paths from entering a wedge from outside through its bottom point.
\er

\noindent
{\bf Proof of Theorem~\ref{T:dualchar}.} First we show that no path in $\overline{\Hi_{\rm cross}(\Wl\cup\Wr)}$ can enter any wedge of
$(\widehat\Wi^{\rm l},\widehat\Wi^{\rm r})$ from outside. If this is false, then there must be some path $\pi\in\Hi_{\rm cross}(\Wl\cup\Wr)$
which enters a wedge $W(\hat r, \hat l)$ from outside. There are two possibilities: either $\pi$ enters $W$ from outside by crossing
one of its two boundaries, which is impossible by Exercise~\ref{E:netcross}; or $\pi$ enters $W$ from outside through its bottom point $z$.
However, by the same argument as why a point of coalescence between two dual Brownian web paths cannot be hit by a forward Brownian web path
(cf.~Theorem~\ref{T:classweb}), no path in $\Wl\cup\Wr$ can enter the bottom point of a wedge $W$, and hence neither can any path in
$\Hi_{\rm cross}(\Wl\cup\Wr)$. This verifies the desired inclusion.

\begin{figure}[tp] 
\centering
\plaat{6cm}
{\includegraphics[width=7cm,height=6cm]{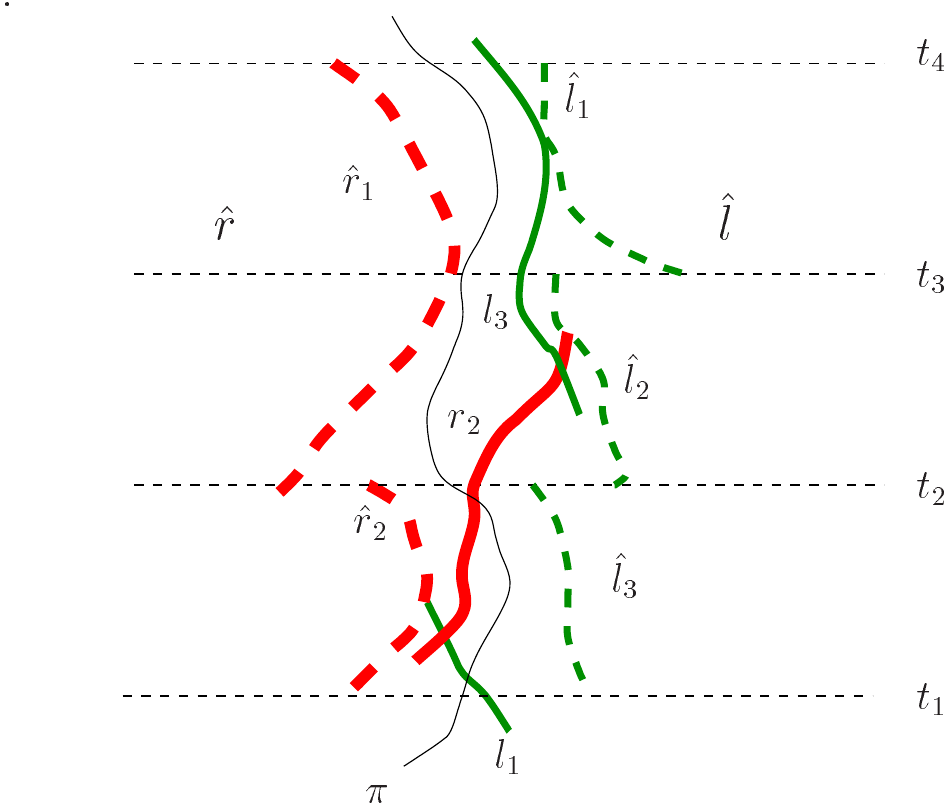}}
\caption{Steering a hopping path in the `fish-trap' $(\hat r,\hat l)$.}
\label{fig:funnel}
\end{figure}

We now show the converse inclusion that any path not entering wedges from outside must be in $\overline{\Hi_{\rm cross}(\Wl\cup\Wr)}$.
Let $\pi$ be such a path. The strategy to approximate $\pi$ by hopping paths is illustrated in Figure~\ref{fig:funnel}.

To approximate $\pi$ in a given time interval, say $[s, t]$, we first partition $[s,t]$ into sub-intervals
of equal length, $[t_{i-1}, t_{i}]$, for $1\leq i\leq N$. Fix an $\eps>0$. From the top time $t_N$, we consider the wedge formed by
$\hat r_1\in \widehat\Wi^{\rm r}$ and $\hat l_1\in \widehat\Wi^{\rm l}$ starting respectively at $\pi(t_N)-\eps$ and $\pi(t_N)+\eps$
at time $t_N$. Note that $\hat r_1$ and $\hat l_1$ cannot meet during the interval $[t_{N-1}, t_N]$, otherwise $\pi$ would be entering
the wedge $W(\hat r_1, \hat l_1)$ from outside. At time $t_{N-1}$, we check whether $\hat r_1<\pi(t_{N-1})-\eps$, and if it is the case,
then we start $\hat r_2$ at $\pi(t_{N-1})-\eps$ at time $t_{N-1}$. Similarly we check whether $\hat l_1>\pi(t_{N-1})+\eps$, and if it is
the case, then start $\hat l_2$ at $\pi(t_{N-1})+\eps$ at time $t_{N-1}$. In any event, we still have a pair of dual left-right paths which
encloses $\pi$ on the interval $[t_{N-2}, t_{N-1}]$, which are within $\eps$ distance of $\pi$ at the top time $t_{N-1}$. This procedure
is iterated until the time interval $[t_0, t_1]$, and it constructs a `fish-trap' of dual left-right paths, which can now be used to construct
a forward hopping path that stays inside the `fish-trap'.

Indeed, a forward path $l\in\Wl$ cannot cross dual paths
$\hat l_i\in\widehat\Wi^{\rm l}$, which form the right boundary of the `fish-trap'. When $l$ hits the left-boundary of the `fish-trap', we
can then hop to a path $r\in \Wr$ until it hits the right-boundary of the `fish-trap'. Iterating this procedure then constructs a hopping
path that stays inside the `fish-trap'. The almost sure equi-continuity of paths in $\Wl \cup\Wr \cup \widehat \Wi^{\rm l}\cup \widehat\Wi^{\rm r}$ ensures that only a finite number of hoppings is needed to reach time $t_N$, and the supnorm distance on the interval $[s,t]$ between $\pi$
and any path inside the `fish-trap' can be made arbitrarily small by choosing $N$ large and $\eps$ small. Therefore $\pi$ can be approximated
arbitrarily well by paths in $\Hi_{\rm cross}(\Wl\cup\Wr)$. For further details, see~\cite[Lemma 4.7]{SS08}.
\qed

\subsection{The mesh construction of the Brownian net}

\begin{figure}[tp] 
\centering
\plaat{5cm}{\includegraphics[width=3.5cm]{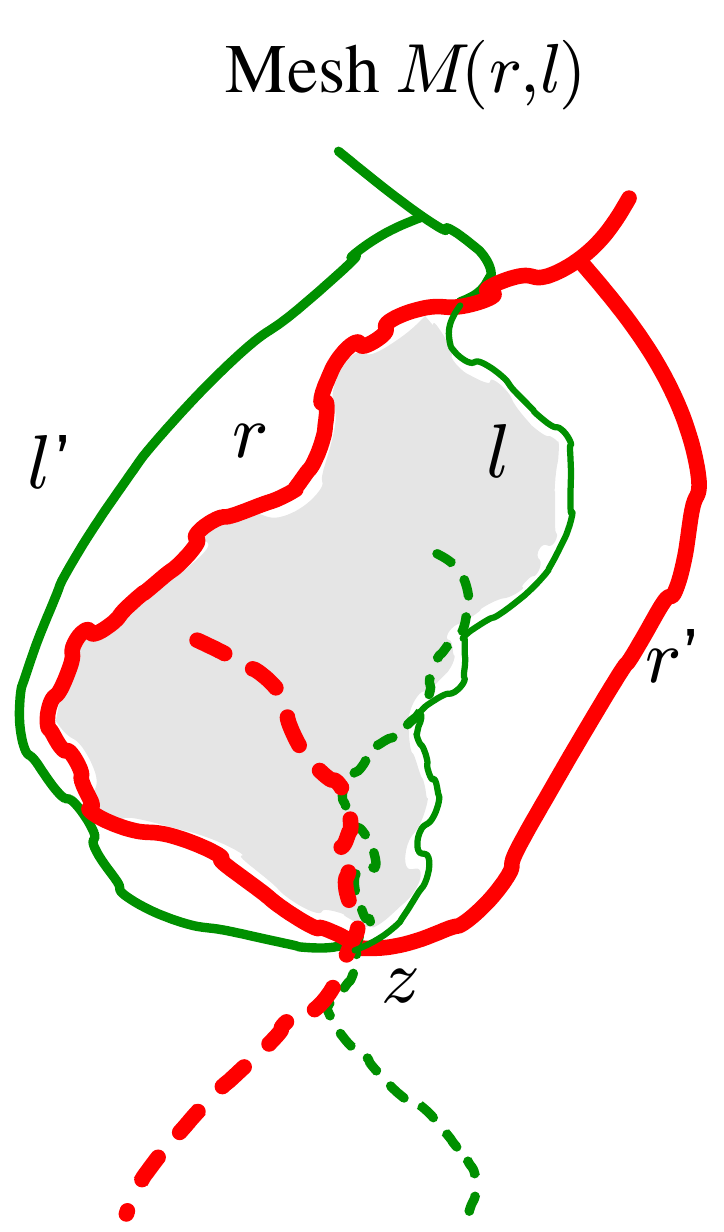}}
\caption{A mesh $M(\rp,\lp)$ with bottom point $z$.}
\label{fig:mesh}
\end{figure}

\bd[Meshes]
Let $(\Wl, \Wr)$ be the standard left-right Brownian web. If for a given $z=(x,t)\in\R^2$, there exist $\lp\in\Wl(z)$ and
$\rp\in\Wr(z)$ such that $\rp(s)<\lp(s)$ on $(t,t+\eps)$ for some $\eps>0$, then denoting $T:=\inf\{s>t:\rp(s)=\lp(s)\}$, we
call the open set $($see Figure~\ref{fig:mesh}$)$
\be
M=M(\rp,\lp):=\big\{(y,s)\in\R^2:t<s<T,\ \rp(s)<y<\lp(s)\big\}
\ee
a {\em mesh} of $(\Wl, \Wr)$ with left and right boundary $\rp$ and $\lp$ and bottom point $z$. A path $\pi$
is said to {\em enter $M$} if there exist $\start_\pi<s<t$ such that $(\pi(s), s)\not\in M$ and $(\pi(t), t)\in M$.
\ed

\bt[Mesh characterization of the Brownian net]\label{T:meshchar}
Let $(\Wl,\Wr)$ be the standard left-right Brownian web. Then almost surely,
\be\label{meshchar}
\Ni=\{\pi\in\Pi: \pi \mbox{ does not enter any mesh of } (\Wl,\Wr)\}
\ee
is the standard Brownian net associated with $(\Wl,\Wr)$, i.e., $\Ni= \overline{\Hi_{\rm cross}(\Wl\cup\Wr)}$.
\et
\br The mesh characterization can also be applied to the Brownian web $\Wi$, where the bottom point of the mesh must be of either
type $(0,2)$, $(1,2)$, or $(0,3)$ in Theorem~\ref{T:classweb}.
\er

We note that there is a subtle difference between {\em a path entering a mesh $M$ from outside}, vs {\em a path
entering a mesh $M$}. In particular, a path entering $M$ (but not entering $M$ from outside) could start inside $\overline{M}$,
hit the boundary of $M$ at a later time, and then move inside $M$. The heart of the proof of Theorem~\ref{T:meshchar} consists
in ruling out such scenarios for Brownian net paths, for which meshes play an essential role. As a by-product, one can show
the following result (see~\cite[Prop.~1.8]{SS08}), which is stronger than the assertions of Exercise~\ref{E:netcross}.

\bp[Containment by left-most and right-most paths]\label{P:contain}
Let $(\Wl,\Wr)$ be the standard left-right Brownian web, and let $\Ni$ be the Brownian net associated with it.
Then almost surely, there exist no $\pi\in\Ni$ and $l\in\Wl$ such that $l(s)\leq\pi(s)$ and $\pi(t)<l(t)$
for some $\start_\pi\vee\start_l<s<t$. An analogue statement holds for right-most paths.
\ep
The proofs of Proposition~\ref{P:contain} and Theorem~\ref{T:meshchar} are fairly involved and we refer the details to \cite[Thm.~1.7 \& Prop.~1.8]{SS08}.

The wedge and mesh characterizations of the Brownian net have the following interesting corollary.
\bp[Brownian net is closed under hopping]\label{P:nethop}
Let $\Ni$ be the standard Brownian net. Then:
\begin{itemize}
\item[\rm (i)] Almost surely for any $\pi_1, \pi_2\in\Ni$ with $\pi_1(t)=\pi_2(t)$ and $\sigma_{\pi_1}, \sigma_{\pi_2}<t$ for some $t$, the path $\pi$ defined by $\pi:=\pi_1$ on $[\sigma_{\pi_1}, t]$ and $\pi:=\pi_2$ on $[t,\infty)$ is in $\Ni$.

\item[\rm (ii)] For any deterministic $t$, almost surely for any $\pi_1, \pi_2\in \Ni$ with $\pi_1(t)=\pi_2(t)$ and $\sigma_{\pi_1}, \sigma_{\pi_2}\leq t$, the path $\pi$ defined by $\pi:=\pi_1$ on $[\sigma_{\pi_1}, t]$ and $\pi:=\pi_2$ on $[t,\infty)$ is in $\Ni$.
\end{itemize}
\ep
\bex
Use the mesh characterization to prove Proposition~\ref{P:nethop}~(i), and the use the wedge characterization together with Proposition~\ref{P:netdensity} below to prove Proposition~\ref{P:nethop}~(ii).
\eex

Here is an even more striking corollary of the mesh characterization \cite[Prop.~1.13]{SS08}.
\bp[Image set property]\label{P:imgset}
Let $\Ni$ be the Brownian net. For $T\in [-\infty, \infty)$, let $\Ni_T:=\{\pi \in \Ni: \sigma_\pi=T\}$ and $N_T:=\{ (\pi(t), t)\in \R^2: \pi\in\Ni_T, t\in (T, \infty)\}$. Then almost surely for any $T\in [-\infty, \infty)$, any path $\pi\in \Pi$ with $\sigma_\pi=T$ and $\{(\pi(t), t) : t> T\} \subset N_T$ is a path in $\Ni$.
\ep
Note that for the Brownian web, this property is easily seen to hold.

\subsection{The branching-coalescing point set}\label{S:bcptset}

Similar to the definition of the coalescing point set from the Brownian web in Section~\ref{S:coptset}, we can define the so-called {\em branching-coalescing point set} from the Brownian net $\Ni$ as  follows.

Given the Brownian net $\Ni$, and a closed set $A\subset \R$, define the {\em brancing-coalescing point set} by
\begin{equation}\label{bcptset}
\xi^A_t := \{ y\in \R: y=\pi(t) \mbox{ for some } \pi \in \Ni(A\times \{0\}) \}, \qquad t\geq 0.
\end{equation}
In words, $\xi^A_t$ is the set of points in $\R$ that lie on some path in $\Ni$ that start from $A$ at time $0$.

Using the wedge characterisation of the Brownian net, we can compute the density of $\xi^\R_t$. As $t\downarrow 0$, we see in Proposition~\ref{P:netdensity} below that the density diverges at the same rate $1/\sqrt{\pi t}$ as for the coalescing point set in
Proposition~\ref{P:webdensity}, which indicates that coalescence plays the dominant role for small times, while the density converges
to the constant 2 as $t\uparrow \infty$, which results from the balance between branching and coalescence for large times.

\bp[Density of branching-coalescing point set]\label{P:netdensity}
Let $\xi^\R_\cdot$ be the branching-coalescing point set defined from the Brownian net
$\Ni$ as in \eqref{bcptset}. Then for all $t>0$ and $a<b$,
\begin{equation}\label{netdensity}
\E[ |\xi^\R_t \cap [a,b]| ] = (b-a)\Big( \frac{e^{-t}}{\sqrt{\pi t}} + 2\Phi(\sqrt{2t})\Big),
\end{equation}
where $\Phi(x) = \frac{1}{\sqrt{2\pi}} \int_{-\infty}^x e^{-\frac{y^2}{2}} {\rm d}y$.
\ep

\bex Prove Proposition~\ref{P:netdensity} by adapting the proof of Proposition~\ref{P:webdensity} and showing that, almost surely, $\xi^\R_t\cap (a, b)\neq\emptyset$ if and only if $\hat r\in \widehat\Wi^{\rm r}$ starting at $(a, t)$, and $\hat l\in\widehat\Wi^{\rm l}$ starting at $(b,t)$, do not meet above time 0.
\eex

\br\label{R:backbone}
Surprisingly, we can even identify the law of $\xi^\R_t$ as $t\uparrow \infty$, which is a Poisson point process on $\R$ with intensity 2. Furthermore, $\xi_\cdot$ is reversible with respect to the law of the Poisson point process on $\R$ with intensity 2. Formulated in terms of the Brownian net, this amounts to the statement that $\Ni(*,-\infty)$, the collection of Brownian net paths started at time $-\infty$ and is called the {\em backbone of the Brownian net} in \cite{SS08}, has the same distribution as $-\Ni(*, -\infty)$, the set of paths obtained by reflecting the graph of each path in $\Ni(*, -\infty)$ across the origin in space-time. These results were established in \cite[Section 9]{SS08}by first observing their analogues for a discrete system of branching-coalescing random walks, and then passing to the continuum limit.
\er

\br Proposition~\ref{P:netdensity} implies that for each $t>0$, almost surely $\xi^\R_t$ is a locally finite point set. However, almost surely there exist a dense set of random times at which $\xi^\R_t$ contains no isolated point, and is in particular uncountable (see~\cite[Prop.~3.14]{SSS09}).
\er

\br As noted in Remark \ref{R:pfaffian}, started from the whole real line, the coalescing point set forms a Pfaffian point process at each time $t>0$. It will be interesting to investigate whether the branching-coalescing point set $\xi^\R_t$ also admits an explicit characterisation as a Pfaffian point process.
\er

\subsection{Special points of the Brownian net}\label{S:netpts}

Similar to the classification of special points for the Brownian web formulated in Theorem~\ref{T:classweb}, we can give an almost sure classification of all points in $\R^2$ according to the configuration of paths in the Brownian net entering and leaving the point. Such an analysis was carried out in \cite{SSS09}, where it was shown that a.s. there are 20 types of points, in contrast to the 7 types for the Brownian web.

Since Brownian net paths must be contained between paths in the left Brownian web and right Brownian web as stated in Proposition~\ref{P:contain}, the classification of special points is in fact carried out mainly for $(\Wl, \Wr)$. First we introduce a notion of equivalence between paths entering and leaving a point, which is weaker than that introduced for the Brownian web in Section~\ref{S:webpts}.

\bd[Equivalence of paths entering and leaving a point]\label{D:patheqv}
We say $\pi_1, \pi_2\in \Pi$ are equivalent paths entering $z = (x, t) \in \R^2$, denoted by $\pi_1\sim^z_{\rm in} \pi_2$, if $\pi_1$ and $\pi_2$ enter $z$ and $\pi_1(t-\eps_n) = \pi_2(t -\eps_n)$ for a sequence $\eps_n\downarrow 0$. We say $\pi_1, \pi_2$ are equivalent paths leaving $z$, denoted by $\pi_1\sim^z_{\rm out} \pi_2$, if $\pi_1$ and $\pi_2$ leave $z$ and $\pi_1(t + \eps_n) = \pi_2(t + \eps_n)$ for a sequence $\eps_n \downarrow 0$.
\ed
When applied to paths in the Brownian web, the above notion of equivalence implies the equivalence introduced in Section~\ref{S:webpts}, which is why we have abused the notation and used the same symbols $\sim^z_{\rm in}$ and $\sim^z_{\rm out}$. Although $\sim^z_{\rm in}$ and $\sim^z_{\rm out}$ are not equivalence relations on the space of all paths $\Pi$, they are easily seen to be equivalence relations on $\Wl\cup\Wr$. This allows us to classify $z$ according to the equivalence classes of paths entering and leaving $z$, which are necessarily ordered.

To denote the type of a point, we first list the incoming equivalence classes of paths from left to right, and then, separated by a comma, the outgoing equivalence classes of paths from left to right. If an equivalence class contains only paths in $\Wl$, resp.\ $\Wr$, we will label it by l, resp.\ r, while if it contains both paths in $\Wl$ and in $\Wr$, we will label it by p, standing for pair. For points with (up to equivalence) one incoming and two outgoing paths, a subscript l, resp.\ r, means that all incoming paths belong to the left one, resp.\ right one, of the two outgoing equivalence classes; a subscript s indicates that incoming paths in $\Wl$ belong to the left outgoing equivalence class, while incoming paths in $\Wr$ belong to the right outgoing equivalence class. If at a point there are no incoming paths in $\Wl\cup\Wr$, then we denote this by o or n, where o indicates that there are no incoming paths in the Brownian net $\Ni$ , while n indicates that there are incoming paths in $\Ni$ (but none in $\Wl\cup\Wr$).

For example, a point is of type $\rm (p, lp)_r$ if at this point there is one equivalence class of incoming paths in $\Wl\cup\Wr$ and there are two outgoing equivalence classes. The incoming equivalence class is of type p while the outgoing equivalence classes are of type l and p, from left to right. All incoming paths in $\Wl\cup\Wr$ continue as paths in the outgoing equivalence class of type p.

Since the dual left-right Brownian web $(\widehat\Wl, \widehat\Wr)$ can be used to define a dual Brownian net $\widehat\Ni$, the type of $z$ w.r.t.\ $(\widehat\Wi^{\rm l}, \widehat\Wi^{\rm r})$ and $\widehat\Ni$ can be defined in the same way, after rotating their graphs in $\R^2$ by $180^o$ around the origin.

\begin{figure}[tp] 
\centering
\plaat{5cm}{\includegraphics[width=15.5cm]{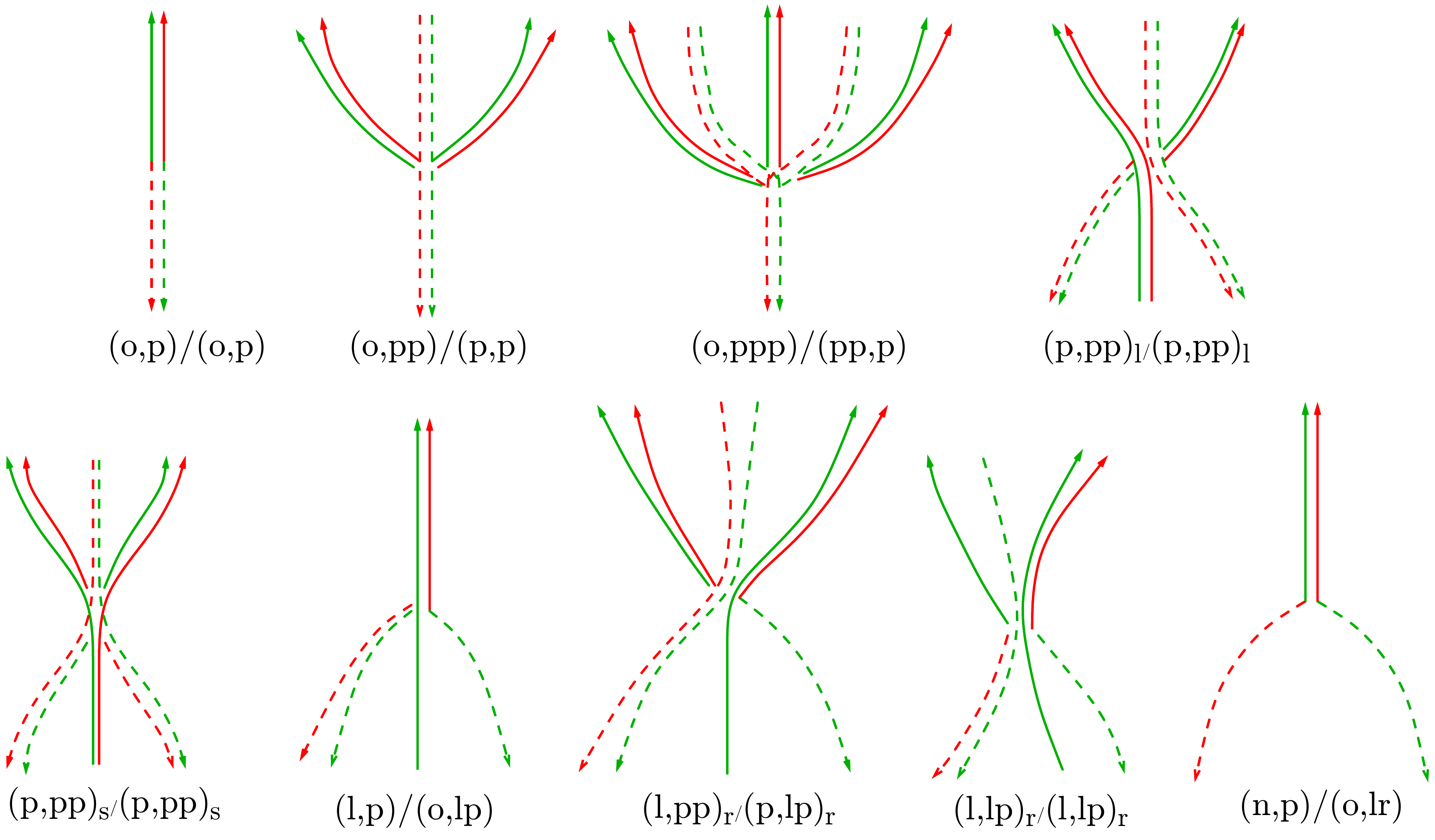}}
\caption{Special points of the Brownian net, modulo symmetry.}
\label{fig:netpoints}
\end{figure}

We now state the classification result from \cite[Theorem 1.7]{SSS09}, while omitting its proof.
\bt[Special points of the Brownian net]\label{T:netpts}
Let $(\Wl, \Wr)$ be the standard left-right Brownian web, let $(\widehat\Wi^{\rm l}, \widehat\Wi^{\rm r})$ be its dual, and let $\Ni$ and $\widehat\Ni$ be the associated Brownian net and its dual. Then almost surely, each point in $\R^2$ is of one of the following 20 types in $\Ni/\widehat\Ni$:
\begin{itemize}
\item[$(1)$]  $\rm (o,p)/(o,p)$, $\rm (o,pp)/(p,p)$, $\rm (p,p)/(o,pp)$, $\rm (o,ppp)/(pp,p)$, $\rm (pp,p)/(o,ppp)$, $\rm (p,pp)_l/(p,pp)_l$, $\rm (p,pp)_r/(p,pp)_r$;

\item[$(2)$] $\rm (p,pp)_s/(p,pp)_s$, {\rm called} separation points;

\item[$(3)$] $\rm (l,p)/(o,lp), (o,lp)/(l,p), (r,p)/(o,pr), (o,pr)/(r,p)$;

\item[$(4)$] $\rm (l,pp)_r/(p,lp)_r, (p,lp)_r/(l,pp)_r, (r,pp)_l/(p,pr)_l, (p,pr)_l/(r,pp)_l$;

\item[$(5)$] $\rm (l,lp)r/(l,lp)_r, (r,pr)l/(r,pr)_l$;

\item[$(6)$] $\rm (n,p)/(o,lr), (o,lr)/(n,p)$;

\end{itemize}
and all of these types occur. For each deterministic time $t\in\R$, almost surely, each point in $\R\times \{t\}$ is of either type $\rm (o, p)/(o, p), (o, pp)/(p, p)$, or $\rm (p, p)/(o, pp)$, and all of these types occur. A deterministic point $(x, t) \in\R^2$ is almost surely of type $\rm (o,p)/(o,p)$.
\et

\br
Note that the points listed in item (1) are analogues of the seven types of points of the Brownian web in Figure~\ref{fig:webspecialpts},
where a path of the Brownian web is replaced by a pair of paths in $(\Wl, \Wr)$. Modulo symmetry, this gives rise to 4 distinct types of points. The types of points listed within each item from (2)--(6) are related to each other by symmetry. Therefore modulo symmetry, there are 9 types of special points for the Brownian net, as illustrated in Figure~\ref{fig:netpoints}.
\er

\br
A basic ingredient in the proof of Theorem~\ref{T:netpts} is the characterisation of the interaction between paths in $(\Wl, \Wr)$ and paths in $(\widehat\Wi^{\rm l}, \widehat\Wi^{\rm r})$. The interaction between paths in $\Wl$ and $\widehat\Wi^{\rm l}$, and similar between paths in $\Wr$ and $\widehat\Wi^{\rm r}$, is given by Skorohod reflection, as mentioned in Remark~\ref{R:skoref}. It turns out that paths in $\Wr$  interact with paths in $\widehat\Wi^{\rm l}$ also via Skorohod reflection, except that given $\hat l\in \widehat\Wi^{\rm l}$, if $r\in \Wr$
initially starts on the left of $\hat l$, then it is Skorohod reflected to the left of $\hat l$ until the collision local time between $r$ and $\hat l$ exceeds an independent exponential random variable, at which time $r$ crosses over to the right of $\hat l$ and is Skorohod reflected to the right of $\hat l$ from that time on, see \cite[Lemma 2.1]{SSS09}.
\er

\br\label{R:structure}
Although Theorem~\ref{T:netpts} classifies points in $\R^2$ mostly according to the configuration of paths in the left-right Brownian web, it can be used to deduce configuration of paths in the Brownian net entering and leaving each type of point. It turns out that $\sim^z_{\rm in}$ is in fact an equivalence relation among paths in $\Ni$, and the same holds for $\sim^z_{\rm out}$ for $z$ of any type other than (o, lr). In particular, when $z$ is not of type (o,lr), each equivalence class of paths in $\Wl\cup\Wr$ entering or leaving $z$ can be enlarged to an equivalence class of paths in $\Ni$. In particular, by Proposition~\ref{P:contain}, all Brownian net paths in an equivalence class of type p, which contains a pair of equivalent paths $(l,r)\in (\Wl, \Wr)$, must be bounded between $l$ and $r$ when sufficiently close to $z$. This applies in particular to points of type $\rm (p,pp)_s$, the {\em separation points}, as well as points of type (pp,p), which we call the {\em meeting points}. See \cite[Section 1.4]{SSS09} for further details.
\er

\bex\label{E:sepcount} Show that the set of separation points of $\Ni$ is a.s.\ countable and dense in $\R^2$.
\eex

\section{Coupling the Brownian web and net}\label{S:couple}

We now introduce a coupling between the Brownian web $\Wi$ and the Brownian net $\Ni$, which is again best motivated from their discrete analogues, the coalescing and branching-coalescing random walks. In particular, this gives the fourth construction of
the Brownian net, the {\em marking construction} mentioned before Section~\ref{S:lrweb} , which was developed by Newman, Ravishankar and Schertzer in \cite{NRS10}.

\begin{figure}[tp] 
\centering
\plaat{5cm}{\includegraphics[width=15.5cm]{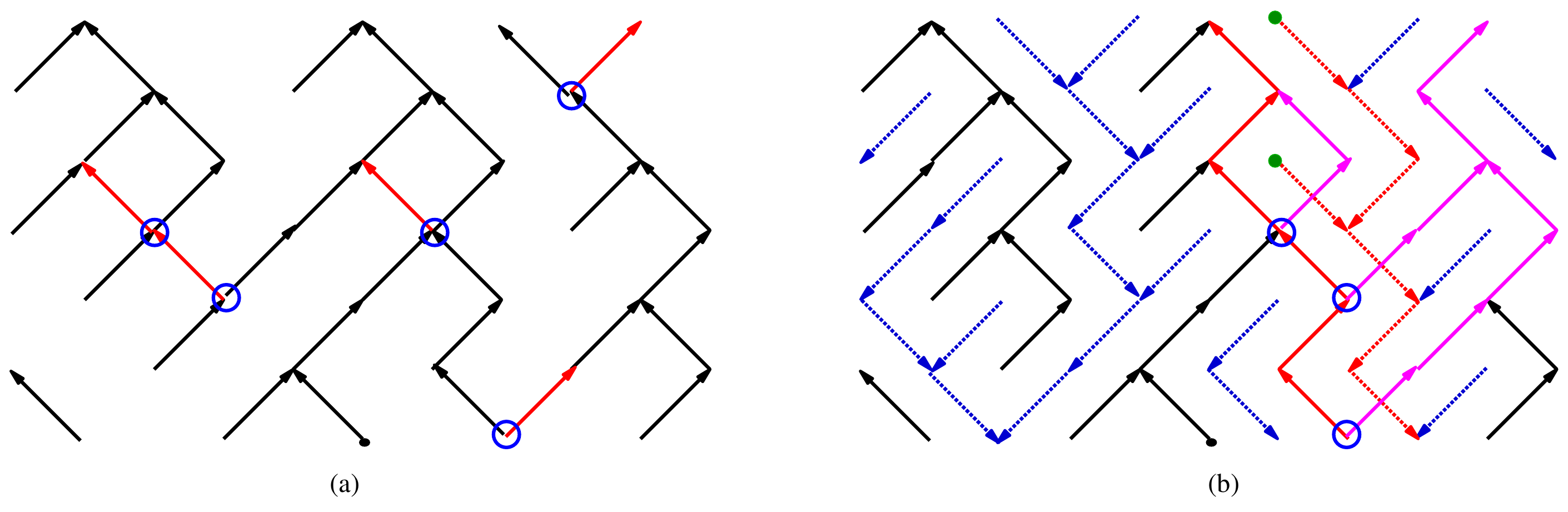}}
\caption{(a) Sampling coalescing random walks from branching-coalescing random walks. (b) Turning coalescing random walks into branching-coalescing random walks. }
\label{fig:couple}
\end{figure}

Given a realisation of branching-coalescing random walks as illustrated in Figure \ref{fig:couple}~(a), where each lattice point in $\Z^2_{\rm even}$ is a branching point with probability $\eps$, we can construct a collection of coalescing random walks by simply forcing the random walk to go either left or right with probability $1/2$ each, independently at each branching point. Interestingly, these branching points have analogues in the continuum limit $\Ni$, which are the {\em separation points}, i.e., points of type $\rm (p,pp)_s$ in Theorem~\ref{T:netpts}.  Therefore given a realisation of the Brownian net $\Ni$, we can sample a Brownian web $\Wi$ by forcing the Brownian web paths to go either left or right with probability $1/2$ each, independently at each separation point. The complication is that each path in the Brownian net will encounter infinitely many separation points on any finite time interval. But fortunately, the number of separation points that is relevant for determining the path's position at a given time $t$ is almost surely locally finite.

Conversely, given a realisation of coalescing random walks as illustrated in Figure~\ref{fig:couple}~(b), we can construct a collection of branching-coalescing random walks by turning each lattice point in $\Z^2_{\rm even}$ independently into a branching point with probability $\eps$. The key observation by Newman et al in~\cite{NRS10} is that, each lattice point in $\Z^2_{\rm even}$ is a point where a random walk path and some dual random walk path are 1 unit of distance apart, and turning that lattice point into a branching point adds a new random walk path that crosses the dual random walk path (see Figure~\ref{fig:couple}~(b), where the dotted lines are the dual paths, and the branching points are circled). In the diffusive scaling limit, most of the added branching points have no effect because they lead to small excursions from the coalescing random walk paths that vanish in the limit. The branching points that have an effect in the limit are points at which a forward random walk path and a dual random walk path, started macroscopically apart, come within 1 unit of distance from each other. In the scaling limit, these become precisely the (1,2) points of the Brownian web in Theorem~\ref{T:classweb}, where a Brownian web path meets a dual Brownian web path. In the discrete setting, these points of close encounter between forward and backward random walk paths are turned independently into branching points with probability $\eps$. In the continuum limit, this leads to Poisson marking of the points of collision between Brownian web paths and dual Brownian web paths, with the intensity measure given by the intersection local time measure between the forward and dual paths. The Brownian net can then be constructed by allowing the Brownian web paths to branch at these Poisson marked (1,2) points, which are precisely the separation points in the resulting Brownian net.

To formulate precisely the coupling between the Brownian web and net motivated by the above heuristic discussions, we first
introduce the necessary background. By Exercise~\ref{E:sepcount}, the set of separation points $S$ is a.s.\ a countable set, and it was pointed out in Remark~\ref{R:structure} that any path $\pi\in \Ni$ entering a separation point $z=(x,t)$ must do so bounded between a pair of equivalent paths $(l, r)\in (\Wl, \Wr)$ entering $z$, and similarly when leaving $z$, it must also be enclosed by one of the two pairs of outgoing equivalent paths $(l_1, r_1), (l_2, r_2)\in (\Wl(z), \Wr(z))$, ordered from left to right. We can then define
\be\label{sgnpi}
{\rm sgn}_\pi(z) := \left\{
\begin{aligned}
-1 \qquad & \mbox{if } l_1\leq \pi\leq r_1 \mbox{  on } [t, \infty), \\
+1 \qquad & \mbox{if } l_2\leq \pi\leq r_2 \mbox{ on } [t, \infty).
\end{aligned}
\right.
\ee
For points of type $(1,2)$ in the Brownian web $\Wi$, we can similarly define
\be\label{sgnw}
{\rm sgn}_\Wi(z) := \left\{
\begin{aligned}
-1 \qquad & \mbox{if z is of type } (1,2)_{\rm l} \mbox{ in } \Wi, \\
+1 \qquad & \mbox{if z is of type } (1,2)_{\rm r} \mbox{ in } \Wi.
\end{aligned}
\right.
\ee
By setting the sign of paths entering each separation point independently to be $\pm1$ with probability $1/2$, we will recover the Brownian web as a subset of the Brownian net.

To construct the Brownian net from the Brownian web, a key object is the local time measure on points of intersection between paths in
$\Wi$ and paths in $\widehat\Wi$.,  i.e., points of type $(1,2)$ in $\Wi$. Its existence was proved in~\cite[Prop.~3.1]{NRS10}, which we quote below (see also \cite[Prop.~3.4]{SSS14}).

\bp{\bf(Intersection local time)}\label{P:refloc}
Let $(\Wi,\hat\Wi)$ be the Brownian web and its dual. Then a.s.\ there exists a
unique measure $\ell$, concentrated on the set of points of type $(1,2)$ in
$\Wi$, such that for each $\pi\in\Wi$ and $\hat\pi\in\hat\Wi$,
\be\label{colmeas}
\begin{aligned}
& \dis\ell\big(\big\{z=(x,t)\in\R^2:
\sig_\pi<t<\hat\sig_{\hat\pi},\ \pi(t)=x=\hat\pi(t)\big\}\big)  \\
=\ & \dis\lim_{\eps\down 0}\eps^{-1}\big|\big\{t\in\R:
\sig_\pi<t<\hat\sig_{\hat\pi},\ |\pi(t)-\hat\pi(t)|\leq\eps\big\}\big|.
\end{aligned}
\ee
The measure $\ell$ is a.s.\ non-atomic and $\si$-finite. We let $\ell_{\rm
  l}$ and $\ell_{\rm r}$ denote the restrictions of $\ell$ to the sets of
points of type $(1,2)_{\rm l}$ and $(1,2)_{\rm r}$, respectively.
\ep
By Poisson marking (1,2) points of $\Wi$ with intensity measure $\ell=\ell_{\rm l}+\ell_{\rm r}$, we obtain a countable set of (1,2) points.
Allowing paths in $\Wi$ to branch at these points then lead to the Brownian net.

We can now formulate the coupling between the Brownian web and net (see \cite[Theorems 4.4 \& 4.6]{SSS14}).
\bt[Coupling between the Brownian web and net]\label{T:coupling}
Let $\Wi$ be the standard Brownian web, and $\Ni$ the standard Brownian net. Let $S$ denote the set of separation points
of $\Ni$. Then there exists a coupling between $\Wi$ and $\Ni$ such that:
\begin{itemize}
\item[{\rm(i)}] Almost surely $\Wi\sub\Ni$, and each separation point $z\in S$ in $\Ni$ is of type $(1,2)$ in  $\Wi$.

\item[{\rm(ii)}] Conditional on $\Ni$, the random variables  $({\rm sgn}_\Wi(z))_{z\in S}$ are i.i.d.\ with $\P[{\rm sgn}_\Wi(z)=\pm1\,|\,\Ni]=1/2$,
and a.s.,
\be\label{webinnet}
\Wi=\{\pi\in\Ni: {\rm sgn}_\pi(z)={\rm sgn}_\Wi(z)\ \forall z\in S
\mbox{ s.t.\ $\pi$ enters }z\}.
\ee

\item[{\rm(iii)}] Conditional on $\Wi$, the sets $S_{\rm l}:=\{z\in S: {\rm sgn}_\Wi(z)=-1\}$ and $S_{\rm r}:=\{z\in S: {\rm sgn}_\Wi(z)=+1\}$ are  independent Poisson point sets with intensities $\ell_{\rm l}$  and $\ell_{\rm r}$, respectively, and a.s.,
\be\label{webtonet}
\Ni = \lim_{\Delta_n\uparrow S} {\rm hop}_{\Delta_n}(\Wi)
\ee
for any sequence of finite sets $\Delta_n$ increasing to $S=S_{\rm l}\cup S_{\rm r}$, where ${\rm hop}_{\Delta_n}(\Wi)$ is
the set of paths obtained from $\Wi$ by allowing paths entering any (1,2) point $z\in \Delta_n$ to continue along either of the two outgoing paths at $z$.
\end{itemize}
\et
\br Theorem~\ref{T:coupling} can be generalised to the case where $\Wl$, $\Wi$, and $\Wr$ may be tilted with different drifts, as long as the drifts remain ordered (see \cite[Theorem 6.15]{SSS14}). We only need to modify the probbility $\P({\rm sgn}_\Wi(z)=\pm 1|\Ni)$ in (ii), while the intensity measures $\ell_{\rm l}$ and $\ell_{\rm r}$ in (iii) need to be multiplied by constants depending on the drifts of $\Wl$,  $\Wi$ and $\Wr$.
\er

\br
Theorem~\ref{T:coupling}~(ii) shows how to sample a Brownian web $\Wi$ from a Brownian net $\Ni$ by forcing the paths to continue either left or right independently at each separation point. What if we sample independently another Brownian web $\Wi'$, conditioned on $\Ni$? How to characterise the joint distribution of $\Wi$ and $\Wi'$? It turns out that $(\Wi, \Wi')$ forms a pair of so-called
{\em sticky Brownian webs}, where Brownian motions in $\Wi$ and $\Wi'$ undergo sticky interaction. Such an object was first introduced in
\cite{HW09b}. For further details,  see~\cite[Sec~3.3 \& Lemma 6.16]{SSS14}.
\er

\noindent
{\bf Proof Sketch.} Theorem~\ref{T:coupling} is proved in \cite[Sec.~6]{SSS14} via discrete approximation, using the fact that a similar coupling as in Theorem~\ref{T:coupling}~(ii) and (iii) holds in the discrete system. The proof is too complex and lengthy to be included here. Instead, we outline below some key ingredients.

The fact that a coupling exists between the Brownian web and net follows by taking the scaling limit of the coupled coalescing and branching-coalescing random walks, where i.i.d.\ signs are assigned to each branching point to determine the path of the coalescing walks.

To show that such a coupling satisfies Theorem~\ref{T:coupling}~(ii), the key is to show that the branching-coalescing random walks, together with the branching points, converge to the Brownian net together with its separation points. Furthermore, under such a convergence, we can match the branching points with the separation points and assign them the same i.i.d.\ signs, such that coalescing random walk paths  converge to paths in the Brownian net as long as the discrete and continuum paths follow the same signs at the respective branching and separation points.

The main difficulty with the above approach is that the set of separation points of $\Ni$ is dense in $\R^2$, so it is unclear in what sense should the branching points converge to the separation points. The solution rests on the observation that, for Brownian net paths starting at some time $S$, their positions at a later time $U$ depend a.s.\ only on the signs of these paths at a locally finite set of separation points  in the time interval $(S,U)$ (called {\em $(S,U)$-relevant separation points}). As we refine our knowledge of the paths at more times in the interval $(S,U)$, more separation points become relevant. Furthermore, the relevant separation points form a locally finite directed graph, called the {\em finite graph representation} of the Brownian net,  which determines a coarse-grained structure of the Brownian net. Therefore a natural notion of convergence is to show that for any $S<U$, the $(S,U)$-relevant separation points and the associated finite graph representation arise as limits of similar structures in the branching-coalescing random walks. Once such a convergence is verified, Theorem~\ref{T:coupling} can then be easily verified.

To show that the coupling between the Brownian web and net satisfies the first statement in Theorem~\ref{T:coupling}~(iii), the key is to show that when restricted to the set of intersection points between a pair of forward and dual Brownian web paths $(\pi, \hat\pi)$, the set of separation points is distributed as a Poisson point process with intensity measure given by the intersection local time measure between $\pi$ and $\hat\pi$. These separation points are in fact relevant separation points with respect to the starting times of $\pi$ and $\hat\pi$, and hence must arise as the limit of relevant branching points in the discrete system. On the other hand, the corresponding relevant branching points form a Bernoulli point process with the intensity measure given by the counting measure on the set of space-time points where a pair of forward and dual coalescing random walk paths approximating $(\pi, \hat\pi)$ come within distance 1. With space-time scaled by $(\eps, \eps^2)$, this counting measure can be shown to converge to the intersection local time measure between $\pi$ and $\hat\pi$, and hence the limit of the branching points (i.e., the set of separation points restricted to the intersections between $\pi$ and $\hat\pi$) must be a Poisson point process with intensity measure given by the intersection local time measure.

To prove the second part of Theorem~\ref{T:coupling}~(iii), i.e., \eqref{webtonet}, note that under the coupling between the Brownian web and net, the Brownian web is embedded in the Brownian net and the Poisson marked (1,2) points of the web are exactly the separation points of the net. Using the coarse-grained structures of the Brownian net given by the finite graph representations, it can then be shown that turning the separation points into branching points for the Brownian web gives the Brownian~net. \qed
\medskip

The proof sketch above shows that the key ingredient in the proof of Theorem~\ref{T:coupling} is the notion of relevant separation points and the associated finite graph representation of the Brownian net. Since they shed new light on the structure of the Brownian net, we will discuss these notions in detail in the rest of this section.

\subsection{Relevant separation points of the Brownian net}\label{S:relsep}

\begin{figure}
\begin{center}
\includegraphics[height=6cm]{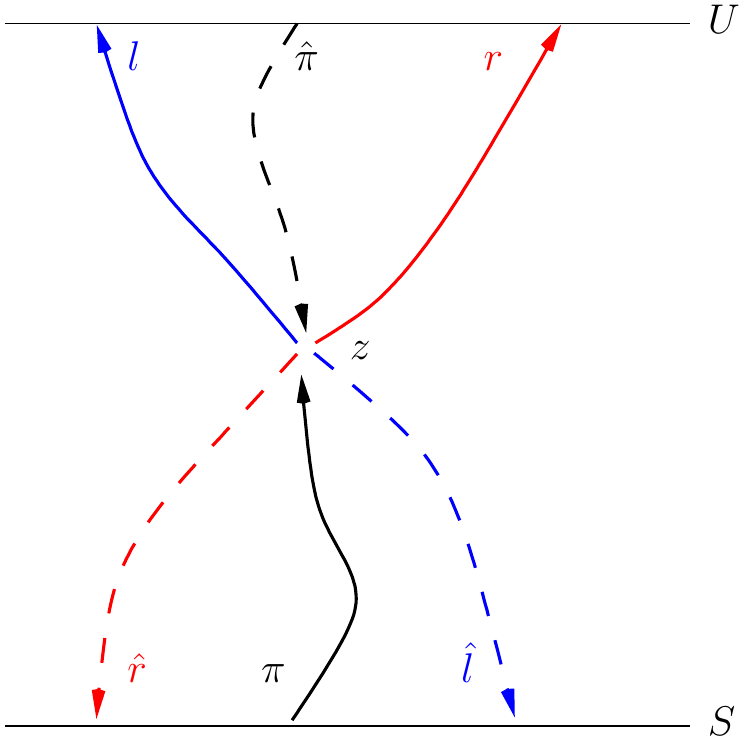}
\caption{An $S,U$-relevant separation point.}\label{fig:relev}
\end{center}
\end{figure}

As noted in Exercise~\ref{E:sepcount}, the set of separation points (i.e., points of type $\rm (p, pp)_s$ in Theorem~\ref{T:netpts}) is
almost surely countable and dense in $\R^2$. In fact almost surely each path in the Brownian net encounters infinitely many separation points in any open time interval. However,  for given deterministic times $S<U$, there is only a locally finite set of separation points which are relevant for deciding where paths in the Brownian net started at time $S$ end up at time $U$. More precisely, as illustrated in Figure~\ref{fig:relev}, we define (cf.~\cite[Sec.~2.3]{SSS09}):

\bd[Relevant separation points]\label{D:relsep}
A separation point $z=(x,t)$ of the Brownian net $\Ni$ is called {\em $(S,U)$-relevant} for some $-\infty\leq S<t<U\leq\infty$, if there exists $\pi\in\Ni$ such that $\sig_\pi=S$ and $\pi(t)=x$, and there exist $l\in\Wl(z)$ and $r\in\Wr(z)$ such that $l<r$ on $(t,U)$.
\ed
Note that because $z$ is a separation point, $l$ and $r$ are continuations of incoming paths at $z$. By Proposition~\ref{P:nethop}, the paths obtained by hopping from $\pi$ to either $l$ or $r$ at $z$ are both paths in the Brownian net, which are distinct on the time interval $(t, U)$. We do not require $l(U)<r(U)$ because a $(S,U)$-relevant separation point defined as above also turns out to be a relevant separation point w.r.t.\ the dual Brownian net $\widehat\Ni$.

\bex
Use the wedge characterization of the Brownian net and the steering argument illustrated in Figure~\ref{fig:funnel} to show that almost surely, for each $-\infty\leq S<U\leq \infty$, a separation point $z=(x,t)$ with $S<t<U$ is $(S,U)$-relevant in $\Ni$ if and only if $-z$ is $(-U,-S)$-relevant in the dual Brownian net $-\widehat\Ni$ rotated $180^o$ around the origin in $\R^2$.
\eex

A crucial property of the $(S,U)$-relevant separation points is that they are almost surely locally finite, which follows from the following density calculation (cf.~\cite[Prop.~2.9]{SSS09}).

\bp{\bf(Density of relevant separation points)}\label{P:relev}
Let $\Ni$ be a standard Brownian net. Then for each deterministic $-\infty\leq S<U\leq\infty$, if $R_{S,U}$
denotes the set of $(S,U)$-relevant separation points, then
\be\label{reldens}
\E\big[\big|R_{S,U}\cap A\big|\big]=2\int_A\Psi(t-S)\Psi(U-t)\,\di x\,\di t \quad \mbox{for all Borel-measurable } A\subset\R\times(S,U),
\ee
where
\be\label{Psib}
\Psi(t):=\frac{e^{-t}}{\sqrt{\pi t}}+2\Phi(\sqrt{2t}), \quad 0<t\leq\infty, \quad \mbox{and}\quad \Phi(x):=\frac{1}{\sqrt{2\pi}}\int_{-\infty}^xe^{-y^2/2}\di y
\ee
In particular, if $-\infty<S$, $U<\infty$, then $R_{S,U}$ is a.s.\ a locally
finite subset of $\R\times[S,U]$.
\ep
{\bf Proof Sketch.}  For $S<s<u<U$, let
$$
E_s:=\{ \pi(s) : \pi \in \Ni, \sigma_\pi=S\} \qquad \mbox{and} \qquad F_u:=\{ \hat\pi(s) : \hat\pi \in \widehat\Ni, \hat\sigma_{\hat\pi}=U\}.
$$
Observe that a.s., each $(S,U)$-relevant separation point $z$ in the strip $\R\times [s,u]$ can be traced back along some Brownian net
path $\pi\in\Ni$ to a position $x\in E_s$. Furthermore, if $(l, r)\in (\Wl, \Wr)$ is the pair of left-right Brownian web paths starting at $(x,s)$, then we must have $F_u\cap (l(u), r(u))\neq \emptyset$, and $z$ is bounded between $l$ and $r$ on the time interval $[s,u]$. Therefore each $(S,U)$-relevant separation point in $\R\times [s,u]$ can be approximated by some $(x,s)$ with $x$ in
$$
Q_{s,u} := \{ x\in E_s:  \ F_u\cap ( l(u), r(u)) \neq \emptyset, \mbox{ where } (l,r) \in (\Wl(x,s), \Wr(x,s))\}.
$$
The density of the set $Q_{s,u}$ can be easily computed. As we partition $(S, U)$ into disjoint intervals $[t_i, t_{i+1})$ of size $1/n$ with $n\uparrow \infty$,  each $(S,U)$-relevant separation point $z$ can be approximated by some $(x_n, t_{i_n})$ with $x_n\in Q_{t_{i_n}, t_{i_n+1}}$ for some $i_n$, such that $(x_n, t_{i_n})\to z$ as $n\to\infty$. The upper bound on the density of $R_{S,U}$ then follows by Fatou's Lemma. For the lower bound, it suffices to consider points in $Q_{t_i, t_{i+1}}$ which are separated from each other by some $\delta>0$ which is sent to zero after taking the limit $n\to\infty$.
\qed

\bex
Prove Proposition~\ref{P:relev} by rigorously implementing the strategy outlined above.
\eex

\subsection{Finite graph representation of the Brownian net}

Having defined $(S,U)$-relevant separation points and established their local finiteness, we are now ready to introduce the {\em finite graph representation}, which is a directed graph with the interior vertices given by the $(S,U)$-relevant separation points, and the directed edges determined by how Brownian net paths go from one relevant separation point to the next. Such a directed graph gives a coarse-grained representation of the Brownian net, and plays a key role in studying the coupling between the Brownian web and net in Theorem~\ref{T:coupling}.

 Let $-\infty<S<U<\infty$ be deterministic times, let $R_{S,U}$ be the set of $(S,U)$-relevant separation points of $\Ni$ and set
\bc
\dis R_S&:=&\dis\R\times\{S\},\\[5pt]
\dis R_U&:=&\dis\big\{(x,U):x\in\R,\ \exists\pi\in\Ni\mbox{ with }
\sig_\pi=S\mbox{ s.t.\ }\pi(U)=x\big\}.
\ec
We will make the set $R:=R_S\cup R_{S,U}\cup R_U$ into a directed graph, with the directed edges representing the Brownian net paths.

First we identify the directed edges leading out from each $z=(x,t)\in R_{S,U}$. By Remark~\ref{R:structure} and Proposition~\ref{P:contain}, any $\pi\in\Ni$ with $\sigma_\pi=S$ that enters a $(S, U)$-relevant separation point $z$ must leave $z$ enclosed by one of the two outgoing pairs of equivalent left-right paths $(l_1, r_1), (l_2, r_2)\in (\Wl(z), \Wr(z))$,  with $l_1\sim^z_{\rm out} r_1$ and $l_2\sim^z_{\rm out}r_2$. Take the pair $(l_1, r_1)$ for instance. If $l_1(U)\neq r_1(U)$, then there must be a last separation point $z'$ along the pair $(l_1, r_1)$, which is also an $(S, U)$-relevant separation point because hopping from $\pi$ to either $l_1$ or $r_1$ at time $t$ still gives a Brownian net path. If $l_1(U)=r_1(U)$, then we just set $z'=(l_1(U), U)=(r_1(U), U) \in R_U$. In either case, we note that all paths that leave $z$ bounded between $l_1$ and $r_1$ must continue to do so until they reach the $z'$. Therefore we draw a directed edge from $z$ to $z'$, denoted by $z\to_{l_1,r_1}z'$, representing all Brownian net paths that go from $z$ to $z'$ while bounded between $l_1$ and $r_1$. Similarly, paths that leave $z$ while bounded between $l_2$ and $r_2$ must all lead to some $z''\in R_{S,U}\cup R_U$, which is represented by a directed edge $z\to_{l_2,r_2}z''$.

Next we identify the directed edges leading out from $z=(x,S)\in R_S$. By Theorem~\ref{T:netpts}, at a deterministic time $S$, almost surely each $z\in R_S$ is of type (o,p), (p,p), or (o, pp). If $z$ is of type (o,p) or (p,p), then a single equivalent pair
$l\sim^z_{\rm out} r$ starts from $z$, and all Brownian net paths leaving $z$ must be bounded between $l$ and $r$ and leads to some
$z'\in R_{S, U}\cup R_U$, which we represent by a directed edge $z\to_{l, r} z'$. Similarly when $z$ is of type (o,pp), two directed edges start from $z$.

Given the above directed graph with vertex set $R_S\cup R_{S, U}\cup R_U$, it is not difficult to see that each Brownian net path from time $S$ to $U$ corresponds to a directed path from $R_S$ to $R_U$, and conversely, each directed path from $R_S$ to $R_U$ can be associated with a family of Brownian net paths from time $S$ to $U$. We summarize the basic properties below (see \cite[Prop.~6.5]{SSS14}).

\begin{figure}
\begin{center}
\includegraphics[width=15.5cm]{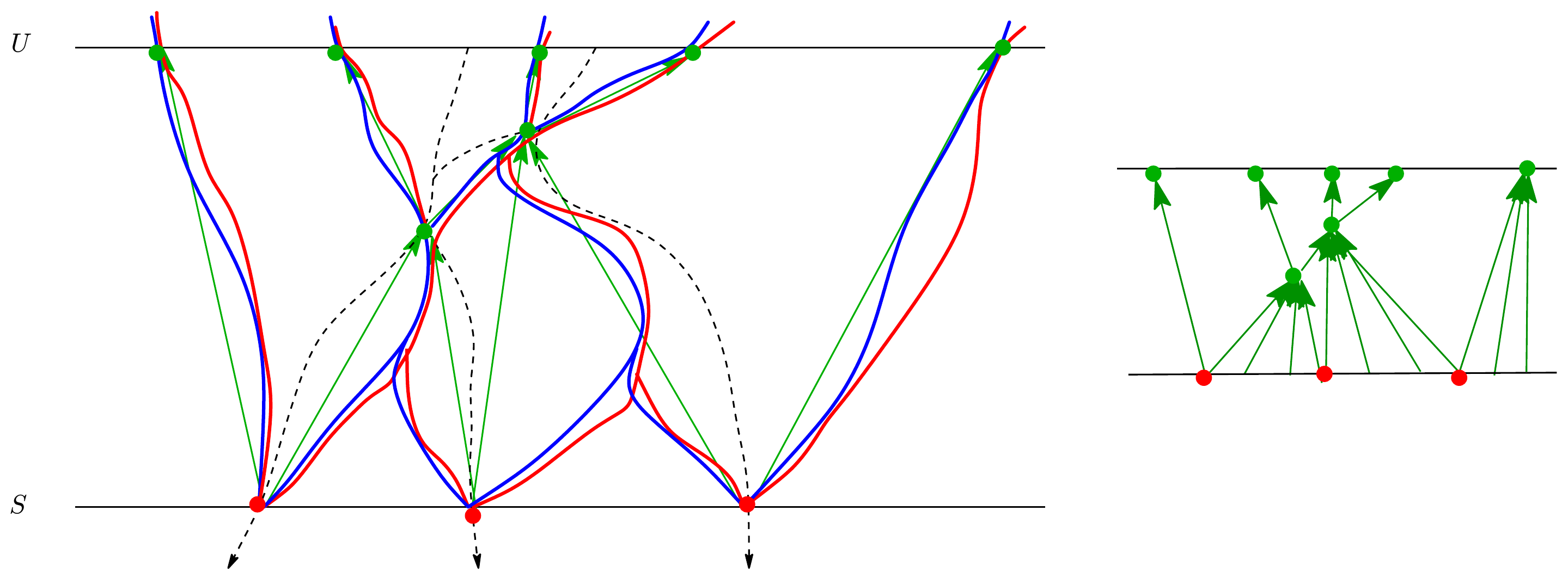}
\caption{Finite graph representation.}\label{fig:fingraph}
\end{center}
\end{figure}

\index{finite graph representation}
\bp{\bf(Finite graph representation)}\label{P:fingraph}
Let $\Ni$ be a Brownian net with associated left-right Brownian web
$(\Wl,\Wr)$, and let $-\infty<S<U<\infty$ be deterministic times.
Let $R:=R_S\cup R_{S,U}\cup R_U$ and directed edges $\to_{l,r}$ be defined as
above. Then, a.s.\ $($see Figure~\ref{fig:fingraph}$)$:
\begin{itemize}
\item[{\bf(a)}] For each $z\in R_S$ that is not of type $\rm(o,pp)$, there
  exist unique $l\in\Wl(z)$, $r\in\Wr(z)$ and $z'\in R$ such that
  $z\to_{l,r}z'$.
\item[{\bf(b)}] For each $z=(x,t)$ such that either $z\in R_{S,U}$ or $z\in
  R_S$ is of type $\rm(o,pp)$, there exist unique $l,l'\in\Wl(z)$,
  $r,r'\in\Wr(z)$ and $z',z''\in R$ such that $l\leq r'<l'\leq r$ on
  $(t,t+\eps)$ for some $\eps>0$, $z\to_{l,r'}z'$ and $z\to_{l',r}z''$. For
  $z\in R_{S,U}$ one has $z'\neq z''$. For $z\in R_S$ of type $\rm(o,pp)$, one
  has $z'\neq z''$ if and only if there exists $\hat\pi\in\widehat\Ni$
  with $\hat\sig_{\hat\pi}=U$ such that $\hat\pi$ enters $z$.
\item[{\bf(c)}] For each $\pi\in\Ni$ with $\sig_\pi=S$, there exist
  $z_i=(x_i,t_i)\in R$ $(i=0,\ldots,n)$ and $l_i\in\Wl(z_i)$, $r_i\in\Wl(z_i)$
  $(i=0,\ldots,n-1)$ such that $z_0\in R_S$, $z_n\in R_U$,
  $z_i\to_{l_i,r_i}z_{i+1}$ and $l_i\leq\pi\leq r_i$ on $[t_i,t_{i+1}]$
  $(i=0,\ldots,n-1)$.
\item[{\bf(d)}] If $z_i=(x_i,t_i)\in R$ $(i=0,\ldots,n)$ and
  $l_i\in\Wl(z_i)$, $r_i\in\Wl(z_i)$ $(i=0,\ldots,n-1)$ satisfy $z_0\in
  R_S$, $z_n\in R_U$, and $z_i\to_{l_i,r_i}z_{i+1}$ $(i=0,\ldots,n-1)$, then
  there exists a $\pi\in\Ni$ with $\sig_\pi=S$ such that $l_i\leq\pi\leq r_i$
  on $[t_i,t_{i+1}]$.
\end{itemize}
\ep

\bex
Use the finite graph representation to show that: conditioned on the Brownian net $\Ni$ with the set of separation points
$S$, almost surely for any $(\alpha_z)_{z\in S}\in \{\pm1\}^S$,
$$
W=\{ \pi \in \Ni: {\rm sgn}_\pi(z) = \alpha_z \, \forall\, z\in S \mbox{ s.t. } \pi \mbox{ enters } z\}
$$
$($cf.~\eqref{webinnet}$)$ is a closed subset of $\Ni$, containing at least one path starting from each $(x,t)\in\R^2$.
\eex

\section{Scaling limits of random walks in i.i.d.\ space-time environment}\label{S:rwre}
We now show how the Brownian web and net can be used to construct the continuum
limits of one-dimensional random walks in i.i.d.\ random space-time environments.

\begin{figure}
\begin{center}
\includegraphics[width=11cm]{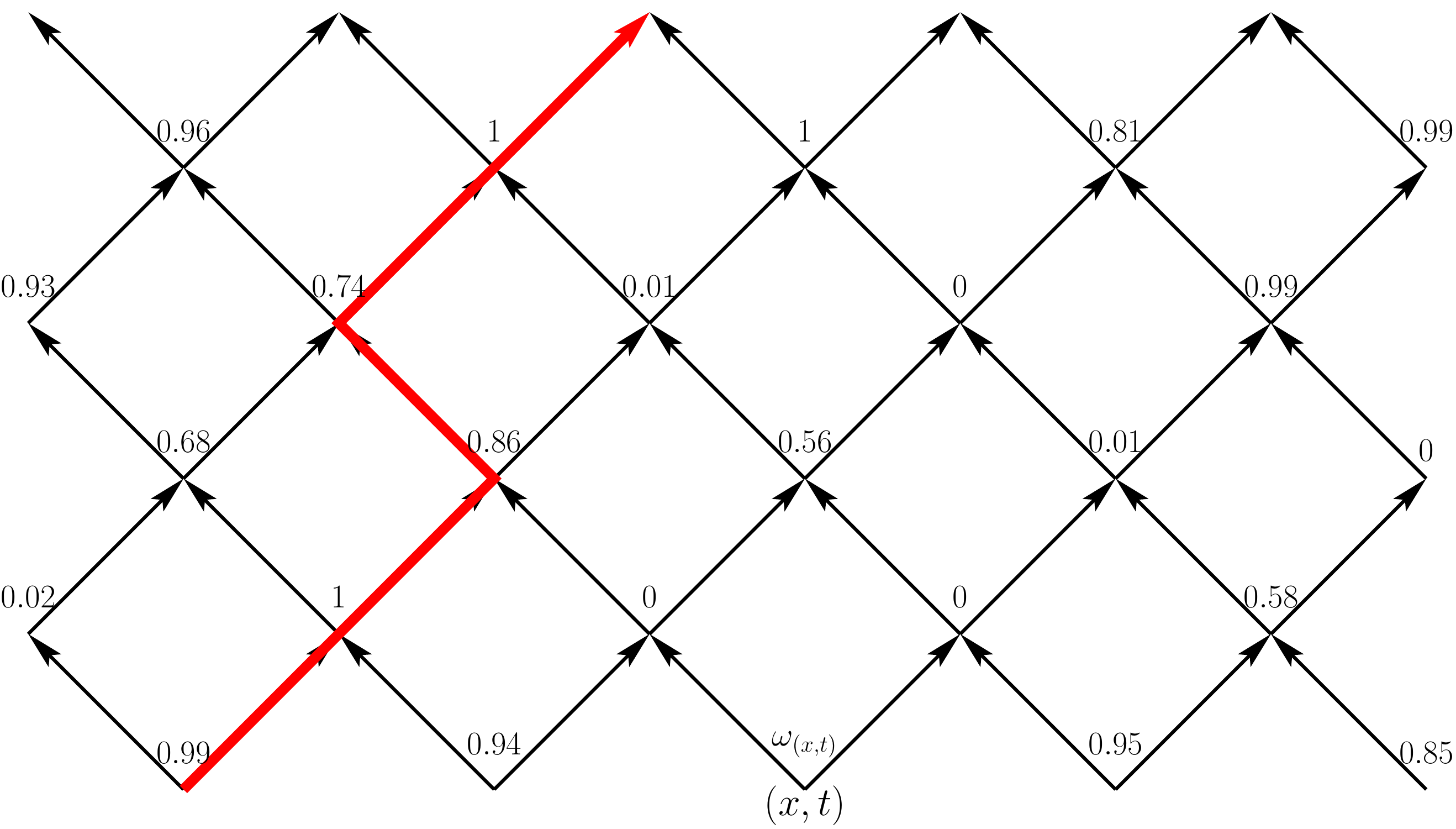}
\caption{Random walk on $\Z^2_{\rm even}$ in a random environment $\omega$.}
\label{fig:disflow}
\end{center}
\end{figure}

A random walk $X$ in an i.i.d.\ random space-time environment is defined as follows. Let $\omega:=(\omega_z)_{z\in \Z^2_{\rm even}}$ be i.i.d.\ $[0,1]$-valued random variables with common distribution $\mu$. We view $\omega$ as a random space-time environment, such that conditioned on $\omega$, if the random walk is at position $x$ at time $t$, then in the next time step the walk jumps to $x+1$ with probability $\omega_{(x,t)}$ and jumps to $x-1$ with probability $1-\omega_{(x,t)}$ (see Figure~\ref{fig:disflow}). Let $P^\omega_{(x,s)}$ denote probability for the walk starting at $x$ at time $s$ in the environment $\omega$, and let $\P$ denote probability for $\omega$.

The question is: if we scale space and time by $\eps$ and $\eps^2$ respectively, is it possible to choose a law $\mu_\eps$ for an environment $\omega^\eps:=(\omega^\eps_z)_{z\in \Z^2_{\rm even}}$ such that the walk in the environment $\omega^\eps$ converges as $\eps\downarrow 0$ to a limiting random motion in a continuum space-time random environment? Different ways of looking at the random walk in random environment suggest that the answer is yes.

One alternative way of characterizing the law of random walk in random environment is to consider the law of the family of random transition probability kernels
\be\label{5.1}
K^\omega_{s,t}(x,y):=P^\omega_{(x,s)}(X(t)=y),  \qquad  (x,s), (y,t) \in \Z^2_{\rm even},  \ s\leq t.
\ee

Still another way is to specify the {\em $n$-point motions}, i.e., the law of $n$ random walks $\vec X:=(X_1(t), \ldots, X_n(t))_{t\geq 0}$ sampled independently from the same environment $\omega$, and then averaged with respect to the law of $\omega$. Note that $\vec X$ is a Markov chain with transition probability kernel
\be\label{5.2}
K^{(n)}_{0,t}(\vec x, \vec y) = \int \prod_{i=1}^n K^\omega_{0,t}(x_i, y_i) \P({\rm d}\omega), \quad (x_i, 0), (y_i, t) \in \Z^2_{\rm even}  \mbox{ for } 1\leq i\leq n.
\ee
Furthermore, the $n$-point motions are {\em consistent} in the sense that the law of any $k$-element subset of $(X_1(\cdot), \ldots, X_n(\cdot))$ is governed by that of the $k$-point motion. Note that  the moments of $\omega_{(0,0)}$ are determined by $(K^{(n)})_{n\in\N}$, since
$$
\E[ \omega_{(0,0)}^n] = K^{(n)}_{0,1}(\vec 0, \vec 1).
$$
Therefore the law of $\omega_{(0,0)}$, as well as that of $\omega$ and $K^\omega$, are uniquely determined by the law of the $n$-point motions $(K^{(n)})_{n\in\N}$.

Evidence for the existence of a continuum limit for random walks in i.i.d.\ space-time random environments came from the convergence of the $n$-point motions, which was first established by Le Jan and Lemaire~\cite{LL04} for i.i.d.\ Beta-distributed random environments, and subsequently extended to general environments by Howitt and Warren~\cite{HW09a}. By the theory of {\em stochastic flow of kernels} developed by Le Jan and Raimond~\cite{LR04}, this implies that the family of random probability kernels $K^\omega$ in \eqref{5.2} also converges in a suitable sense to a continuum limit. Motivated by these results, it was then shown in \cite{SSS14} that under the same convergence criterion as in \cite{HW09a}, not only the $n$-point motions converge, but also the random environments themselves converge to a continuum space-time limit which can be constructed explicitly from the Brownian web and net.

In Section~\ref{S:HW} below, we will first briefly recall the theory of {\em stochastic flows of kernels} from \cite{LR04}, and then review Howitt and Warren's convergence result~\cite{HW09a} for the $n$-point motions of random walks in i.i.d.\ space-time random environments. In Section~\ref{S:HWenv}, we will then show how to construct the limiting continuum space-time environment from the Brownian web and net,  which arises naturally if one looks at the discrete environments in the right way. Lastly in Section~\ref{S:HWprop}, we will discuss some properties of the continuum random motion in random environment.

\subsection{Stochastic flows of kernels and the Howitt-Warren flows}\label{S:HW}

In \eqref{5.1} and \eqref{5.2}, we saw that the family of random transition probability kernels $K^\omega_{s,t}(x,\cdot)$ and the family of consistent $n$-point motions provide alternative characterizations of a random walk in i.i.d.\ space-time random environment. It turns out that in general, without knowing the existence of any underlying random environment, there is a correspondence between a consistent family of $n$-point motions on a Polish space, and a family of random probability kernels called a {\em stochastic flow of kernels}. This was the main result of Le Jan and Raimond in \cite[Theorem 2.1]{LR04}, which motivated the study of concrete examples of consistent $n$-point motions in~\cite{LR04b, LL04, HW09a}. We now recall the notion of a {\em stochastic flow of kernels}.

For any Polish space $E$, let $\Bi(E)$ denote the Borel $\sigma$-field on $E$ and let $\Mi_1(E)$ denote the space of probability measures on $E$, equipped with the weak topology and Borel $\sigma$-algebra. A {\em random probability kernel}, defined on some probability space $(\Omega,\Fi,\P)$, is a measurable function $K:\Omega\times E\to\Mi_1(E)$. Two random probability kernels $K,K'$ are said to equal in finite dimensional distributions if for each $x_1,\ldots,x_n\in E$, the $n$-tuple of random probability measures $\big(K(x_1,\,\cdot\,),\ldots,K(x_n,\,\cdot\,)\big)$ is equally distributed with $\big(K'(x_1,\,\cdot\,),\ldots,K'(x_n,\,\cdot\,)\big)$. Two or more random probability kernels are called independent if their finite-dimensional distributions are independent.

\bd[Stochastic flow of kernels]\label{D:stochflow}
A {\em stochastic flow of kernels} on a Polish space $E$ is a collection $(K_{s,t})_{s\leq t}$ of random probability kernels on $E$ such that
\begin{itemize}

\item[{\rm(i)}] For all $s\leq t\leq u$ and $x\in E$, a.s.\ $K_{s,s}(x,A)=\de_x(A)$ and $\dis\int_EK_{s,t}(x, {\rm d}y)K_{t,u}(y,A)=K_{s,u}(x,A)$ for all $A\in\Bi(E)$.

\item[{\rm(ii)}] For each $t_0<\cdots<t_n$, the random probability kernels $(K_{t_{i-1},t_i})_{i=1,\ldots,n}$ are independent.

\item[{\rm(iii)}] $K_{s,t}$ and $K_{s+u,t+u}$ are equal in finite-dimensional distributions for each real $s\leq t$ and~$u$.\end{itemize}
\ed
We have omitted two weak continuity conditions on $(K_{s,t}(x, \cdot)_{s\leq t, x\in E}$ from \cite[Def.~2.3]{LR04}, which are automatically satisfied by a version of $K$ if the $n$-point motions defined via \eqref{5.2} is a family of Feller processes.

\br\label{R:CK} Although we motivated the notion of stochastic flow of kernels from random walks in i.i.d.\ space-time random environments, the kernels $K$ in Definition~\ref{D:stochflow} cannot be associated with an underlying random environment unless there exists a version of $K$ such that Definition \ref{D:stochflow}~(i) is strengthened to:
$$
{\rm (i')} \quad A.s.,\  \forall\, s\leq t\leq u \mbox{ and } x\in E, \ \ K_{s,s}(x,\cdot)=\de_x(\cdot) \mbox{ and } \dis\int_EK_{s,t}(x, {\rm d}y)K_{t,u}(y,\cdot)=K_{s,u}(x,\cdot).
$$
In general, it is not known whether such a version always exists.
\er

Le Jan and Raimond showed in \cite[Theorem 2.1]{LR04} that every consistent family of Feller $n$-point motions corresponds to a stochastic flow of kernels. Using Dirichlet form construction of Markov processes, they then constructed as an example in~\cite{LR04b} a consistent family of $n$-point motions on the circle which are a special type of sticky Brownian motions. Subsequently, Le Jan and Lemaire showed in~\cite{LL04} that the $n$-point motions of random walks in i.i.d.\ Beta distributed space-time random environments converge to the $n$-point motions constructed in \cite{LR04b}. Howitt and Warren~\cite{HW09a} then found the general condition for the convergence of $n$-point motions of random walks in i.i.d.\ space-time random environments, and they characterized the limiting $n$-point motions, which are also sticky Brownian motions, in terms of well-posed martingale problems.

We next recall Howitt and Warren's result \cite{HW09a}, or rather, a different formulation of their result as presented in \cite[Appendix A]{SSS14} with discrete instead of continuous time random walks, and with a reformulation of the martingale problem characterizing the limiting sticky Brownian motions.

\bt[Convergence of $n$-point motions to sticky Brownian motions]\label{T:HWconv} For $\eps>0$, let $\mu_\eps$ be the common law of an i.i.d.\ space-time random environment
$(\omega_z^\eps)_{z\in \Z^2_{\rm even}}$, such that
\be\label{mucon}
\begin{aligned}
{\rm(i)} \quad & \dis\eps^{-1}\int(2q-1)\mu_\eps(\di q) \ \aston{\eps} \  \bet,\\
{\rm(ii)} \quad & \dis\eps^{-1}q(1-q)\mu_\eps(\di q)\ \Aston{\eps}\  \nu(\di q)
\end{aligned}
\ee
for some $\beta\in\R$ and finite measure $\nu$ on $[0,1]$, where $\Rightarrow$ denotes weak convergence.

Let $\vec X^\eps:=(X^\eps_1, \ldots, X^\eps_n)$ be $n$ independent random walks in $\omega^\eps$ with  $\eps \vec X^\eps(0) \to \vec X(0)= (X_1(0), \ldots, X_n(0))\in \R^n$ as $\eps\downarrow 0$. Then $(\eps \vec X^\eps(\eps^2 t))_{t\geq 0}$ converges weakly to
a family of sticky Brownian motions $(\vec X(t))_{t\geq 0}=(X_1(t), \ldots, X_n(t))_{t\geq 0}$, whose law is the unique solution of
the following {\em Howitt-Warren martingale problem with drift $\beta$ and characteristic measure $\nu$}:
\begin{itemize}
\item[\rm (i)] $\vec X$ is a continuous, square-integrable semimartingale with initial condition $\vec X(0)$;

\item[\rm (ii)] The covariance process between $X_i$ and $X_j$ is given by
\be\label{HWcov}
\langle X_i,X_j\rangle(t)=\int_0^t1_{\{X_i(s)=X_j(s)\}}\di s,  \qquad t\geq 0,\ i,j=1,\ldots,n;
\ee

\item[\rm (iii)] For each non-empty $\Delta \subset \{1, \ldots, n\} $, let
\be\label{fgDe}
f_\Delta(\vec x):=\max_{i\in\Delta}x_i\qquad\mbox{and}\qquad
g_\Delta(\vec x):=\big|\{i\in\Delta:x_i=f_\Delta(\vec x)\}\big| \qquad(\vec x\in\R^n).
\ee
Then
\be\label{MP2b}
f_\De\big(\vec X(t)\big)-\int_0^t\bet_+\big(g_\De(\vec X(s))\big)\di s
\ee
is a martingale with respect to the filtration generated by $\vec X$, where
\be\label{betplus}
\beta_+(1) := \bet\quad\mbox{and} \quad \beta_+(m):= \beta+2\int\nu({\rm d} q)\sum_{k=0}^{m-2}(1-q)^k \quad \mbox{for } m\geq 2.
\ee
\end{itemize}
\et

\br
The Howitt-Warren sticky Brownian motions evolve independently when they are apart, and experience sticky interaction when they meet. In particular, when $n=2$, $X_1(t)-X_2(t)$ is a Brownian motion with stickiness at the origin, which is just a time changed Brownian motion such that its local time at the origin has been turned into real time, modulo a constant multiple that determines the stickiness. More generally, for Howitt-Warren sticky Brownian motions started at $X_1(0)=\cdots=X_n(0)$, the set of times with $X_1(t)=X_2(t)=\cdots=X_n(t)$ is a nowhere dense set with positive Lebesgue measure. The measure $\nu$ determines a two-parameter family of constants $\theta_{k,l}=\int q^{k}(1-q)^{l} \frac{\nu({\rm d}q)}{q(1-q)}$, $k,l\geq 1$, which can be regarded as the rate (in a certain excursion theoretic sense) at which $(X_1,\cdots, X_n)$ split into two groups, $(X_1,\cdots, X_k)$ and $(X_{k+1},\cdots, X_{k+l})$, with $k+l=n$.
\er

It is easily seen that the $n$-point motions defined by the Howitt-Warren martingale problem in Theorem~\ref{T:HWconv} form a consistent family, and it is Feller by \cite[Prop.~8.1]{HW09a}. Therefore by the afore-mentioned result of Le Jan and Raimond \cite[Thm~2.1]{LR04}, there exists a stochastic flow of kernels $(K_{s,t})_{s\leq t}$ on $\R$, unique in finite-dimensional distributions, such that the $n$-point motions of $(K_{s,t})_{s\leq t}$ in the sense of (\ref{5.2}) are given by the unique solutions of the Howitt-Warren martingale problem. Therefore we define

\bd[Howitt-Warren flow]\label{D:HWflow}
We call the stochastic flow of kernels, whose $n$-point motions solve the Howitt-Warren martingale problem in Theorem~\ref{T:HWconv} for some $\bet\in\R$ and finite measure $\nu$ on $[0,1]$, the {\em Howitt-Warren flow} with {\em drift} $\bet$ and {\em characteristic measure} $\nu$.
\ed
\br
We single out three special classes of Howitt-Warren flows: (1) $\nu=0$, for which the $n$-point motions are coalescing Brownian motions, and the flow is known as the {\em Arratia flow}; (2) $\nu({\rm d}x) = a {\rm d}x$ for
some $a>0$, which we will call the {\em Le Jan-Raimond flow} since it was first constructed in~\cite{LR04b} via Dirichlet forms and subsequently shown in~\cite{LL04} to arise as limits of random walks in i.i.d.\ Beta-distributed environments; (3) $\nu = a\delta_0 + b\delta_1$ for some $a, b\geq 0$ with $a+b>0$, called the {\em erosion flow}, which was studied in \cite{HW09b}.
\er

As noted in Remark~\ref{R:CK}, it is not known a priori whether there exists a version of the Howitt-Warren flow $K$ such that almost surely, $(K_{s,t}(x, \cdot))_{s<t, x\in \R}$ are truly transition probability kernels of a random motion in a random environment. It is then natural to ask whether such an underlying random environment indeed exists for the Howitt-Warren flows, and if yes, whether the environment can be explicitly characterized as in the discrete case. This is where the Brownian web and Brownian net enter the picture.

\subsection{The space-time random environment for the Howitt-Warren flows}\label{S:HWenv}

How to construct a continuum space-time random environment such that the Howitt-Warren flow with drift $\beta$ and characteristic measure $\nu$ is indeed the family of transition probability kernels of a random motion in this random environment? The answer again lies in discrete approximation.

{\bf Special case}: $\nu$ satisfies $\int q^{-1}(1-q)^{-1}\nu(\di q)<\infty$. In this case, the continuum random environment can be
constructed from the Brownian net. For $\eps>0$, define a probability measure $\mu_\eps$ on $[0,1]$ by
\be\label{nubar}\ba{l}
\mu_\eps:=b\eps \bar\nu+\ffrac{1}{2}(1-(b+c)\eps)\de_0
+\ffrac{1}{2}(1-(b-c)\eps)\de_1\\[5pt]
\dis\quad\mbox{where}\quad
b:=\int\frac{\nu(\di q)}{q(1-q)},\quad
c:=\bet-\int(2q-1)\frac{\nu(\di q)}{q(1-q)},\quad
\bar\nu(\di q):=\frac{\nu(\di q)}{bq(1-q)}.
\ec
When $\eps$ is sufficiently small, such that $1-(b+|c|)\eps\geq 0$, $\mu_\eps$ is a probability measure on $[0,1]$ and is easily seen to satisfy \eqref{mucon} as $\eps\downarrow 0$. Therefore by Theorem~\ref{T:HWconv}, $\mu_\eps$ determines the law of an i.i.d.\ random environment $\omega^\eps:=(\omega^\eps_z)_{z\in \Z^2_{\rm even}}$, whose associated $n$-point motions converge to that of the Howitt-Warren flow with drift $\beta$ and characteristic measure $\nu$.

\begin{figure}[t]
\begin{subfigure}{.5\textwidth}
  \centering
  \includegraphics[width=.9\linewidth]{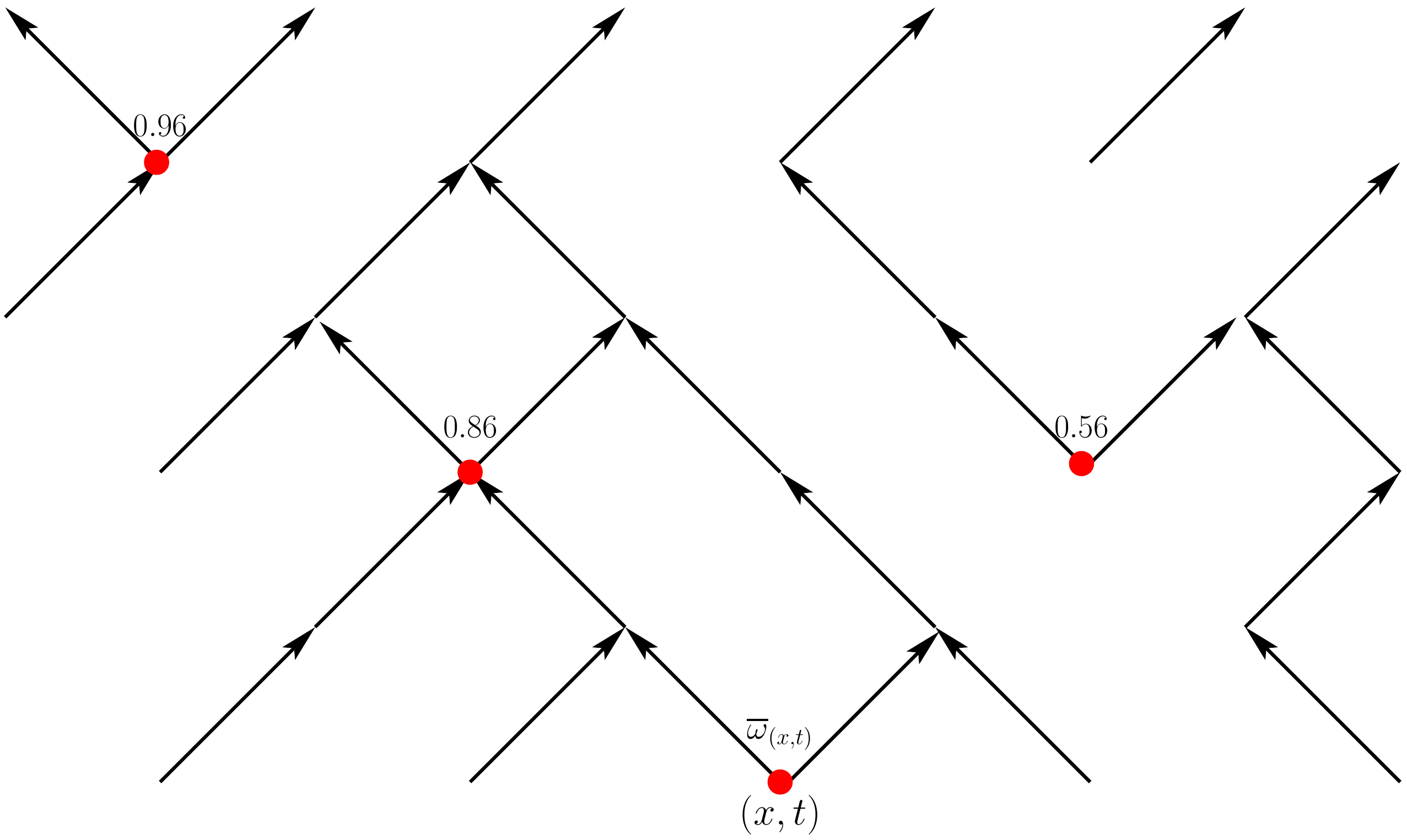}
  \caption{}
  \label{fig:netenviron}
\end{subfigure}%
\begin{subfigure}{.5\textwidth}
  \centering
  \includegraphics[width=.9\linewidth]{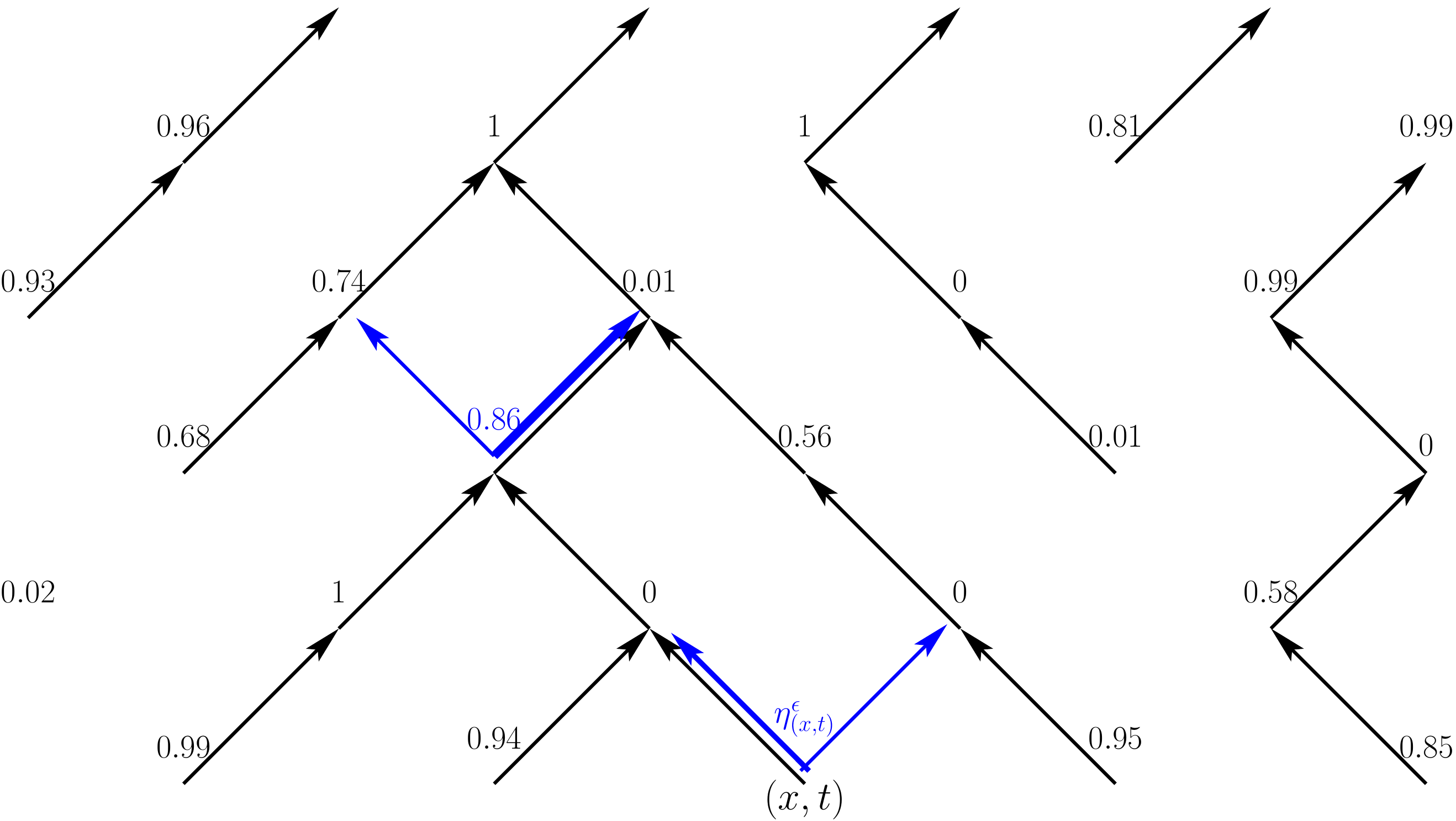}
  \caption{}
  \label{fig:webenviron}
\end{subfigure}
\caption{Representation of the random environment $(\omega^\eps_z)_{z\in\Z^2_{\rm even}}$ in terms of: (a) a marked discrete net $(N^\eps,\bar\omega^\eps)$; (b) a marked discrete web $(W^\eps,\eta^\eps)$.}
\label{fig:netwebenviron}
\end{figure}

Note that for small $\eps$, most of the $\omega^{\eps}_z$ are either zero or one. We can thus encode $\omega^\eps$ in two steps. First we identify a collection of branching-coalescing random walks determined by $\omega^\eps$, where at each $z=(x,t)\in \Z^2_{\rm even}$,
the walk moves to $(x+1, t+1)$ if $\omega^\eps_z=1$, moves to $(x-1, t+1)$ if $\omega^\eps_z=0$, and is a branching point if $\omega^\eps_z\in (0,1)$ since the walk has strictly positive probability of moving to either $(x+1, t+1)$ or $(x-1, t+1)$. Given the set of branching-coalescing random walk paths, which we denote by $N^\eps$ and call it a {\em discrete net}, we can then specify the value of $\omega^\eps_z$ at each branching point by sampling i.i.d.\ random variables $\bar\omega^\eps_z$ with common law $\bar\nu$ (see Figure~\ref{fig:netenviron}). The pair $(N^\eps, \bar\omega^\eps)$ then gives an alternative representation of the random environment $\omega$, wherein a walk must navigate along $N^\eps$, and when it encounters a branching point $z$, it jumps either left or right with probability $1-\bar\omega^\eps_z$, resp.~$\bar\omega^\eps_z$.

The above setup is essentially the same as the coupling between branching-coalescing random walks and coalescing random walks discussed in Section~\ref{S:couple}, which corresponds to taking $(\bar\omega^\eps_z)_{z\in \Z^2_{\rm even}}$ to be i.i.d.\ $\{0,1\}$-valued random variables in the current setting. As we rescale space and time by $\eps$ and $\eps^2$ respectively, we note that $N^\eps$ converges to a variant of the Brownian net $\Ni_{\beta_-, \beta_+}$, constructed from a pair of left-right Brownian webs as in Section~\ref{S:net} with respective drifts
\be\label{betapm}
\beta_- =  \beta-2\int\nu(\di q)(1-q)^{-1}, \qquad  \beta_+ = \beta+2\int\nu(\di q)q^{-1}.
\ee
The branching points of $N^\eps$ converge to the set of separation points of $\Ni_{\beta_-, \beta_+}$, denoted by $S$. Since the law of $\bar\omega^\eps_z$ at the branching points is $\bar\nu$, independent of $\eps$, we should then assign i.i.d.\ random variables $\omega_z$ with the same law $\bar\nu$ to each separation point $z\in S$. The pair $(\Ni_{\beta_-, \beta_+}, (\omega_z)_{z\in S})$ then gives the desired continuum space-time random environment, wherein a random motion $\pi$ must navigate along $\Ni_{\beta_-, \beta_+}$, and independently at each separation point $z$ it encounters, it chooses ${\rm sgn}_\pi(z)$ (see \eqref{sgnpi}) to be $+1$ with probability $\omega_z$ and $-1$ with probability $1-\omega_z$.

We next give a precise formulation of how the subclass of Howitt-Warren flows with $\int q^{-1}(1-q)^{-1}\nu(\di q)<\infty$ can be obtained from the random environment $(\Ni_{\beta_-, \beta_+}, (\omega_z)_{z\in S})$ (cf.~\cite[Theorem 4.7]{SSS14}). For its proof, see \cite{SSS14}. To construct a version of the Howitt-Warren flow which a.s.\ satisfies the Chapman-Kolmogorov equation (condition (i') in Remark~\ref{R:CK}), we will sample a collection of coalescing paths $\Wi$ given the environment $(\Ni_{\beta_-, \beta_+}, (\omega_z)_{z\in S})$.

\bt[Constructing Howitt-Warren flows in a Brownian net] \label{T:HWnet}
Let $\beta\in\R$ and let $\nu$ be a finite measure on $[0,1]$ with $\int q^{-1}(1-q)^{-1}\nu(\di q)<\infty$. Let $\Ni_{\beta_-, \beta_+}$ be a
Brownian net with drifts $\beta_-,\beta_+$ defined as in \eqref{betapm}, and let $S$ be its set of separation points. Conditional on $\Ni_{\beta_-, \beta_+}$, let $\omega:=(\omega_z)_{z\in S}$ be i.i.d.\ $[0,1]$-valued random variables with law $\bar\nu$
defined as in \eqref{nubar}. Conditional on $(\Ni_{\beta_-, \beta_+},\omega)$, let $(\alpha_z)_{z\in S}$ be independent $\{-1,+1\}$-valued random variables such that $\P[\alpha_z=1\,|\,(\Ni_{\beta_-, \beta_+},\omega)]=\omega_z$. Then
\be\label{webinnet4}
\Wi:=\{\pi\in\Ni_{\beta_-, \beta_+}: {\rm sgn}_\pi(z)=\alpha_z\ \forall z\in S \mbox{ s.t.\ $\pi$ enters }z\}.
\ee
is distributed as a Brownian web with drift $\beta$. For any $z=(x,t) \in\R^2$, if $z$ is of type $(1,2)$, then let $\pi_z^\up\in \Wi(z)$ be any path in $\Wi$ entering $z$ restricted to the time interval $[t, \infty)$; otherwise let $\pi_z^\up$ be the rightmost path in
$\Wi(z)$. Then
\bc\label{HWconst2}
K^\up_{s,t}(x,\cdot):=\P\big[\pi^\up_{(x,s)}(t)\in \cdot\,\big|\,(\Ni_{\beta_-, \beta_+},\omega)\big], \qquad s\leq t,\ x\in\R
\ec
defines a version of the Howitt-Warren flow with drift $\bet$ and characteristic measure $\nu$, which satisfies condition {\rm (i')} in Remark~\ref{R:CK}.
\et

\medskip

{\bf General case:} $\nu$ is any finite measure on $[0,1]$. Let $(\mu_\eps)_{\eps>0}$ satisfy \eqref{mucon} and let $\omega^\eps:=(\omega^\eps_z)_{z\in \Z^2_{\rm even}}$ be an i.i.d.\ random environment with common law $\mu_\eps$. Without assuming
$\int q^{-1}(1-q)^{-1}\nu(\di q)<\infty$, it may no longer be possible to capture the continuum random environment by a Brownian net, because either $\beta_-$ or $\beta_+$  in \eqref{betapm} or both can be infinity. Instead, we follow the alternative view on the coupling between branching-coalescing and coalescing random walks in Section~\ref{S:couple}, where we first sample the coalescing walks and then introduce the branching points. The same approach can be applied here to encode the random environment $\omega$ in two steps.

First we construct a collection of coalescing random walks by sampling coalescing walks in the same environment $\omega$ and then average over the law of $\omega^\eps$. More precisely, from each $z=(x,t)\in \Z^2_{\rm even}$, the walk moves to either $(x+1, t+1)$ or $(x-1, t+1)$, represented by $\alpha^\eps_z:=1$ or $\alpha^\eps_z=-1$,  with respective probability $\int_0^1 q \mu_\eps({\rm d}q)$ and $\int_0^1 (1-q) \mu_\eps({\rm d}q)$. Let $W^\eps$ denote this collection of coalescing random walks, which we call a {\em discrete web} and has a natural dual $\widehat W^\eps$ as illustrated in Figure~\ref{fig:webdual}. Next, we identify the law of $\omega^\eps$ conditioned on $W^\eps$.

Note that conditioned on $W^\eps$, $(\omega^\eps_z)_{z\in \Z^2_{\rm even}}$ are independent with conditional distribution
\be\label{mulr}
\mu_\eps^{\rm l} :=\frac{(1-q)\mu_\eps({\rm d}q)}{\int (1-q)\mu_\eps({\rm d}q)} \quad \mbox{ if } \alpha^\eps_z=-1; \qquad
\mu_\eps^{\rm r} := \frac{q\mu_\eps({\rm d}q)}{\int q\mu_\eps({\rm d}q)} \quad \mbox{ if } \alpha^\eps_z=1.
\ee
Therefore conditioned on $W^\eps$, if we sample i.i.d.\ random variables $(\eta^\eps_z)_{\{z:\alpha_z^\eps=-1\}}$ with common law $\mu_\eps^{\rm l}$ and i.i.d.\ random variables $(1-\eta^\eps_z)_{\{z: \alpha^\eps_z=1\}}$ with common law $\mu_\eps^{\rm r}$, then $(W^\eps, \eta^\eps)$ provides an alternative representation for the environment $\omega^\eps$, where at each $z\in\Z^2_{\rm even}$, a walk in the random environment follows the same jump as in $W^\eps$ with probability $1-\eta^\eps_z$, and jumps in the opposite direction with probability $\eta^\eps_z$ (see Figure~\ref{fig:webenviron}).

The above setup is an extension of the coupling between coalescing and branching-coalescing random walks discussed in Section~\ref{S:couple}, which corresponds to taking $(\eta^\eps_z)_{z\in \Z^2_{\rm even}}$ to be i.i.d.\ $\{0,1\}$-valued random variables with mean
$\eps$. As we rescale space and time by $\eps$ and $\eps^2$ respectively, we note that $W^\eps$ converges to a Brownian web $\Wi_0$ with drift $\beta$. On the other hand, $\eta^\eps$ can be seen to converge to a marked Poisson point process on the set of intersection points between paths in $\Wi_0$ and paths in its dual $\widehat \Wi_0$, i.e., the $(1,2)$ points of $\Wi_0$. Indeed, we can regard $(\eta^\eps_z)_{\{z:\alpha_z^\eps=-1\}}$ as a marked point process on the set of space-time points in $\Z^2_{\rm even}$ where a random walk path in $W^\eps$ is exactly one unit of distance to the left of a dual random walk path in $\widehat W^\eps$. As noted in Section~\ref{S:couple}, when space-time is rescaled by $(\eps, \eps^2)$ and measure is rescaled by $\eps$, the counting measure on these points of collision between forward and dual random walk paths converge to the intersection local time measure $\ell_{\rm l}$ on points of type $(1,2)_{\rm l}$ in the Brownian web. Since for every $u\in (0,1)$, by \eqref{mucon}, the mark $\eta^\eps_z$ at each $z\in \Z^2_{\rm even}$ with $\alpha^\eps_z=-1$ satisfies
$$
\P(\eta^\eps_z\in (u, 1]) = \frac{\int_u^1(1-q)\mu_\eps({\rm d}q)}{\int_0^1 (1-q)\mu_\eps({\rm d}q)}
= \frac{\eps \int_u^1 q^{-1} \eps^{-1}q(1-q)\mu_\eps({\rm d}q)}{\int_0^1 (1-q)\mu_\eps({\rm d}q)}
\sim \eps\int_u^1 \frac{2}{q}\nu({\rm d}q) \quad \mbox{as } \eps\downarrow 0,
$$
it follows that $(\eta^\eps_z)_{\{z:\alpha_z^\eps=-1\}}$ converges to a marked Poisson point
process $\Mi_{\rm l}\subset\R^2\times [0,1]$ with intensity measure $\ell_{\rm l}({\rm d}z)\otimes \frac{2}{q} 1_{\{q>0\}}\nu({\rm d}q)$. Similarly, $(1-\eta^\eps_z)_{\{z:\alpha_z^\eps=+1\}}$ converges to a marked Poisson point process $\Mi_{\rm r}\subset\R^2\times [0,1]$ with intensity measure $\ell_{\rm r}({\rm d}z)\otimes \frac{2}{1-q} 1_{\{q<1\}}\nu({\rm d}q)$.

The triple $(\Wi_0, \Mi_{\rm l}, \Mi_{\rm r})$ then gives the desired continuum space-time random environment, wherein a random motion $\pi$ must navigate along paths in $\Wi_0$, and independently at each marked $(1,2)_{\rm l}$ point $(z, \eta_z)\in \Mi_{\rm l}$, or marked $(1,2)_{\rm r}$ point $(z, 1-\eta_z)\in \Mi_{\rm r}$ it encounters, with probability $\eta_z$, $\pi$ chooses orientation $-{\rm sgn}_{\Wi_0}(z)$, i.e., $\pi$ switches to the second outgoing path in $\Wi_0(z)$ instead of continuing along the incoming path. This description is correct when $\nu(\{0\})=\nu(\{1\})=0$. However when $\nu(\{0\})>0$, we have excluded $\ell_{\rm l}({\rm d}z)\otimes \frac{2}{q} 1_{\{q=0\}}\nu({\rm d}q)$ from the intensity measure of $\Mi_{\rm l}$, which has a non-negligible effect on paths sampled in the random environment.
To understand this effect, we can approximate $\nu(\{0\})\delta_0({\rm d}q)$ by $\nu(\{0\}) \delta_h({\rm d}q)$ with $h\downarrow 0$.
This leads to Poisson marking $(1,2)_{\rm l}$ points of $\Wi_0$ with intensity measure $\frac{2}{h} \nu(\{0\})\ell_{\rm l}({\rm d}z)$, and whenever a random motion $\pi$ in this environment encounters such a Poisson point $z$, it chooses its orientation at $z$ to be $-{\rm sgn}_{\Wi_0}(z)$ with probability $h$. The net effect is that, at each point $z$ of a Poisson point process $B_{\rm l}$ with intensity measure $2 \nu(\{0\})\ell_{\rm l}({\rm d}z)$, $\pi$ chooses its orientation at $z$ to be $-{\rm sgn}_{\Wi_0}(z)$ with probability 1. As $h\downarrow 0$, the resulting effect is that, when we sample an independent motion $\pi'$ in the same random environment, then an independent copy of $B_{\rm l}$, call it $B_{\rm l}'$, must be sampled such that whenever $\pi'$ encounters some $z\in B_{\rm l}'$, it chooses its orientation at $z$ to be  $-{\rm sgn}_{\Wi_0}(z)$ with probability $1$.

We now give a precise formulation how the Howitt-Warren flow can be obtained from a random environment $(\Wi_0, \Mi_{\rm l}, \Mi_{\rm r})$ (cf.~\cite[Theorem 3.7]{SSS14}). As in Theorem~\ref{T:HWnet}, we will sample a collection of coalescing Brownian motions $\Wi$ given the environment $(\Wi_0, \Mi_{\rm l}, \Mi_{\rm r})$.

\bt{\bf(Construction of Howitt-Warren flows)}\label{T:HWweb}
Let $\beta\in\R$ and let $\nu$ be a finite measure on $[0,1]$. Let $\Wi_0$ be a Brownian web with drift $\beta$. Let $\Mi$ be a marked Poisson point process on $\R^2\times [0,1]$ with intensity measure
$$
\ell_{\rm l}({\rm d}z)\otimes \frac{2}{q} 1_{\{q>0\}}\nu({\rm d}q) + \ell_{\rm r}({\rm d}z)\otimes \frac{2}{1-q} 1_{\{q<1\}}\nu({\rm d}q).
$$
Conditional $(\Wi_0, \Mi)$, let $\alpha_z$ be independent $\{-1,+1\}$-valued random variables with $\P[\alpha_z=+1|(\Wi_0,\Mi)]=\omega_z$ for each $(z, \omega_z)\in \Mi$, and let
$$
A:=\{z: (z, \omega_z)\in \Mi, \alpha_z\neq {\rm sgn}_{\Wi_0}(z)\}.
$$
Let $B$ be an independent Poisson point set with intensity $2\nu_{\rm l}(\{0\})\ell_{\rm l}+2\nu_{\rm r}(\{1\})\ell_{\rm r}$. Define
\be\label{webswitch}
\Wi := \lim_{\Delta_n\up A\cup B} {\rm switch}_{\Delta_n}(\Wi_0)
\ee
for any sequence of finite sets $\Delta_n \up A\cup B$, where ${\rm switch}_{\Delta_n}(\Wi_0)$ is the set of paths obtained from $\Wi_0$ by redirecting all paths in $\Wi$ entering any $z\in \Delta_n$ in such a way that $z$ of type $(1,2)_{\rm l}$ in $\Wi$ becomes type
$(1,2)_{\rm r}$ in $\Wi_0$ and vice versa. Then $\Wi$ is equally distributed with $\Wi_0$, and with $\pi^\up_z$ defined as in Theorem~\ref{T:HWnet},
\bc\label{HWconst}
K^\up_{s,t}(x,\cdot):=\P\big[\pi^\up_{(x,s)}(t)\in \cdot\,\big|\,(\Wi_0,\Mi)\big] \qquad s\leq t,\ x\in\R,
\ec
defines a version of the Howitt-Warren flow with drift $\bet$ and characteristic measure $\nu$, which satisfies condition {\rm (i')} in Remark~\ref{R:CK}.
\et

\bex Use Theorem~\ref{T:coupling} to show that $\Wi$ as defined in \eqref{webswitch} does not depend on the choice of $\Delta_n\up A\cup B$, and $\Wi$ is equally distributed with $\Wi_0$.
\eex

\subsection{Properties of the Howitt-Warren flow}\label{S:HWprop}

The construction of the continuum space-time random environment underlying the Howitt-Warren flows allows the study of almost sure properties of the flow $K^\up :=(K^\up_{s,t}(x, \cdot))_{s<t, x\in\R}$ (cf.~\cite[Section 2]{SSS14}). In particular, we can study almost sure path properties of the following measure-valued process induced by the Howitt-Warren flow, called the {\em Howitt-Warren process}:
\be\label{HWprocess}
\rho_t({\rm d}y) = \int K^\up_{0,t} (x, {\rm d}y) \rho_0({\rm d}x),
\ee
where $\rho_0$ can be taken to be any locally finite measure on $\R$ in the class
\be
\Mi_{\rm g}(\R) := \Big\{ \rho : \int e^{-cx^2} \rho({\rm d}x) <\infty \mbox{ for all } c>0 \Big\},
\ee
where $\rho_n\in \Mi_{\rm g}(\R)$ is defined to converge to $\rho \in \Mi_{\rm g}(\R)$ if $\int f(x)\rho_n({\rm d}x) \to \int f(x) \rho({\rm d}x)$ for all $f\in C_c(\R)$, which is just convergence in the vague topology, plus $\int e^{-cx^2} \rho_n({\rm d}x) \to \int e^{-cx^2} \rho({\rm d}x)$ for all $c>0$.

A first consequence is that $(\rho_t)_{t\geq 0}$ is a Markov process with continuous sample path, and continuous dependence on the initial condition and starting time, which makes it a Feller process.

\bt[Howitt-Warren process]\label{T:Feller}
Let $(K^\uparrow_{s,t}(x, \cdot))_{s\leq t, x\in\R}$ be the version of Howitt-Warren flow with drift $\beta$ and characteristic measure $\nu$,
defined in Theorem~\ref{T:HWnet} or \ref{T:HWweb}. Let $(\rho_t)_{t\geq 0}$ be the Howitt-Warren process defined from $K^\up$ as in \eqref{HWprocess} with $\rho_0\in \Mi_{\rm g}$. Then
\begin{itemize}
\item[\rm (i)] $(\rho_t)_{t\geq 0}$ is an $\Mi_{\rm g}(\R)$-valued Markov process with almost sure continuous sample paths;

\item[\rm (ii)] If $(\rho^{\li n\re}_t)_{t\geq s_n}$ are Howitt-Warren processes defined from $K^\up$  with deterministic initial condition $\rho^{\li n\re}_{s_n}$ at time $s_n$, with $s_n\to 0$, then for any $t>0$ and $t_n\to t$,
\be\label{rhonconv}
\rho^{\li n\re}_{s_n}\Asto{n}\rho_0\quad\mbox{implies}\quad
\rho^{\li n\re}_{t_n}\Asto{n}\rho_t\qquad a.s.,
\ee
where $\Rightarrow$ denotes convergence in $\Mi_{\rm g}(\R)$.
\end{itemize}
If $\int q^{-1} (1-q)^{-1} \nu({\rm d}q)<\infty$, then the above statements hold with $\Mi_{\rm g}(\R)$ replaced by $\Mi_{\rm loc}(\R)$, the space of locally finite measures on $\R$ equipped with the vague topology.
\et

We can also identify almost surely the support of $\rho_t$ for all $t\geq 0$.
\bt[Support of Howitt-Warren process]\label{T:supp}
Let $(\rho_t)_{t\geq 0}$ be a Howitt-Warren process with drift $\bet$ and characteristic measure, and the initial condition $\rho_0$ has compact support. Let $\beta_-,\beta_+$ be defined as in \eqref{betapm}.
\begin{itemize}
\item[{\bf(a)}] If $-\infty<\bet_-<\bet_+<\infty$, then a.s.\ for all $t>0$, the support of $\rho_t$ satisfies
\be\label{supbraco}
{\rm supp}(\rho_t)= \{\pi(t):  \pi \in \Ni_{\beta_-, \beta_+}, \pi(0)\in {\rm supp}(\rho_0) \},
\ee
where $\Ni_{\beta_-, \beta_+}$ is the Brownian net with drift parameters $\beta_-, \beta_+$ as in Theorem~\ref{T:HWnet}.

\item[{\bf(b)}] If $\bet_-=-\infty$ and $\bet_+<\infty$, then a.s.\ ${\rm
  supp}(\rho_t)=(-\infty,r_t]\cap\R$ for all $t>0$, where $r_t:=\sup({\rm
  supp}(\rho_t))$. An analogue statement holds when $\bet_->-\infty$ and
$\bet_+=\infty$.
\item[{\bf(c)}] If $\bet_-=-\infty$ and $\bet_+=\infty$, then a.s.\ ${\rm
  supp}(\rho_t)=\R$ for all $t>0$.
\end{itemize}
\et

Theorem~\ref{T:supp}~(a) shows that when the Howitt-Warren flow can be constructed from a Brownian net, then at deterministic
times, $\rho_t$ is almost surely atomic. This result can be extended to general Howitt-Warren flows. However, almost surely there exist random times when $\rho_t$ contains no atoms, and the only exceptions are the Arratia flow with $\nu=0$, and the {\em erosion flows}, which have characteristic measures of the form $\nu= a\delta_0+ b\delta_1$ with $a, b\geq 0$ and $a+b>0$.

\bt[Atomicness vs non-atomicness] \label{T:atom}
Let $(\rho_t)_{t\geq 0}$ be a Howitt-Warren process with drift $\beta$ and characteristic measure $\nu$.
\begin{itemize}
\item[{\bf(a)}] For each $t>0$, $\rho_t$ is a.s.\ purely atomic.
\item[{\bf(b)}] If $\int_{(0,1)}\nu({\rm d} q)>0$, then a.s.\ there exists a dense set of random times $t>0$ when $\rho_t$ contains no atoms.
\item[{\bf(c)}] If $\int_{(0,1)}\nu({\rm d} q)=0$, then a.s.\ $\rho_t$ is purely atomic at all $t>0$.
\end{itemize}
\et

We can also study ergodic properties of the Howitt-Warren process. In particular, for any Howitt-Warren process other than the measure-valued process generated by the coalescing Arratia flow, there is a unique spatially ergodic stationary law for $(\rho_t)_{t\geq 0}$, which
is also the weak limit of $\rho_t$ as $t\to\infty$ if the law of the initial condition $\rho_0$ is spatially ergodic with finite mean density.

\bt[Ergodic properties] \label{T:erg}
Let $(\rho_t)_{t\geq 0}$ be a Howitt-Warren process with drift $\beta\in\R$, characteristic measure $\nu \neq 0$, and initial law $\Li(\rho_0)$.
\begin{itemize}
\item[\rm (i)] If $\Li(\rho_0)$ is ergodic w.r.t.\ $T_a \rho_0(\cdot)=\rho_0(a+\cdot)$ for all $a\in \R$, and $\E[\rho_0([0,1])]=1$, then as $t\to\infty$, $\Li(\rho_t)$ converges weakly  to a limit $\Lambda_1$ which is also ergodic with respect to $T_a$ for all $a\in\R$. Furthermore, if $\Li(\rho_0)=\Lambda_1$, then $\E[\rho_0([0,1])]=1$ and $\Li(\rho_t)=\Lambda_1$ for all $t>0$.

\item[\rm (ii)] If $\Li(\rho_0)$ is ergodic w.r.t.\ $T_a$ for all $a\in\R$ and $\E[\rho_0([0,1])]=\infty$, then as $t\to\infty$, $\Li(\rho_t)$ has no sub-sequential weak limits supported on $\Mi_{\rm loc}(\R)$, the set of locally finite measures on $\R$.
\end{itemize}
\et
\br
When $\nu({\rm d}x)=1_{[0,1]}(x){\rm d}x$, it is known from \cite[Prop.~9(b)]{LR04b} that $\Lambda_1$ is the law of a random measure
$\rho^*= \sum_{(x,u)\in {\cal P}} u\delta_x$ for a Poisson point process $\cal P$ on $\R\times [0,\infty)$ with intensity measure ${\rm d}x \times  u^{-1} e^{-u}{\rm d}u$.
\er

The proof of the above results can be found in \cite{SSS14}.

\br\label{R:HWfluc}
Howitt-Warren flows are the continuum analogues of discrete random walks in i.i.d.\ space-time random environments, and they share the same fluctuations on large space-time scales. In particular, for any Howitt-Warren process $(\rho_t)_{t\geq 0}$ (assuming drift $\beta=0$ and $\rho_0({\rm d}x)={\rm d}x$ for simplicity),  the rescaled current process $({\cal I}(nt, \sqrt{n}x)/n^{1/4})_{t>0, x\in\R}$, where
$$
{\cal I}(t,x) := \int_\R \int_R 1_{\{u<0\}} 1_{\{v>x\}}K^\uparrow_{0,t}(u, {\rm d}v)  \rho_0({\rm d}u) - \int_\R \int_R 1_{\{u>0\}} 1_{\{v<x\}}K^\uparrow_{0,t}(u, {\rm d}v)  \rho_0({\rm d}u),
$$
converges to a universal Gaussian process as $n\to\infty$. This was shown in~\cite{Y14}, and the same universal (Edwards-Wilkinson) fluctuations have been established earlier for random walks in i.i.d.\ space-time random environments (see e.g.~\cite{S10}). In a different direction, the transition probabilities of a random walk in i.i.d.\ space-time random environments are believed to have the same universal (Tracy-Widom GUE) fluctuations as the point-to-point partition functions of a directed polymer. This was verified recently in~\cite{BC15} for special Beta-distributed random environments, and one expects similar results to hold for the Howitt-Warren flow $K^\uparrow_{0,t}(0, {\rm d}x)$.
\er

\section{Convergence to the Brownian web and net}\label{S:conv}

In this section, we give general convergence criteria for the Brownian web that were originally formulated in \cite{FINR04, NRS05},
and simplified criteria when paths do not cross each other.  We will discuss strategies for verifying these criteria and focus in particular on the convergence of coalescing random walks to the Brownian web. Lastly we will formulate a set of convergence criteria for the Brownian net and verify them for branching-coalescing simple random walks.

\subsection{General convergence criteria for the Brownian web}\label{SS:conv}

We give here general convergence criteria for a sequence of random variables $(X_n)_{n\in\N}$, taking values in the space of compact sets of paths $\Hi$, to converge in distribution to the Brownian web $\Wi$. When $X_n$ consists of non-crossing paths, these criteria can be greatly simplified, which will be discussed in the next subsection.

First we formulate criterion which ensures tightness for the laws of $(X_n)_{n\in\N}$. We then formulate criteria which ensure that any subsequential limit of $(X_n)_{n\in\N}$ contains a copy of the Brownian web (lower bound), but nothing more (upper bound).
\medskip

\noindent
{\bf Tightness:}
To understand the tightness criterion we will formulate, let us first see what should be the tightness criterion for a sequence of path-valued random variables $(Y_n)_{n\in\N}$, i.e., $Y_n\in \Pi$ with the space of paths $\Pi$ defined as in Section~\ref{SS:hausdorff}. Due to the compactification of $\R^2$ in Section~\ref{SS:hausdorff}, it suffices to show that when restricted to any finite space-time window $\Lambda_{L,T}=[-L, L]\times [-T, T]$, the law of the random paths $(Y_n)_{n\in\N}$ are tight in the sense that for any $\eps>0$,
there exists a modulus of continuity $\phi: [0,1]\to \R$, which is increasing with $\phi(\delta)\downarrow 0$ as $\delta\down 0$,
such that uniformly in $n\in\N$,
\be\label{modcon1}
\P\big( \forall\, s \mbox{ with } (Y_n(s), s) \in \Lambda_{L,T} \mbox{ and } t\in [s, s+1], \ |Y_n(t)-Y_n(s)|\leq \phi(t-s) \big) \geq 1-\eps.
\ee
The modulus of continuity $\phi$ allows the construction of an equicontinuous, and hence compact set of paths. To construct a $\phi$ that
satisfies \eqref{modcon1}, it suffices to show that for any $\eta>0$,
\be\label{modcon2}
\lim_{\delta\downarrow 0} \limsup_{n\to\infty} \P(\exists\,  s<t<s+\delta \mbox{ with } (Y_n(s), s) \in \Lambda_{L,T}, \ s.t.\  |Y_n(t)-Y_n(s)| >\eta) =0.
\ee
Indeed, fix a sequence $\eta_m\downarrow 0$. Then for each $m\in\N$, by \eqref{modcon2}, we can find $\delta_m>0$ sufficiently small such that uniformly in $n\in\N$, the probability in \eqref{modcon2} is bounded by $\eps/2^m$. We can then define $\phi(h):=\eta_m$ for
$h\in (\delta_{m-1}, \delta_m]$, which is easily seen to satisfy \eqref{modcon1}.

When we consider a sequence of random compact sets of paths $(X_n)_{n\in\N}$, the tightness criterion is similar, except that we need to control the modulus of continuity uniformly for all paths in $X_n$ (cf.~Exercise~\ref{E:comsub}). Therefore condition \eqref{modcon2} should be modified to require showing that
for any finite $\Lambda_{L, T}$ and any $\eta>0$,
\be\label{modcon3}
\lim_{\delta\downarrow 0} \limsup_{n\to\infty} \P(\exists\,  \pi\in X_n \mbox{ and } s<t<s+\delta \mbox{ with } (\pi(s), s) \in \Lambda_{L,T}, \ s.t.\  |\pi(t)-\pi(s)| >\eta) =0.
\ee
To control the modulus of continuity of all paths in $X_n$ simultaneously, it is convenient to divide $\Lambda_{L,T}$ into
$16LT/\delta \eta$ sub-rectangles of dimension $\eta/2\times \delta/2$, and bound the event in \eqref{modcon3} by the union of events
where $\Lambda_{L,T}$ in \eqref{modcon3} is replaced by one of the $16LT/\delta \eta$ sub-rectangles of $\Lambda_{L,T}$. This leads to the following tightness criterion as formulated in \cite[Prop.~B1]{FINR04}.

\begin{figure}
\begin{center}
\includegraphics[width=11cm]{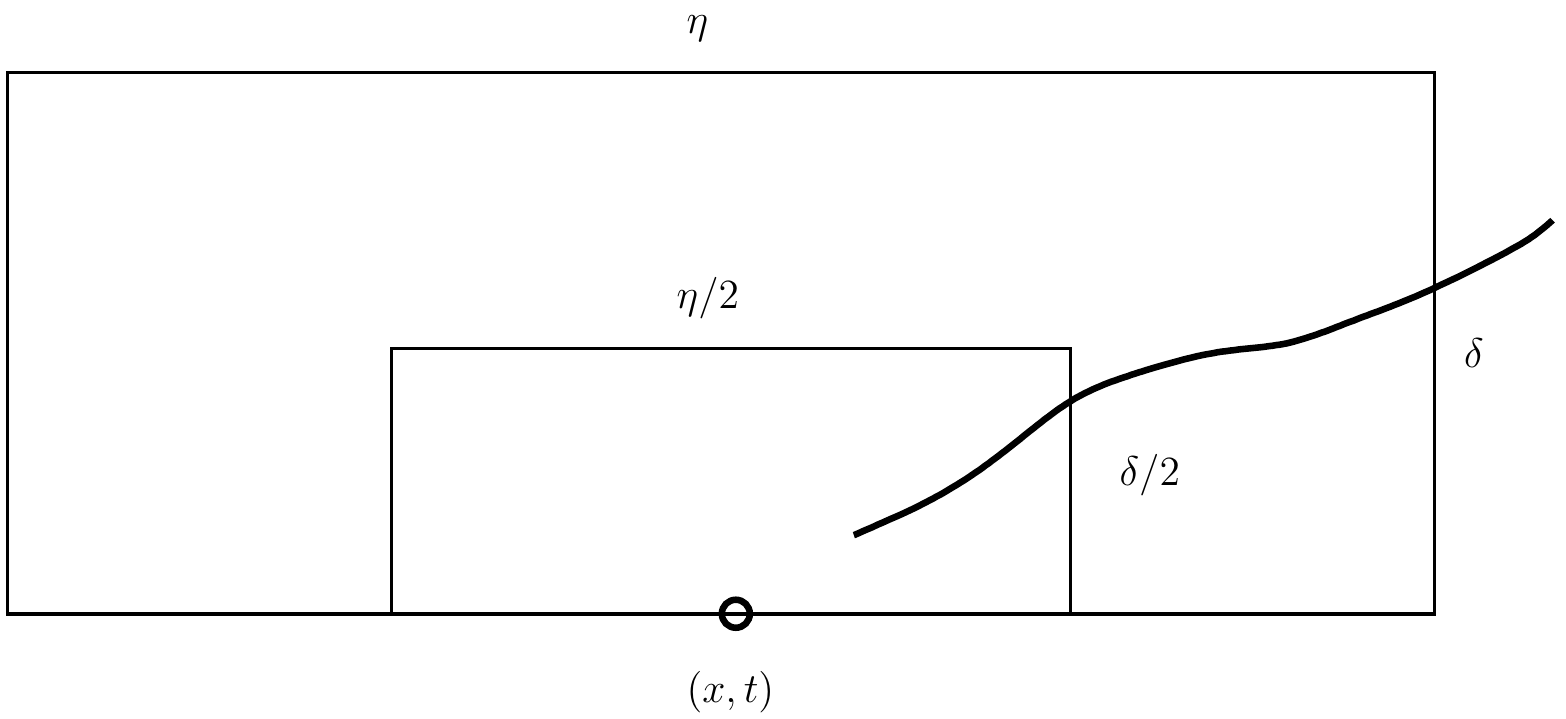}
\caption{A path causing the event $A_{\delta, \eta}(x,t)$ to occur.}
\label{fig:tight}
\end{center}
\end{figure}

\bp[Tightness criterion] \label{P:tight1} The law of a sequence of $(\Hi, \Bi_\Hi)$-valued random variables $(X_n)_{n\in\N}$ is tight if
$$
 {\rm (T)} \hspace{2cm} \forall \,  L,T\in  (0,\infty), \quad  \lim_{\delta\downarrow 0} \delta^{-1} \limsup_{n\to\infty} \sup_{(x,t)\in [-L,L]\times [-T, T]} \P(X_n \in A_{\delta, \eta}(x,t)) = 0, \hspace{4cm}
$$
where $A_{\delta, \eta}(x,t) \in \Bi_\Hi$ consists of compact sets of paths $K\in \Hi$, such that $K$ contains some path
which intersects the rectangle $[x-\eta/4, x+\eta/4]\times [t, t+\delta/2]$, and at a later time, intersects the left or right boundary of the bigger rectangle $[x-\eta/2, x+\eta/2]\times [t, t+\delta]$ $($see Figure~\ref{fig:tight}$)$.
\ep
\medskip

\noindent
{\bf Lower bound:} Assuming that $(X_n)_{n\in\N}$ is a tight sequence of $\Hi$-valued random variables, we then need criterion to ensure that any subsequential weak limit of $(X_n)_{n\in\N}$ contains almost surely a random subset which is distributed as the standard Brownian web. The following criterion serves this purpose, which is a form of convergence in finite-dimensional distributions.
$$
{\rm (I)}  \quad
\begin{aligned}
&\mbox{There exists $\pi_{n,z}\in X_n$ for each $z\in \R^2$, such that for any deterministic $z_1, \ldots, z_k\in\R^2$,}\\
&\mbox{$(\pi_{n,z_i})_{1\leq i\leq k}$ converge in distribution to coalescing Brownian motions starting at $(z_i)_{1\leq i\leq k}$.}
\end{aligned}
$$
If $(X_n)_{n\in\N}$ is a tight sequence that satisfies condition (I), then for any deterministic countable dense set $\Di\subset\R^2$,
we note that by going to a further subsequence if necessary, any subsequential limit $\cal X$ of $(X_n)_{n\in\N}$ contains a collection of coalescing Brownian motions $\{\pi_z\}_{z\in \Di}$ starting from each $z\in\Di$. Since $\Wi:=\overline{\{\pi_z: z\in\Di\}}$ is a standard Brownian web
by Theorem~\ref{T:webchar}, we obtain the desired lower bound that a.s.\ $\Wi\subset \cal X$.
\bigskip

\noindent
{\bf Upper bound:} Assuming that $(X_n)_{n\in\N}$ is a tight sequence of $\Hi$-valued random variables that satisfies condition (I), then it only remains to formulate criterion to ensure that any subsequential weak limit $\cal X$ of $(X_n)_{n\in\N}$ a.s.\ contains no more path
than the Brownian web $\Wi\subset \cal X$. There are several approaches to this problem, depending partly on whether paths in $X_n$ can cross each other or not.

One way is to control the expectation of the following family of counting random variables
\be\label{eta}
\eta_{\cal X}(t, h; a, b) = \big|\{ \pi(t+h) : \pi \in {\cal X}, \pi(t)\in [a,b]  \}\big|, \qquad t\in\R, h>0, a<b,
\ee
which considers all paths in $\cal X$ that intersect $[a,b]$ at time $t$ and counts the number of distinct positions these paths occupy at time $t+h$. Thanks to the image set property of the Brownian web $\Wi$, which is a special case of Prop.~\ref{P:imgset} for the Brownian net, it is easily seen that if $\cal X$ contains strictly more paths than $\Wi$, than we must have $\eta_{\cal X}(t,h; a,b)> \eta_\Wi(t, h; a,b)$ for some rational $t, h, a,b$. Therefore to show $\cal X=\Wi$ a.s., it suffices to show that
\be\label{den1}
\E[\eta_{\cal X}(t,h; a,b)] = \E[\eta_\Wi(t, h; a,b)] \qquad \mbox{for all } t\in \R, h>0, a<b.
\ee

The following sufficient criteria have been formulated in \cite{FINR04} to ensure that \eqref{den1} holds for any subsequential limit $\cal X$ of a sequence of $\Hi$-valued random variables $(X_n)_{n\in\N}$:
$$
\begin{aligned}
&{\rm (B1')} \quad \quad \forall\, h_0>0, \ \limsup_{n\to\infty}   \sup_{h> h_0}\sup_{a, t\in \R} \P\big[\eta_{X_n}(t,h; a, a+\eps)\geq 2\big]\aston{\eps} 0,\\
&{\rm (B2')} \quad \quad \forall\, h_0>0, \ \frac{1}{\eps}\limsup_{n\to\infty}   \sup_{h> h_0}\sup_{a, t\in \R}
\P\big[X_n(t,h;a, \eps) \neq X^-_n(t,h;a, \eps) \cup X^+_n(t,h;a, \eps)\big]\aston{\eps} 0,
\end{aligned}
\qquad
$$
where $X_n(t,h;a,\eps)\subset \R$ is the set of positions occupied at time $t+h$ by paths in $X_n$ which intersect the interval $[a, a+\eps]$ at time $t$, while $X_n^-(t, h; a, \eps)$ (resp.\  $X_n^+(t, h; a, \eps)$) is the subset of $X_n(t,h,; a, \eps)$ induced by paths in $X_n$ which occupy the leftmost (resp.\ rightmost) position at time $t$ among all paths in $X_n$ that intersect $[a, a+\eps]$ at time $t$.
Condition (B1') ensures that for each deterministic point $z\in \R^2$, any subsequential limit $\cal X$ contains a.s.\ at most one path starting from $z$. Condition (B2') ensures that $\eta_{\cal X}(t,h; a,b)$ can be approximated by partitioning $[a,b]$ into equal-sized intervals and considering only paths in $\cal X$ starting at the boundaries of these intervals. Condition (B1') and (I) together ensure that paths in $\cal X$ starting at these boundary points are distributed as coalescing Brownian motions. Condition \eqref{den1} then follows.
\medskip

We thus have the following convergence criteria for the Brownian web~\cite[Theorem~5.1]{FINR04}.
\bt[Convergence criteria A] \label{T:webconvA} Let $(X_n)_{n\in\N}$ be a sequence of $(\Hi, \Bi_\Hi)$-valued random variables satisfying conditions {\rm (T), (I), (B1'), (B2')}. Then $X_n$ converges in distribution to the standard Brownian web $\Wi$.
\et

Condition (B2') turns out to be difficult to verify when paths in $X_n$ can cross each other. This is in particular the case for non-nearest neighbor coalescing random walks on $\Z$, which led to the formulation in~\cite{NRS05} of an alternative criterion in place of (B2'). The observation is that instead of $\eta_{\cal X}(t, h; a, b)$, we can also consider the alternative family of counting random variables
\be\label{etahat}
\hat \eta_{\cal X}(t, h; a, b) = \big|\{ \pi(t+h) \cap (a,b) : \pi \in {\cal X}, \pi(t)\in\R  \}\big|, \qquad t\in\R, h>0, a<b,
\ee
which considers all paths in $\cal X$ starting before or at time $t$ and counts the number of distinct positions these paths occupy in the interval $(a,b)$ at time $t+h$. Given ${\cal X}\supset \Wi$, to show that ${\cal X}=\Wi$ a.s., it also suffices to show that
\be\label{den2}
\E[\hat \eta_{\cal X}(t,h; a,b)] = \E[\hat \eta_\Wi(t, h; a,b)] \qquad \mbox{for all } t\in \R, h>0, a<b.
\ee
This leads to the following alternative convergence criteria in \cite[Theorem 1.4]{NRS05}.
\bt[Convergence criteria B] \label{T:webconvB}
Let $(X_n)_{n\in\N}$ be a sequence of $(\Hi, \Bi_\Hi)$-valued random variables satisfying conditions
{\rm (T), (I), (B1')} and the following condition
$$
{\rm (E)} \quad \mbox{For any subsequential weak limit } {\cal X}, \  \E[\hat \eta_{\cal X}(t,h; a,b)] = \E[\hat \eta_\Wi(t, h; a,b)] \quad \forall\, t\in \R, h>0, a<b.
$$
Then $X_n$ converges in distribution to the standard Brownian web $\Wi$ as $n\to\infty$.
\et
Condition (E) is actually much easier to verify than it appears. It turns out to be enough to establish the density bound
\be\label{Edensbd}
\limsup_{n\to\infty}  \E[\hat \eta_{X_n}(t,h; a,b)] <\infty \qquad \forall\,  t\in \R, h>0, a<b,
\ee
which by Fatou's Lemma implies that any subsequential weak limit $\cal X$ satisfies $\E[\hat \eta_{\cal X}(t,\eps; a,b)]<\infty$ for all
$t\in \R$, $\eps>0$ and $a<b$. In particular,
\be
{\cal X}^{t-}(t+\eps) :=\{ \pi (t+\eps) : \pi \in {\cal X}^{t-}\} \quad \mbox{with} \quad {\cal X}^{t-}:=\{\pi \in {\cal X}: \sigma_\pi \leq t\},
\ee
the set of positions at time $t+\eps$ generated by paths in $\cal X$ starting before or at time $t$, is almost surely a locally finite subset of $\R$. As a random subset of $\R$, ${\cal X}^{t-}(t+\eps)$ arises as the limit of $X_n^{t-}(t+\eps)$,\footnote{In practice, ${\cal X}^{t-}(t+\eps)$ may contain positions which arise from limits of paths in $X_n$ that start at times $t_n\downarrow t$. Therefore we should consider
instead $X_n^{(t+\delta)-}(t+\eps)$ for some $\delta>0$, so that its limit contains ${\cal X}^{t-}(t+\eps)$.} and we can use Skorohod's representation theorem~\cite{B99} to couple them such that almost surely, $X_n^{t-}(t+\eps)$ converges to ${\cal X}^{t-}(t+\eps)$ w.r.t.\ the Hausdorff metric on subsets of $\R$. For Markov processes such as coalescing random walks, we expect that the law of paths in $X^{t-}_n$ restricted to the time interval $[t+\eps, \infty)$ depends only on their positions $X^{t-}_n(t+\eps)$ at time $t+\eps$, and furthermore, conditions (B1') and (I) can be applied to these restricted paths conditioned on $X^{t-}_n(t+\eps)$. Therefore given $X_n^{t-}(t+\eps)$ converging to ${\cal X}^{t-}(t+\eps)$, we can apply conditions (B1') and (I) to conclude that, ${\cal X}^{t-}$ restricted to the time interval $[t+\eps, \infty)$ is a collection of coalescing Brownian motions starting from the locally finite set ${\cal X}^{t-}(t+\eps)$ at time $t+\eps$, and hence
$$
\E[\hat\eta_{\cal X}(t, h; a,b)] \leq \E[\hat \eta_\Wi(t+\eps, h-\eps; a,b)] = \frac{b-a}{\sqrt{\pi (h-\eps)}}.
$$
Sending $\eps\downarrow 0$ then establishes condition (E). For non-nearest neighbor coalescing random walks, the above strategy was carried out in~\cite{NRS05}.

\subsection{Convergence criteria for non-crossing paths}

When $(X_n)_{n\in\N}$ almost surely consists of paths that do not cross each other, i.e., $X_n$ contains no paths $\pi_1, \pi_2$ with $(\pi_1(s)-\pi_2(s))(\pi_1(t)-\pi_2(t))<0$ for some $s<t$, tightness in fact follows from condition (I)~\cite[Prop.~B2]{FINR04}.

\bp[Tightness criterion for non-crossing paths] \label{P:tight2}
If for each $n\in\N$, $X_n$ is an $\Hi$-valued random variable consisting almost surely of paths that do not cross each other, and $(X_n)_{n\in\N}$ satisfies condition {\rm (I)}, then $(X_n)_{n\in\N}$ is a tight family.
\ep
This result holds because when paths do not cross, the modulus of continuity of all paths in $X_n$ can be controlled by the modulus of continuity of paths in $X_n$ starting at a grid of space-time points, similar to Step (2) in the proof sketch for Theorem~\ref{T:webchar}.

When $X_n$ consists of non-crossing paths, conditions (B1') and (B2') can also be simplified to
$$
\begin{aligned}
&{\rm (B1)} \hspace{2cm}  \qquad \forall\, h>0, \ \ \limsup_{n\to\infty}   \sup_{a, t\in \R} \P\big[\eta_{X_n}(t,h; a, a+\eps)\geq 2\big]\aston{\eps} 0,\\
&{\rm (B2)} \hspace{2cm}  \qquad \forall\, h>0, \ \ \frac{1}{\eps}\limsup_{n\to\infty}   \sup_{a, t\in \R}
\P\big[\eta_{X_n}(t, h; a, a+\eps) \geq  3 \big]\aston{\eps} 0.
\end{aligned}
\hspace{3cm}
$$
Recall from the discussions before Theorem~\ref{T:webconvA} that condition (B1) is to ensure that for each deterministic $z\in \R^2$, any subsequential limit $\cal X$ of $(X_n)_{n\in\N}$ almost surely contains at most one path starting from $z$.  This property is easily seen to be implied by condition (I) (by the same argument as for Theorem~\ref{T:webchar}~(a))  when $X_n$ consists of non-crossing paths, which implies that ${\cal X}$ also consists of non-crossing paths. Therefore condition (B1) also becomes redundant, which leads to the following simplification of Theorems~\ref{T:webconvA} (cf.~\cite[Theorem 2.2]{FINR04}) and \ref{T:webconvB}.

\bt[Convergence criteria C]\label{T:webconvC} Let $(X_n)_{n\in\N}$ be a sequence of $(\Hi, \Bi_\Hi)$-valued random variables which a.s.\ consist of non-crossing paths. If $(X_n)_{n\in\N}$ satisfies conditions {\rm (I)}, and either {\rm (B2)} or {\rm (E)}, then $X_n$ converges in distribution to the standard Brownian web.
\et

Condition (B2) is often verified by applying the FKG positive correlation inequality~\cite{FKG70}, together with a bound on the distribution of the time of coalescence between two paths (see e.g.~Section~\ref{SS:srweb} below). However, FKG is a strong property that is not satisfied by most models. In such cases, verifying condition (B2) can be difficult.  Besides checking condition (E), another alternative is to use the dual (a.k.a.\ wedge) characterization of the Brownian web, as noted in Remark~\ref{R:dualchar}, to upper bound any subsequential weak limit of $(X_n)_{n\in\N}$.\footnote{Recently a new approach to verify condition (B2) was proposed in \cite{SS15}, using a Lyapunov function type criterion on the gaps between three non-crossing paths, assuming that the gaps evolve jointly as a Markov process.
}

Indeed, if $(X_n)_{n\in\N}$ consists of non-crossing paths and satisfy condition (I), then we can construct a collection of dual paths $\widehat X_n$ which almost surely do not cross paths in $X_n$, and the starting points of paths in $\widehat X_n$ become dense in $\R^2$ as $n\to\infty$. The tightness of $(X_n)_{n\in\N}$ is easily seen to imply the tightness of $(\widehat X_n)_{n\in\N}$, and any subsequential weak limit $({\cal X}, \widehat {\cal X})$ of $(X_n, \widehat X_n)$ must satisfy the property that: for any deterministic countable dense set $\Di\subset \R^2$, ${\cal X}(\Di)$ is distributed as a collection of coalescing Brownian motions, which by the non-crossing property a.s.\ uniquely determines $\widehat {\cal X}(\Di)$, which is distributed as a collection of dual coalescing Brownian motions. By the wedge characterization of the Brownian web in Remark~\ref{R:dualchar}, ${\cal X}$ equals $\overline{{\cal X}(\Di)}$, a standard Brownian web, if no path in ${\cal X}$ enters any wedge $W(\hat \pi_1, \hat \pi_2)$ of $\widehat {\cal X}(\Di)$ from outside (defined as in \eqref{wedge} with $\widehat {\cal X}(\Di)$ replacing both $\widehat \Wi^{\rm l}$ and $\widehat \Wi^{\rm r}$). This leads to
$$
{\rm (U)}  \quad
\begin{aligned}
&\mbox{For each $n\in\N$, there exists $\widehat X_n \in \widehat \Hi$ whose paths a.s. do not cross those of $X_n$ and whose} \\
&\mbox{starting points are dense in $\R^2$ as $n\to\infty$, s.t. for any subsequential weak limit $({\cal X}, \widehat {\cal X})$ of} \\
&\mbox{$(X_n, \widehat X_n)$ and any deterministic countable dense $\Di\subset\R^2$, a.s.\ paths in $\cal X$ do not enter any}\\
&\mbox{wedge of $\widehat {\cal X}(\Di)$ from outside.}
\end{aligned}
$$
We then have the following convergence result.

\bt[Convergence criteria D]\label{T:webconvD} Let $(X_n)_{n\in\N}$ be a sequence of $(\Hi, \Bi_\Hi)$-valued random variables consisting of non-crossing paths and which satisfy conditions {\rm (I)} and {\rm (U)}. Then $X_n$ converges in distribution to the standard Brownian web.
\et
\br\label{R:verifyU}
To verify condition {\rm (U)}, it suffices to show that paths in $X_n$ do not enter wedges of $\widehat X_n$ from outside, and for any deterministic $z_1, z_2\in\R^2$, not only do there exist paths $\hat \pi_{n,1}, \hat \pi_{n, 2}\in \widehat X_n$ which converge to dual coalescing Brownian motions starting at $z_1$ and $z_2$ (which follows from condition (I) for $(X_n)_{n\in\N}$ and the non-crossing between paths in $X_n$ and $\widehat X_n$), but also the time of coalescence between $\hat \pi_{n,1}$ and $\hat \pi_{n,2}$ converges to that of the coalescing Brownian motions. This ensures that the wedge $W(\hat \pi_{n,1}, \hat\pi_{n,2})$ converges and no path in the limit can enter the wedge through its bottom point. The latter can also be accomplished by showing that no limiting forward and dual paths can spend positive Lebesgue time together, as carried out in \cite{RSS15}.
\er

\subsection{Convergence of coalescing simple random walks to the Brownian web}\label{SS:srweb}

We now illustrate how the convergence criteria in Theorems~\ref{T:webconvC} and \ref{T:webconvD} for non-crossing paths can be verified for the discrete time coalescing simple random walks on $\Z$ (cf.~\cite[Theorem~6.1]{FINR04}).

Let $X$ denote the collection of discrete time coalescing simple random walk paths on $\Z$, with one walk starting from every space-time lattice site $z\in \Z^2_{\rm even}$. It is an easy exercise to show that $X$ is a.s.\ precompact in the space of paths $\Pi$, and with a slight abuse of notation, we will henceforth denote the closure of $X$ in $(\Pi, d)$ also by $X$.

For each $\eps\in (0,1)$, let $S_\eps: \R^2\to \R^2$ denote the diffusive scaling map
\be\label{Seps}
S_\eps(x,t) = (\eps x, \eps^2 t).
\ee
For a path $\pi\in \Pi$, let $S_\eps\pi$ denote the path whose graph is the image of the graph of $\pi$ under $S_\eps$.  For a set of paths $K$, define $S_\eps K:=\{ S_\eps \pi: \pi\in K\}$.

\bt\label{T:SRWweb}
Let $X$ be the collection of coalescing simple random walk paths on $\Z$ defined as above. Then as $\eps\downarrow 0$, $X_\eps:=S_\eps X$
converges in distribution to the standard Brownian web $\Wi$.
\et
{\bf Proof sketch.} We show how the various conditions in Theorems~\ref{T:webconvC} and \ref{T:webconvD} can be verified.
\medskip

\noindent
{(I)}:  This condition follows by Donsker's invariance principle. Indeed, coalescing random walks can be constructed from independent random walks by the same procedure as the inductive construction of coalescing Brownian motions from independent Brownian motions in the proof for Theorem~\ref{T:webchar}. Furthermore, this construction is a.s.\ continuous w.r.t.\ the independent Brownian motions.
Therefore (I) follows from Donsker's invariance principle for independent random walks and the Continuous Mapping Theorem.
\medskip

\noindent
{(B2)}: We will verify this condition by applying the FKG inequality~\cite{FKG70}. By translation invariance in space-time, it suffices to
show that for any $t>0$,
$$
\lim_{\delta\downarrow 0} \delta^{-1} \limsup_{\eps\downarrow 0} \P(\eta_{X_\eps}(0,t; 0,\delta)\geq 3)=0.
$$
Assume w.l.o.g.\ that $t=1$. Formulated in terms of $X$, and let $A:=\delta^{-1}$ and $\sqrt{n}=\delta \eps^{-1}$, it is equivalent to showing that
\be\label{FKG1}
\lim_{A\to\infty} A \limsup_{n\to\infty} \P(\eta_{X}(0, A^2n; 0, \sqrt{n})\geq 3)=0.
\ee
For $i\in \Z$, let $\pi_i$ denote the random walk starting at $2i$ at time $0$, and let $\tau_{i, j}$ denote the first meeting (coalescence) time between $\pi_i$ and $\pi_j$. Then by a decomposition according to the first index $k\in\N$ with $\tau_{k-1, k}>A^2n$, we have
$$
\P\big(|\{\pi_i(A^2n) : 0\leq i\leq \sqrt{n}\}|\geq 3\big) = \sum_{k=1}^{\sqrt{n}-1} \P(\tau_{0, k-1}\leq A^2n, \tau_{k-1, k}>A^2 n, \tau_{k, \sqrt{n}}>A^2n).
$$
Let us restrict to the event $E_k:=\{\tau_{0, k-1}\leq A^2n, \tau_{k-1, \sqrt{n}}>A^2n\}$ and condition on $\pi_{k-1}$ and $\pi_{\sqrt n}$. Note that the event $\{\tau_{k-1, k}>A^2n\}$ is increasing (while $\{\tau_{k, \sqrt{n}}>A^2n\}$ is decreasing) w.r.t.\ the increments of $\pi_k$,  $(\pi_k(i)-\pi_k(i-1))_{1\leq i\leq A^2n}\in \{\pm1\}^{A^2n}$, where the product space $\{\pm 1\}^{A^2n}$ is equipped with the partial order $\prec$ such that $(a_i)_{1\leq i \leq A^2n} \prec (b_i)_{1\leq i\leq A^2n}$ if $a_i\leq b_i$ for  all $i$. Furthermore, on the event $\{\tau_{k-1, k}, \tau_{k, \sqrt{n}}>A^2 n\}$, $\pi_k$ is distributed as an independent random walk with i.i.d.\ increments, whose law on $\{\pm 1\}^{A^2n}$
satisfies the FKG inequality~\cite{FKG70}. This implies that the events $\{\tau_{k-1, k}>A^2 n\}$ and $\{\tau_{k, \sqrt{n}}>A^2 n\}$ are negatively correlated under the law of $\pi_k$ with i.i.d.\ increments. Denoting $\P_k$ for probability for $\pi_k$, we then have
$$
\begin{aligned}
\sum_{k=1}^{\sqrt{n}-1} \P(\tau_{0, k-1}\leq A^2 n, \tau_{k-1, k}>A^2 n, \tau_{k, \sqrt n}>A^2 n)
& =  \sum_{k=1}^{\sqrt{n}-1} \E\big[1_{E_k} \P_k(\tau_{k-1, k}>A^2 n, \tau_{k, \sqrt n}>A^2 n)\big] \\
& \leq \sum_{k=1}^{\sqrt{n}-1} \E\big[1_{E_k}\P_k(\tau_{k-1, k}>A^2 n)\cdot \P_k(\tau_{k, \sqrt n}>A^2 n)\big]  \\
& \leq \sum_{k=1}^{\sqrt{n}-1} \P(\tau_{k-1, k}>A^2 n) \P(\tau_{k, \sqrt n}>A^2 n) \\
& \leq \sqrt{n} \frac{C}{\sqrt{A^2 n}}\cdot \frac{C\sqrt n}{\sqrt{A^2 n}}\leq \frac{C^2 }{A^2},
\end{aligned}
$$
where for the second inequality, we used the fact that $\P_{k}(\tau_{k-1, k}>A^2 n)$ and $\P_{k}(\tau_{k, \sqrt n}>A^2 n)$ are respectively
functions of $\pi_{k-1}$ and $\pi_{\sqrt n}$, which are distributed as independent random walks on the event $E_k$, and in the last line, we used that
\begin{eqnarray}
\P(\tau_{0,1}> A^2n) & \leq & \frac{C}{A\sqrt{n}},  \label{hittime} \\
\P(\tau_{k,\sqrt n}> A^2n) & \leq & \P(\exists\, k< i\leq \sqrt{n}: \tau_{i-1, i}\geq A^2n) \,  \leq\,  \frac{C \sqrt n}{A\sqrt{n}} = \frac{C}{A} , \label{hittime2}
\end{eqnarray}
which hold for any random walk on $\Z$ with finite variance~\cite[Prop.~32.4]{S76}. Condition \eqref{FKG1} then follows.
\medskip

\noindent
(U): Recall from Figure~\ref{fig:webdual} that the collection of coalescing simple random walk paths $X$ uniquely determines a collection of dual coalescing simple random walk paths, which we denote by $\widehat X$. Clearly no path in $X$ can enter
any wedge of $\widehat X$ from outside. Furthermore, for any $z_1, z_2\in \R^2$ and any choice of paths $\hat \pi_{\eps, 1}, \hat \pi_{\eps, 2}\in S_\eps \widehat X$ which converge to dual coalescing Brownian motions starting at $z_1, z_2$, it is easily seen that the time of coalescence between $\hat \pi_{\eps, 1}$ and $\hat \pi_{\eps, 2}$ also converges to that between the limiting Brownian motions. Condition (U) then follows.
\medskip

We will show how condition (E) is verified for general random walks in the next section.
\qed

\subsection{Convergence of general coalescing random walks to the Brownian web}
Let $X$ denote the collection of coalescing random walk paths on $\Z^2$ with one walk starting from each site in $\Z^2$, where the
increments are i.i.d.\ with distribution $\mu$ with zero mean and finite variance $\sigma^2 := \sum_{x\in \Z} x^2 \mu(x)$, such
that the walks are irreducible and aperiodic. Let $S^\sigma_\eps :\R^2 \to \R^2$ be the diffusive scaling map
\be\label{Seps2}
S^\sigma_\eps(x,t) = (\eps \sigma^{-1}x, \eps^2 t),
\ee
and the action of $S^\sigma_\eps$ on paths and sets of paths are defined as $S_\eps$ in \eqref{Seps}.

We have the following convergence result for general coalescing random walks on $\Z$, where paths may cross (see~\cite[Theorem~1.5]{NRS05} and \cite[Theorem~1.2]{BMSV06}).

\bt\label{T:RWweb}
Let $X$ be the collection of coalescing random walk paths on $\Z$ defined as above. If the random walk increment distribution
$\mu$ has zero mean, variance $\sigma^2$, and finite $r$-th moment $\sum_{x\in \Z} |x|^r \mu(x)<\infty$ for some $r>3$, then as $\eps\downarrow 0$, $S^\sigma_\eps X$ converges in distribution to a standard Brownian web $\Wi$.
\et
\br {\rm It was pointed out in \cite[Remark 4.1]{NRS05} that if $\sum_{x\in\Z} |x|^r\mu(x)=\infty$ for some $r<3$, then tightness of the collection of coalescing random walks is lost due to the presence of arbitrarily large jumps originating from every space-time window on the diffusive scale.  }
\er

\noindent
{\bf Proof sketch.} For $X$ consisting of random walk paths that can cross each other, it is not known how to verify condition (B2')
in Theorem~\ref{T:webconvA}. We will instead apply Theorem~\ref{T:webconvB} and sketch how the conditions therein can be verified.
Further details can be found in \cite{NRS05}.

\begin{figure}
\begin{center}
\includegraphics[width=16cm]{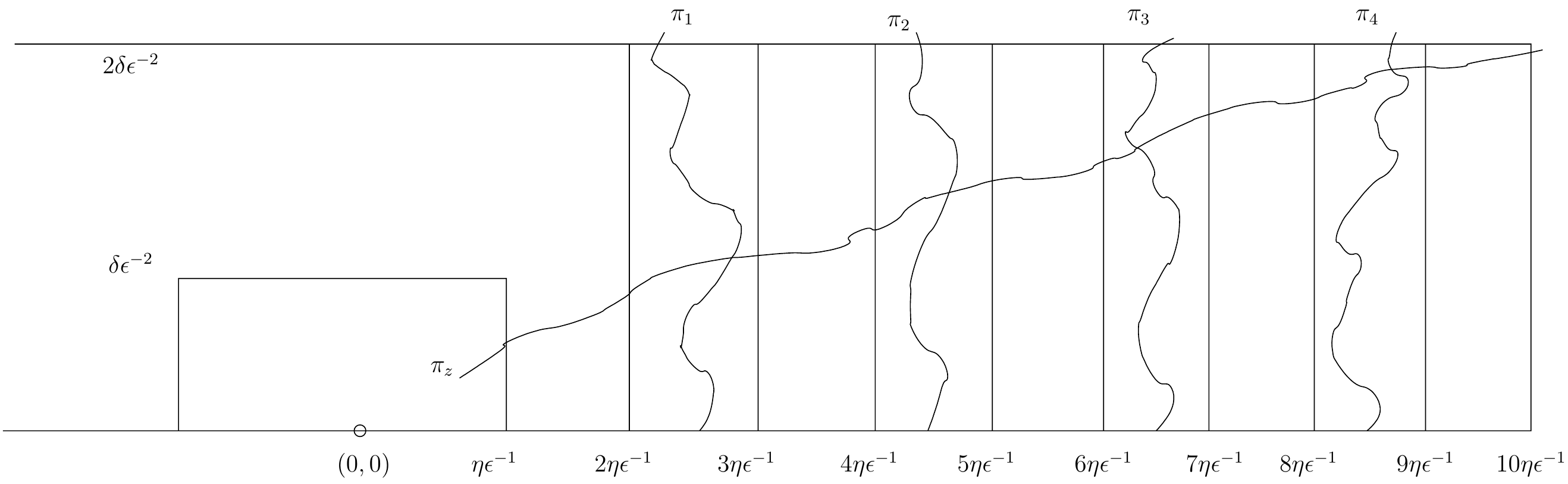}
\caption{The random walks $\pi_1, \pi_2, \pi_3, \pi_4$ start respectively at $2.5\eta\eps^{-1}$, $4.5\eta\eps^{-1}$, $6.5\eta\eps^{-1}$ and $8.5\eta\eps^{-1}$ at time $0$ and each stays within an interval of size $\eta\eps^{-1}$ up to time $2\delta\eps^{-2}$. The walk $\pi_z$ starts from $z$ inside the rectangle $[-\eta\eps^{-1}, \eta\eps^{-1}]\times [0, \delta\eps^{-2}]$ and crosses the right side of the rectangle $[-10\eta\eps^{-1}, 10\eta\eps^{-1}]\times [0, 2\delta\eps^{-2}]$ before time $2\delta\eps^{-2}$ without meeting $(\pi_i)_{1\leq i\leq 4}$.
}
\label{fig:tight2}
\end{center}
\end{figure}

To verify condition (T), by translation invariance and reformulation in terms of the set of unscaled random walk paths $X$, it suffices to show that for any $\eta>0$,
\be\label{RWTight}
\lim_{\delta\downarrow 0} \delta^{-1} \limsup_{\eps\downarrow 0} \P(X\in A_{2\delta\eps^{-2}, 20\eta\eps^{-1}}(0,0)) = 0,
\ee
where the event $E:=\{X\in A_{2\delta\eps^{-2}, 20\eta\eps^{-1}}(0,0)\}$ (see Figure~\ref{fig:tight}) is the event that $X$ contains some path which intersects the rectangle $R_1:=[-\eta\eps^{-1}, \eta\eps^{-1}]\times [0, \delta\eps^{-2}]$, and at a later time, intersects the left or right boundary of the bigger rectangle $R_2:=[-10\eta\eps^{-1}, 10\eta\eps^{-1}]\times [0, 2\delta\eps^{-2}]$. The event $E$ can occur either due to a random walk which starts outside the rectangle $R_1$ and has a jump that crosses $R_1$ horizontally, the probability of which can be easily shown to be negligible; or a random walk starts from some $z\in R_1\cap\Z^2$ and crosses the left, resp.\ the right side of $R_2$, which events we denote respectively by $E_z^-$ and $E_z^+$.

We can bound $\P(\cup_{z\in R_1\cap \Z^2} E_z^+)$ as follows. Let $(\pi_i)_{1\leq i\leq 4}$ be four random walks in $X$, starting respectively at $2.5\eta\eps^{-1}$, $4.5\eta\eps^{-1}$, $6.5\eta\eps^{-1}$ and $8.5\eta\eps^{-1}$ at time $0$. Let $B$ denote the event that each of these walks stay confined in a centered interval of size $\eta\eps^{-1}$ up to time $2\delta \eps^{-2}$ (see Figure~\ref{fig:tight2}). For $1\leq i\leq 4$, let $\tau_i$ denote the first time $n\in\N$ with
$\pi_z(n)\geq \pi_i(n)$, let $\tau_5$ denote the first time $\pi_z$ crosses the right side of $R_2$, and let $C_i$ denote the event that $\pi_z$ does not meet $\pi_i$ before time $\tau_5$. We can then bound
$$
\begin{aligned}
\lim_{\delta\downarrow 0} \delta^{-1}\limsup_{\eps\downarrow 0} \P(\cup_{z\in R_1\cap \Z^2}E_z^+) & =
\lim_{\delta\downarrow 0} \delta^{-1}\limsup_{\eps\downarrow 0} \P(B^c) + \lim_{\delta\downarrow 0} \delta^{-1}\limsup_{\eps\downarrow 0} \P(B \cap \cup_{z\in R_1\cap \Z^2}E_z^+) \\
& = \lim_{\delta\downarrow 0} \delta^{-1}\limsup_{\eps\downarrow 0} \P(B \cap \cup_{z\in R_1\cap \Z^2}E_z^+)  \\
& \leq \lim_{\delta\downarrow 0} \delta^{-1}\limsup_{\eps\downarrow 0} 2\eta\delta\eps^{-3} \max_{z\in R_1\cap \Z^2}\P(B \cap E_z^+),
\end{aligned}
$$
where in the second line we used Donsker's invariance principle and properties of Brownian motion. When the random walk increments are bounded by some $K<\infty$, it is easy to see that $\P(B \cap E_z^+)\leq C \eps^4$ uniformly in $z\in R_1\cap E_z^+$, so that the limit above equals zero.

Indeed, by successively conditioning on $\pi_z$ and $(\pi_i)_{1\leq i\leq 4}$ up to the stopping times $\tau_4$, $\tau_3$, $\tau_2$ and $\tau_1$, we note that $\pi_z$ comes within distance $2K$ of $\pi_i$ at time $\tau_i$ for each $1\leq i\leq 4$, and $\pi_z$ and $\pi_i$ must separate by a distance of at least $\eta\eps^{-1}$ before time $\tau_{i+1}$ without meeting. By the strong Markov property, these events are conditionally independent, and each event has a probability of order $\eps$ by \cite[Lemma 2.4]{NRS05}. It then follows that $\P(B\cap E_z^+)\leq C \eps^4$. When the random walks have unbounded increments, it is necessary to control the overshoot $\pi_z(\tau_i)-\pi_i(\tau_i)$, so that the probability of $\pi_z$ overshooting more than one $\pi_i$ in one jump is negligible when taken union over all starting positions $z\in R_1\cap \Z^2$. This can be done when the random walk increments have finite 5th moment~\cite{NRS05}. To relax to finite $r$-th moment for some $r>3$, a multi-scale argument is needed to take advantage of the coalescence and reduction of the random walks instead of the crude union bound as above. This was carried out in~\cite{BMSV06}.

Condition (I) can be verified by a similar argument as that for coalescing simple random walks in Theorem~\ref{T:SRWweb}. The complication is that random walk paths can cross without coalescing. However, for random walks with finite second moments, at the time when their paths cross, the distance between the two walks is of order one uniformly w.r.t.\ their starting positions, and hence it takes another time interval of order one for the two walks to coalesce. Therefore the crossing time and the coalescence time between any pair of walks are indistinguishable in the diffusive scaling limit, and forcing coalescence of the random walks at crossing times gives a good approximation, for which Donsker's invariance principle and the Continuous Mapping Theorem can be applied.

Condition (B1') amounts to showing that
\be
\lim_{\eta\downarrow 0} \P(\eta_X(0, n; 0, \eta \sqrt{n}) \geq 2) = 0,
\ee
which follows from the same bounds as in \eqref{hittime}--\eqref{hittime2}, since on the event $\{\eta_X(0, n; 0, \eta \sqrt{n})\geq 2\}$, we must have
$\tau_{i-1, i}\geq n$ for some $1\leq i\leq \eta \sqrt{n}$, where $\tau_{i-1, i}$ is the first meeting time between the two walks in $X$ starting
respectively at $i-1$ and $i$ at time $0$.

The key to verifying condition (E) is the density bound \eqref{Edensbd}, which formulated in terms of the set of unscaled random walks $X$ becomes
\be\label{Edensbd2}
\P(0\in \xi^\Z(n)) \leq \frac{C}{\sqrt n} \qquad \mbox{uniformly in } n\in\N,
\ee
where  $\xi^A(n):=\{\pi(n) : \pi \in X, \pi(0)\in A\}$ for $A\subset \Z$.  By translation invariance, for any $L\in\N$,
$$
\begin{aligned}
\P(0\in \xi^\Z(n)) & = \frac{1}{L} \E\big[\, \big| \xi^\Z(n) \cap [0, L)\big|\,\big] \\
& \leq \frac{1}{L} \sum_{k\in\Z}\E\big[\, \big| \xi^{\Z\cap [kL, (k+1)L) }(n) \cap [0, L)\big|\,\big] \\
& =  \frac{1}{L} \sum_{k\in\Z}\E\big[\, \big| \xi^{\Z\cap [0, L) }(n) \cap [kL, (k+1)L)\big|\,\big] \\
& = \frac{1}{L} \E\big[\, \big| \xi^{\Z\cap [0, L) }(n) \big|\,\big] \leq \frac{1}{L}\E\Big[1+\sum_{i=1}^{L-1} 1_{\{\tau_{i-1, i}>n\}}\Big] \leq \frac{1}{L} + \frac{C}{\sqrt n},
\end{aligned}
$$
where in the last inequality, we applied \eqref{hittime}. Letting $L\to\infty$ then gives \eqref{Edensbd2}. The rest of the proof of condition (E) then follows the line of argument sketched after Theorem~\ref{T:webconvB}.
\qed

\subsection{Convergence to the Brownian net}\label{SS:netconv}

Convergence to the Brownian net so far has only been established by Sun and Swart~\cite{SS08} for branching-coalescing simple random walk paths with asymptotically vanishing branching probability (see Figure~\ref{fig:netdual}), and recently by Etheridge, Freeman and Straulino~\cite{EFS15} for the genealogies of a spatial Lambda-Fleming-Viot process. We identify below the key conditions and formulate them as  convergence criteria for the Brownian net, which can be applied to random sets of paths with certain non-crossing properties. Finding effective and verifiable convergence criteria for random sets of paths which do not satisfy the non-crossing condition (C) below remains a major challenge.

Given a sequence of $(\Hi, \Bi_\Hi)$-valued random variables $(X_n)_{n\in\N}$, we first impose a non-crossing condition.
$$
{\rm (C)}  \quad
\begin{aligned}
&\mbox{There exist subsets of {\em non-crossing} paths } W_n^{\rm l}, W_n^{\rm r}\subset X_n, \mbox{ such that no path }\pi\in X_n \mbox{ crosses} \\
&\mbox{any } l\in W_n^{\rm l} \mbox{ from right to left, i.e.}, \pi(s)> l(s) \mbox{ and } \pi(t)<l(t) \mbox{ for some } s<t, \mbox{ and no path}\\
& \pi\in X_n \mbox{ crosses any } r\in W_n^{\rm r} \mbox{ from left to right.}
\end{aligned}
$$
The second condition is an analogue of condition (I), which ensures that $(W^{\rm l}_n, W^{\rm r}_n)_{n\in\N}$ is a tight family by Prop.~\ref{P:tight2}, and any subsequential weak limit contains a copy of the left-right Brownian web $(\Wl, \Wr)$ defined as in Theorem~\ref{T:lrweb}.
$$
{\rm (I_\Ni)}  \quad
\begin{aligned}
&\mbox{There exist $l_{n,z}\in W_n^{\rm l}$ and $r_{n,z}\in W_n^{\rm r}$ for each $z\in \R^2$, such that for any deterministic}  \\
&\mbox{$z_1, \ldots, z_k\in\R^2$, $(l_{n,z_1}, \ldots, l_{n, z_k}, r_{n, z_1}, \ldots, r_{n, z_k})$ converge in distribution to a collection} \\
&\mbox{of left-right coalescing Brownian motions starting at $(z_i)_{1\leq i\leq k}$, as in Theorem~\ref{T:lrweb}.}
\end{aligned}
$$
The tightness of $(W^{\rm l}_n, W^{\rm r}_n)_{n\in\N}$ and condition (C) imply that $(X_n)_{n\in\N}$ is also a tight family, since (C) implies that almost surely, the modulus of continuity of paths in $X_n$ (cf.~\eqref{modcon3}) can be bounded by the modulus of continuity of paths in $W^{\rm l}_n \cup W^{\rm r}_n$, whose starting points become dense in $\R^2$ as $n\to\infty$ by condition ${\rm (I_\Ni)}$.

The next condition is
$$
{\rm (H)}  \quad
\begin{aligned}
&\mbox{A.s., $X_n$ contains all paths obtained by hopping among paths in $W_n^{\rm l}\cup W_n^{\rm r}$ at crossing}\\
&\mbox{times, defined as in the hopping construction of the Brownian net in Theorem~\ref{T:netchar}}.
\end{aligned}
$$
Condition (H) ensures that any subsequential weak limit of $(X_n)_{n\in\N}$ contains not only a copy of the left-right Brownian web
$(\Wl, \Wr)$, but also a copy of the Brownian net $\Ni$ constructed by hopping among paths in $\Wl\cup\Wr$ at crossing times.

Lastly, we formulate the analogue of condition (U) in Theorem~\ref{T:webconvD}, which gives an upper bound on any subsequential weak limit of $(X_n)_{n\in\N}$ via the dual wedge characterization of the Brownian net given in Theorem~\ref{T:dualchar}.
$$
{\rm (U_\Ni)}  \quad
\begin{aligned}
&\mbox{There exist $\widehat W_n^{\rm l}, \widehat W_n^{\rm r} \in \widehat \Hi$, whose starting points are dense in $\R^2$
as $n\to\infty$, such that a.s.}\\
&\mbox{paths in $W_n^{\rm l}$ and $\widehat W_n^{\rm l}$ (resp.\  paths in $W_n^{\rm r}$ and $\widehat W_n^{\rm r}$) do not cross, and for any subsequential}\\
&\mbox{weak limit $({\cal X}, W^{\rm l}, W^{\rm r}, \widehat W^{\rm l}, \widehat W^{\rm r})$ of } (X_n, W_n^{\rm l}, W_n^{\rm r}, \widehat W^{\rm l}_n, \widehat W^{\rm r}_n) \mbox{ and any deterministic countable}\\
& \mbox{dense } \Di\subset\R^2, \mbox{a.s.\ paths in ${\cal X}$ do not enter any wedge of } (\widehat W^{\rm l}(\Di), \widehat W^{\rm r}(\Di)) \mbox{ from outside.}
\end{aligned}
$$
Condition (H) implies that $(W^{\rm l}, W^{\rm r})$ contains a copy of the left-right Brownian web $(\Wl, \Wr)$, and the non-crossing property implies that $(\widehat W^{\rm l}(\Di), \widehat W^{\rm r}(\Di))$ coincides with $(\widehat \Wi^{\rm l}(\Di), \widehat\Wi^{\rm r}(\Di)$ for the dual left-right Brownian web. By Theorem~\ref{T:dualchar}, the assumption that no path in $\cal X$ enters any wedge of $(\widehat \Wi^{\rm l}(\Di), \widehat\Wi^{\rm r}(\Di)$ from outside then implies that  $\cal X$ is contained in the Brownian net constructed from $(\Wl,\Wr)$, which is the desired upper bound on $\cal X$.

We thus have the following convergence result.

\bt[Convergence criteria for the Brownian net]\label{T:netconvD} Let $(X_n)_{n\in\N}$ be a sequence of $(\Hi, \Bi_\Hi)$-valued random variables which satisfy conditions ${\rm (C),  (I_\Ni), (H)}$ and ${\rm (U_\Ni)}$ above. Then $X_n$ converges in distribution to the standard Brownian net $\Ni$.
\et
\br\label{R:UNi}
To verify condition ${\rm (U_\Ni)}$, it suffices to show that paths in $X_n$ do not enter wedges of $(\widehat W_n^{\rm l}, \widehat W_n^{\rm r})$ from outside, and when a sequence of pairs of paths in $(\widehat W_n^{\rm l}, \widehat W_n^{\rm r})$ converge to a pair of dual left-right coalescing Brownian motions, the associated first meeting times between the pair also converge, so that the associated wedges converge.
\er

The list of conditions in Theorem~\ref{T:netconvD} can be verified for the collection of branching-coalescing simple random walks on $\Z$,
which answers {\bf Q.1} at the start of Section~\ref{S:net}.

\bt\label{T:bcwtn}
For $\eps\in (0,1)$, let $X_\eps$ be the collection of branching-coalescing simple random walk paths on $\Z^2_{\rm even}$ with branching probability $\eps$, with walks starting from every site of $\Z^2_{\rm even}$ (see Figure~\ref{fig:netdual}). Let $S_\eps$ be defined as in \eqref{Seps}.  Then $S_\eps X_\eps$ converges in distribution to the standard Brownian net $\Ni$ as $\eps\downarrow 0$.
\et
{\bf Proof sketch.}
Conditions (C) and (H) hold trivially for the branching-coalescing simple random walks. Condition ${\rm (I_\Ni)}$ was verified in \cite[Section 5]{SS08} using the fact that a single pair of leftmost and rightmost random walk paths solve a discrete analogue of the SDE \eqref{coupsde} for a pair of left-right coalescing Brownian motions, and furthermore, the time when the pair of discrete paths meet converge to the continuum analogue. By Remark~\ref{R:UNi}, it is then easily seen that condition ${\rm (U_\Ni)}$ also holds, and hence Theorem~\ref{T:netconvD} can be applied.
\qed

\section{Survey on related results}\label{S:disc}

In this section, we survey interesting results connected to the Brownian web and net that have not been discussed so far,  including alternative topologies, models whose scaling limits are connected to the Brownian web and net (including population genetic models,
true self-avoiding walks, planar aggregation, drainage networks, supercritical oriented percolation), and relation between Brownian web and net, critical planar percolation and Tsirelson's theory of noise~\cite{T04a, T04b}.

\subsection{Alternative topologies}

We review here several alternative choices of state spaces and topologies for the Brownian web, and compare them with the paths topology of Fontes et al~\cite{FINR04} introduced in Section~\ref{SS:hausdorff}. In particular, we will review the {\em weak flow topology}
of Norris and Turner~\cite{NT15}, the {\em tube topology} of Berestycki, Garban and Sen~\cite{BGS13}, and the {\em marked metric measure spaces} used by Greven, Sun and Winter~\cite{GSW15}.

Another natural extension of the paths topology  is to consider the space of compact sets of c\`adl\`ag paths equipped with the Hausdorff topology, where the space of c\`adl\`ag paths is equipped with the Skorohod metric after compactification of space-time as done in Figure~\ref{fig:comp}. Such an extension has been carried out by Etheridge, Freeman and Straulino~\cite{EFS15} in the study of scaling limits of spatial Lambda-Fleming-Viot processes.

\subsubsection{Weak flow topology} \label{SS:flow}

In \cite{NT15}, Norris and Turner formulated a topology for stochastic flows, which includes the Arratia flow generated by coalescing Brownian motions.

Recall that a flow on a space $E$ is a two-parameter family of functions
$(\phi_{s,t})_{s\leq t}$ from $E$ to $E$, which satisfies the flow condition $\phi_{t, u}\circ \phi_{s,t} = \phi_{s,u}$ for any $s\leq t\leq u$. It is easily seen that the Brownian web almost surely defines a flow on $\R$. Indeed, if $(x,s)\in \R^2$ is a $(1,2)$ point in the Brownian web with one incoming and two outgoing paths, then $(\phi_{s,t}(x))_{t\geq s}$ should be defined as the continuation of the incoming Brownian path so that the flow condition is satisfied, and otherwise, $(\phi_{s,t}(x))_{t\geq s}$ can be defined to be any of the Brownian paths starting at $(x,s)$. Note that there is no unique definition of this flow. Furthermore, it is not clear how to define a suitable topology on the space of flows in order to prove convergence of flows, as well as weak convergence of stochastic flows.

Norris and Turner addressed these issues by introducing the notion of a {\em weak flow}, which on $\R$, is a family $(\phi_{s,t})_{s\leq t}$ with the properties that
\begin{itemize}
\item[\rm (i)] For all $s\leq t$, $\phi_{s,t}\in \Di$, the space of non-decreasing functions, so that the flow lines $(\phi_{s,t}(x))_{t\geq s}$, for $(x,s)\in\R^2$, do not cross.

\item[\rm (ii)] If $\phi_{s,t}^+$ (resp.~$\phi_{s,t}^-$) is the right (resp.~left)-continuous versions of $\phi_{s,t}$, then
$$
\phi_{t,u}^- \circ \phi_{s,t}^- \leq \phi_{s,u}^- \leq \phi_{s,u}^+ \leq \phi_{t,u}^+ \circ \phi_{s,t}^+ \qquad \mbox{for all } s<t<u.
$$
\end{itemize}
One can endow the space $\Di$ with a metric $d_\Di$ such that if $f, g\in \Di$, and $f^\times, g^\times$ denote the graphs of $f$ and
$g$ rotated clockwise by $\pi/4$ around the origin, then
\be
d_\Di(f,g) = \sum_{n=1}^\infty 2^{-n} \big(1\wedge \sup_{|x|<n} |f^\times(x)-g^\times(x)|\big).
\ee
Note that when $f$ and $g$ are distribution functions, $|f^\times-g^\times|_\infty$ is just the L\'evy metric on probability measures.
The space of {\em continuous weak flows} $C^\circ(\R, \Di)$ then consists of continuous maps $\phi: \{(s,t)\in \R^2: s\leq t\} \to \Di$ with $\phi_{s,s}={\rm id}$ for all $s\in\R$, and it is endowed with the topology of uniform convergence on bounded subsets of $\{(s,t)\in\R^2: s<t\}$.

To allow for non-continuous weak flows, Norris and Turner also introduced the space of {\em c\'adl\'ag weak flows} $D^\circ (\R, \Di)$,
which consists of maps $\phi$ from the space of bounded intervals to $\Di$, such that $\phi_\emptyset ={\rm id}$, $\phi_{(s,t)}\to \phi_\emptyset$ as $t\downarrow s$ and $\phi_{(t,u)}\to \phi_\emptyset$ as $t\uparrow u$, while $\phi_{\{t\}}$ captures the jump discontinuity
of the flow at time $t$. One can then equip $D^\circ (\R, \Di)$ with a Skorohod-type topology.

It was shown in \cite{NT15} that $C^\circ (\R, \Di)$ and $D^\circ (\R, \Di)$ are Polish spaces, the Arratia flow is a continuous weak flow, and a family of Poisson local disturbance flows converge in distribution to the Arratia flow. Actually \cite{NT15} considered only weak flows on the circle $\mathbb S$. To have non-decreasing maps, a map $f: \mathbb S\to \mathbb S$ is lifted to a map from $\R$ to $\R$ by identifying $\mathbb S$ with $[0,2\pi)$ and extending $f$ outside $[0,2\pi)$ by setting $f(x+2\pi)=f(x)+2\pi$. The extension to weak flows on $\R$ has subsequently been carried out by Ellis in \cite{E10}.

Compared with the paths topology introduced in Section \ref{SS:hausdorff}, the advantage of the weak flow topology is that it is more natural for studying stochastic flows, and it allows for discontinuity in the flow lines (paths). The limitation is that it is restricted to flows with non-crossing paths. Proving weak convergence to the Arratia flow is very simple: it suffices to verify condition (I) in Section~\ref{SS:conv},
which ensures that every weak limit point contains the coalescing Brownian motions. There is no need to upper bound the limiting set of paths, unlike the paths topology of Section~\ref{SS:hausdorff}, because the weak flow topology effectively discards all paths starting from the same space-time point other than the leftmost and rightmost paths, because of the metric $d_\Di$ on $\Di$.

\subsubsection{Tube topology} \label{SS:tube}

Just as the paths topology of Fontes et al~\cite{FINR04}  was inspired by a similar topology introduced by Aizenman et al~\cite{A98, AB99}
for two-dimensional percolation, Berestycki, Garban and Sen took inspiration from another topology for two-dimensional percolation, the {\em quad topology} of Schramm and Smirnov~\cite{SS11}, and introduced in~\cite{BGS13} the {\em tube topology} for sets of continuous paths in $\R^d$. The basic idea is that the configuration of a set of paths is captured by the set of space-time tubes these paths cross.

Let $T^*:=([T^*], \partial_0 T^*, \partial_1 T^*):=([0,1]^d \times [0,1], [0,1]^d \times \{0\}, [0,1]^d\times \{1\})$ be the unit tube in $\R^{d+1}$,
with $\partial_0 T^*$ and $\partial_1 T^*$ being the lower and upper faces of the tube. A tube $T=([T], \partial_0T, \partial_1 T)$ in $\R^{d+1}$ is then defined to be the image of $T^*$ under a homeomorphism $\phi: \R^{d+1}\to \R^{d+1}$ with the property that $\partial_0T\subset \R^d \times \{t_0\}$ and $\partial_1 T \subset \R^d \times \{t_1\}$ for some $t_0<t_1$, and $[T]\subset \R^d\times [t_0, t_1]$. We call $t_0$ the bottom time of $T$, and $t_1$ the top time. The space of all tubes, denoted by ${\cal T}$, can then be equipped with the Hausdorff metric $d_{\cal T}$, with
\be
d_{\cal T}(T_1, T_2) := d_{\rm Haus}([T_1], [T_2]) + d_{\rm Haus}(\partial_0 T_1, \partial_0 T_2) + d_{\rm Haus}(\partial_1 T_1, \partial_1 T_2),
\ee
where $d_{\rm Haus}$ is the Hausdorff distance between subsets of $\R^{d+1}$. The metric space $({\cal T}, d_{\cal T})$ is then separable.

Given a continuous path $\pi: [t, \infty)\to\R^d$ starting at time $t$, it is said to cross a tube $T$ with lower face $\partial_0 T\subset \R^d\times \{t_0\}$ and upper face $\partial_1 T\subset \R^d\times \{t_1\}$ for some $t_0<t_1$, if $t\leq t_0$, $(\pi(t_0), t_0)\in \partial_0 T$,
$(\pi(t_1), t_1)\in \partial_1 T)$, and $(\pi(s), s)\in [T]$ for all $s\in [t_0, t_1]$. Given a set $K$ of continuous paths in $\R^d$, one can then identify the set of all tubes ${\rm Cr}(K)\subset \cal T$ which are crossed by some path in $K$. Furthermore, if $K$ is a compact set of continuous path w.r.t.\ a metric on path space defined in the same way as in Section~\ref{SS:hausdorff}, then ${\rm Cr}(K)$ is in fact a closed subset of $\Ti$. Therefore a random compact set of paths can be identified with a random closed subset of ${\cal T}$. The state
space of closed subsets of the metric space ${\cal T}$ can then be equipped with the Fell topology, which makes it compact (see e.g.~\cite[Appendix B]{M05}).

One can actually narrow down the state space further. Observe that we can define a partial order $\leq$ on $\cal T$, where we denote $T_1\leq T_2$ if whenever $T_2$ is crossed by some path $\pi$, $T_1$ is also crossed by $\pi$. We denote $T_1<T_2$ if there are open neighbourhoods $U_i\subset\Ti$ around $T_i$, $i=1,2$, such that $T_1'\leq T_2'$ for all $T_1'\in U_1$ and $T_2'\in U_2$. A set of tubes $S\subset \cal T$ is then called {\em hereditary} if $T\in S$ implies that $T'\in S$ for all $T'<T$. Note that the set of tubes ${\rm Cr}(K)$ induced by a set of paths $K$ is always hereditary. Therefore the state space for random compact sets of paths can be taken to be the space of {\em closed hereditary subsets} of  ${\cal T}$, denoted by ${\cal H}$, equipped with the Fell topology. It can be shown that ${\cal H}$ is closed
under the Fell topology, and hence compact. This gives the tube topology defined in \cite[Section 2]{BGS13}, and the Brownian web can
be realised as a random variable taking values in ${\cal H}$.

The main advantages of the tube topology are: (i) Tightness for a family of ${\cal H}$-valued random variables comes for free because the state space ${\cal H}$ is compact. For instance, under the tube topology, Theorem~\ref{T:RWweb} on the convergence of general coalescing random walks to the Brownian web can be established under the optimal finite second moment assumption, since the higher moment assumption is only used to establish tightness under the path topology. (ii) The tube topology, being actually a weaker topology than the path topology of Fontes et al~\cite{FINR04} introduced in Section~\ref{SS:hausdorff}, makes it much easier to construct coalescing flows which do not satisfy the non-crossing property of the Arratia flow. In particular, the coalescing Brownian flow on the Sierpinski gasket was constructed in \cite{BGS13} using the tube topology, and an invariance principle was established.

To characterize a probability measure on  ${\cal H}$, it turns out to be sufficient to determine the probability of the joint crossing of any finite collection of tubes chosen from a  deterministic countable dense subset of ${\cal T}$. To prove the weak convergence of a sequence of probability measures on $\cal H$, it is sufficient to find a large enough set of tubes $\hat {\cal T}\subset \cal T$ such that the the probability of the joint crossing of any finite subset of $\hat {\cal T}$ converges to a limit, see \cite[Prop.~2.12]{BGS13}. The strategy for verifying this convergence criterion for coalescing flows is similar to
the proof of condition (E) in Theorem~\ref{T:webconvB}. First one approximates a given tube $T$ with bottom time $t_0$ by tubes $T_\delta$ with bottom time $t_0+\delta$ for $\delta>0$, then one applies a coming down from infinity result to show that among paths which
intersect the lower face of $T$, only finitely many remain in the lower face of $T_\delta$ at time $t_0+\delta$, and lastly one uses condition (I) to control the joint distribution of the remaining finite collection of coalescing paths.

The tube topology is weaker than the paths topology introduced in Section~\ref{SS:hausdorff}, because the mapping from the space of compact sets of continuous paths (with the paths topology) to the space of closed hereditary sets of tubes ${\cal H}$ (with the Fell topology) is continuous as shown in \cite[Lemma A.2]{BGS13}. In particular, given a set of paths $K$, restricting a path $\pi\in K$ to
a later starting time and adding it to $K$ has no effect under the tube topology. In a sense, the tube topology is also insensitive toward the behavior of the paths near their starting times. In particular, a sequence of paths $\pi_n$, which start at time $0$ and converge to a path $\pi$ uniformly on $[\delta, \infty)$ for any $\delta>0$, will converge to $\pi$ under the tube topology. But $\pi_n$ may have wild oscillations in the interval $[0, 1/n)$ which prevent it from converging in the paths topology.

\subsubsection{Marked metric measure spaces}\label{SS:mmm}

In~\cite{GSW15}, Greven, Sun and Winter treats the (dual) Brownian web as a stochastic process taking values in the space of spatially-marked metric measure spaces. The notion of a $V$-marked metric probability measure space was introduced by Depperschmidt, Greven and Pfaelhuber in \cite{DGP11}, which is simply a complete separable metric space $(X,d)$, together with a Borel probability measure $\mu$ on the product space $X\times V$ for some complete separable metric space of marks $V$. The probability measure $\mu$ can be regarded as a sampling measure, and each $V$-marked metric measure ($V$-mmm) space can be uniquely determined (up to isomorphism) by the joint distribution of $(v_i)_{i\in\N}$ and the distance matrix $(d(x_i, x_j))_{i,j\in\N}$, where $(x_i, v_i)_{i\in\N}$ are i.i.d.\ samples drawn from $X\times V$ with common distribution $\mu$. A sequence of $V$-mmm spaces is then said to converge if the associated random vector of marks $(v_i)_{i\in\N}$ and the random distance matrix $(d(x_i, x_j))_{i,j\in\N}$ converge in finite-dimensional distribution. By truncations in the mark space $V$, one can also extend the notion of $V$-mmm spaces to the case where $\mu$ is only required to be finite on bounded sets when projected from $X\times V$ to $V$, which was done in \cite{GSW15}.

To see how the Brownian web fits in the framework of $V$-mmm spaces, note that the coalescing Brownian motions in the (dual) Brownian web can be regarded as the space-time genealogies of a family of individuals. More precisely, if $\widehat\Wi$ is the dual Brownian web introduced in Section~\ref{SS:webdual}, then for each $z=(x,t)\in\R^2$, interpreted as an individual at position $x$ at time $t$, $\hat\pi_z\in \widehat\Wi$ determines its spatial genealogy, with $\hat\pi_z(s)$ (for $s\leq t$) being the spatial location of its ancestor at time $s$. Coalescence of two paths $\hat\pi_{z_1}$ and $\hat\pi_{z_2}$ then signify the merging of the two genealogy lines. If we consider all individuals indexed by $\R$ at a given time $t$, then we can measure the genealogical distance between individuals, i.e.,
\be\label{ultrametric}
d((x,t), (y,t)) :=  2(t-\hat\tau_{(x,t), (y,t)}),
\ee
where $\hat\tau_{(x,t), (y,t)}$ is the time of coalescence between $\hat\pi_{(x,t)}$ and $\hat\pi_{(y,t)}\in \widehat\Wi$. It is easily seen that $d$ is in fact an ultra-metric, i.e.,
$$
d((x_1,t), (x_2, t)) \leq d((x_1,t), (x_3, t)) \vee d((x_2,t), (x_3, t)) \qquad \mbox{for all } x_1, x_2, x_3\in \R.
$$
The collection of individuals indexed by $\R$ at time $t$ then form a metric space with metric $d$. The mark of an individual is just its spatial index, so that we can identify the metric space with the mark space $V:=\R$, and a natural sampling measure is the Lebesgue measure on $\R$. There is one complication, namely that a.s.\ there exists $x\in \R$ such that multiple paths in $\widehat\Wi$ start from $(x,t)$. The correct interpretation is that each such path encodes the genealogy of a distinct individual, which happen to occupy the same location. Therefore the space of population should be enriched from $\R$ to take into account these individuals. Such an enrichment does not affect the Lebesgue sampling measure, since it follows from Theorem~\ref{T:classweb} that a.s.\ there are only countably many points of multiplicity in $\widehat\Wi$ at a given time $t$. If we let $(X_t, d_t, \mu_t)$ denote the $\R$-mmm space induced by $\widehat\Wi$ on $\R$ at time $t$, then $\widehat\Wi$ a.s.\ determines the stochastic process $(X_t, d_t, \mu_t)_{t\in\R}$, and conversely, it can be seen
that $(X_t, d_t, \mu_t)_{t\in\R}$ a.s.\ determines $\widehat\Wi$, and hence the two can be identified.

We can also identify $\widehat\Wi$ with a single $\R^2$-mmm space instead of the $\R$-mmm space-valued process $(X_t, d_t, \mu_t)_{t\in\R}$. Namely, we can consider all points in $\R^2$ and extend the genealogical distance in \eqref{ultrametric} between individuals at the same time to individuals at different times:
\be\label{ultrametric2}
d((x_1,t_1), (x_2,t_2)) :=  t_1+t_2-2\hat\tau_{(x_1,t_1), (x_2,t_2)},
\ee
where $\hat\tau_{(x_1,t_1), (x_2,t_2)}$ is the time of coalescence between $\hat\pi_{(x_1, t_1)}$ and $\hat\pi_{(x_2, t_2)}\in \widehat\Wi$. We then obtain a metric space whose elements can be identified with the mark space $V:=\R^2$ (with suitable enrichment of $\R^2$ to take into account points of multiplicity in $\widehat\Wi$), on which we equip as a sampling measure the Lebesgue measure on $\R^2$.
We remark that one can also enlarge the mark space and let the whole genealogy line $\hat\pi_z\in\widehat\Wi$ be the mark for a point $z\in\R^2$.

For models arising from population genetics, $V$-mmm space is a natural space for the spatial genealogies, and the dual Brownian web determines the genealogies of the so-called {\em continuum-sites stepping-stone model}. Because $V$-mmm spaces are characterized via sampling, proving convergence essentially reduces to condition (I) formulated in Section~\ref{SS:conv}. However, proving tightness for $V$-mmm space-valued random variables is a non-trivial task and can be quite involved (see e.g.~\cite{GSW15}).
Also the Brownian net cannot be characterized using $V$-mmm spaces, since there is no natural analogue of the metric $d$ in \eqref{ultrametric2} when paths can branch.

\subsection{Other models which converge to the Brownian web and net}

We review here various models which have scaling limits that are connected to the Brownian web and net. These include population genetic models such as the voter model and the spatial Fleming-Viot processes, stochastic Potts models, the true self-avoiding walks, planar aggregation models, drainage networks, and supercritical oriented percolation in dimension $1+1$.

\subsubsection{Voter model and spatial Fleming-Viot processes}\label{SS:voter}

Voter model and spatial Fleming-Viot processes are prototypical population genetic models where the spatial genealogies of the population are coalescing random walks, which converge to the Brownian web under diffusive scaling.

\begin{figure}[t]
\begin{subfigure}{.5\textwidth}
  \centering
  \includegraphics[width=.9\linewidth]{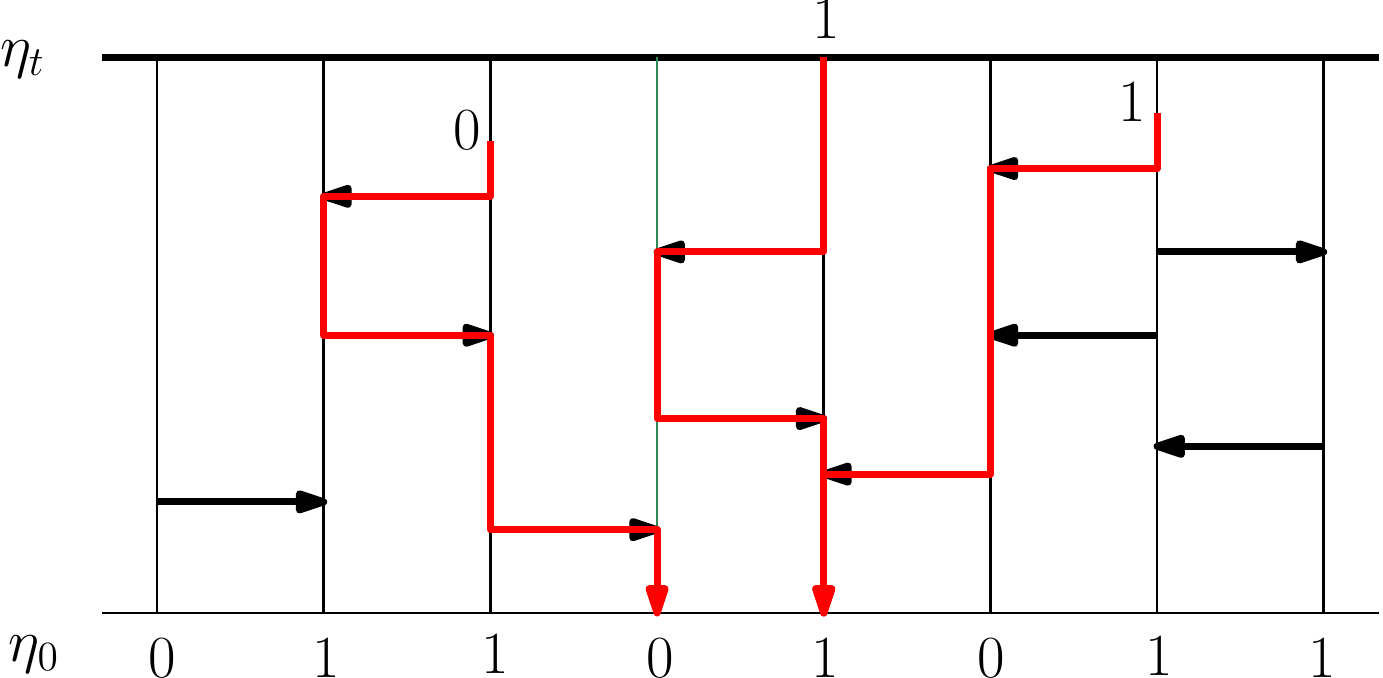}
  \caption{}
  \label{fig:voter}
\end{subfigure}%
\begin{subfigure}{.5\textwidth}
  \centering
  \includegraphics[width=.9\linewidth]{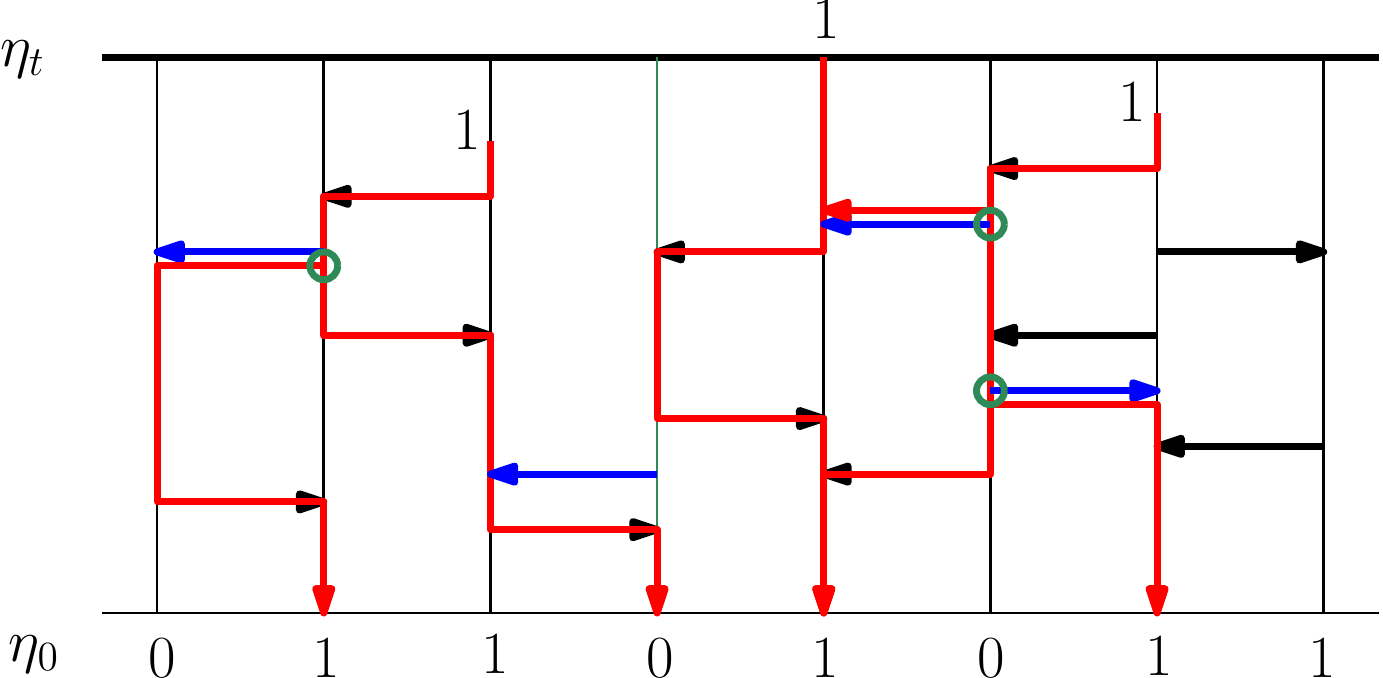}
  \caption{}
  \label{fig:bvoter}
\end{subfigure}
\caption{Harris' graphical construction of: (a) the voter model; (b) the biased voter model.}
\label{fig:harris}
\end{figure}

Arratia~\cite{A79, A81} first conceived the Brownian web in studying the scaling limit of the voter model on $\Z$, which is an interacting
particle system $(\eta_t)_{t\geq 0}$ with $\eta_t \in \Omega:=\{0,1\}^\Z$, modelling the opinions of a collection of individuals indexed by $\Z$. Independently for each $x\in\Z$, at exponential rate $1$, a resampling event occurs where the voter at $x$ picks one of its two neighbors with equal probability and changes its opinion $\eta_t(x)$ to that of the chosen neighbor. The resampling events can be represented by a Poisson point process of arrows along the edges over time, with each arrow from $x$ to $y$ at time $t$ signifying the voter at $x$ changing its opinion at time $t$ to that of $y$ (see Figure~\ref{fig:voter}). This is known as Harris' graphical construction. To identify $\eta_t(x)$, one just needs to trace the genealogy of where the opinion $\eta_t(x)$ comes from backward in time, i.e., follow the arrows backward in time, until an ancestor $y$ is reached at time $0$ so that $\eta_t(x)=\eta_0(y)$. The genealogy line $\hat\pi_{(x,t)}$
for $\eta_t(x)$ is then a continuous time random walk running backward in time. Furthermore, for multiple space-time points $((x_i,t_i))_{1\leq i\leq k}$, their joint genealogy lines $(\hat\pi_{(x_i, t_i)})_{1\leq i\leq k}$ is a collection of coalescing random walks. This
is known as the {\em duality} between the voter model and coalescing random walks (see \cite{HL75, L04} for more details). Taking diffusive scaling limit of the joint genealogies of all voter opinions at all possible times then leads to what we now call the (dual) Brownian web.

The spatial Fleming-Viot process $(\xi_t)_{t\geq 0}$ is a measure-valued process which extends the voter model by allowing a continuum
of individuals at each site $x\in \Z$, represented by a probability measure on the type space $[0,1]$. Individuals in the population migrate
on $\Z$ as independent random walks, and resampling (one individual changes its type to that of another) takes place between every pair of individuals at the same site with exponential rate $1$. This model has recently been studied in \cite{GSW15}, and the joint spatial genealogies are a collection of coalescing random walks running backward in time, where each pair of walks at the same site coalesce with exponential rate $1$. It is not difficult to see that the diffusive scaling limit of these genealogies also gives the (dual) Brownian web.

It is natural to ask whether the (dual) Brownian web also determines the genealogies of a population genetic model on $\R$. The answer is affirmative, and such a continuum model has been studied before and known as the {\em continuum-sites stepping-stone model} (CSSM), see e.g.\ Donnelly et al~\cite{DEF+00}. The Brownian web effectively gives a Harris graphical construction of the CSSM. Convergence of the voter model as a measure-valued process on $\R\times \{0,1\}$ to the CSSM has been established in \cite{AS11},
and convergence of the genealogies of the spatial Fleming-Viot process to that of the CSSM as stochastic processes taking values in the space of $\R$-mmm spaces (see Section~\ref{SS:mmm}) has been established in \cite{GSW15}. In both results, the use of the
sampling measure allows one to establish convergence under the optimal finite second moment assumption on the increments of the
underlying coalescing random walks, in contrast to Theorem~\ref{T:RWweb}.

We remark that coalescing random walks on $\Z$ also model the evolution of boundaries between domains of different spins in the $0$-temperature limit of the stochastic Potts model, and the Brownian web has been used to study aging in this model~\cite{FINS01}.

\subsubsection{Biased voter model}\label{SS:bvoter}

We have seen that the voter model is dual to coalescing random walks. It turns out that branching-coalescing random walks is dual to the {\em biased voter-model}, also known as Williams-Bjerknes model~\cite{WB71, S77}, which modifies the voter model by adding a selective bias so that type $1$ is favoured over type $0$. More precisely, each voter independently undergoes a second type of resampling event with rate $\eps$, where the voter at $x$ changes its type to that of the chosen neighbor $y$ only if $y$ is of type $1$. The Harris graphical construction can then be modified by adding a second independent Poisson point process of (selection) arrows along the edges over time, where a selection arrow is used only if it points to an individual of type $1$ at that time (see Figure~\ref{fig:bvoter}). To determine $\eta_t(x)$, we then trace its genealogy backward in time by following the resampling arrows, and when a selection arrow is encountered, we follow both potential genealogies (by either ignoring or following the selection arrow). The potential genealogies then form a collection of branching-coalescing random walks with branching rate $\eps$. It is not difficult to see that $\eta_t(x)=1$ if and only if $\eta_0(y)=1$ for at least one ancestor $y$ that can be reached by the potential genealogies. This establishes the duality between the biased voter model and branching-coalescing random walks. If we consider a sequence of biased voter models with selection rate $\eps\downarrow 0$, while rescaling space-time by $(\eps, \eps^2)$, then their genealogies converge to the Brownian net.

We remark that branching-coalescing random walks on $\Z$ also model evolving boundaries between domains of different spins in the low temperature limit of the stochastic Potts model, where a new domain of spins can nucleate at the boundary of two existing domains~\cite{MNR13}.

\subsubsection{True self-avoiding walks and true self-repelling motion}

Arratia first conceived what we now call the Brownian web in his unfinished manuscript~\cite{A81}, and the subject lay dormant until T\'oth and Werner~\cite{TW98} discovered a surprising connection between the Brownian web and the {\em true self-avoiding walk} on $\Z$ with bond repulsion~\cite{T95}. A special case of the true self-avoiding walk is defined as follows. Let $l_0(\cdot)$ be an integer-valued function defined on the edges of $\Z$, which can be regarded as the initial condition for the edge local time of the walk, i.e., how many times each edge has been traversed. Given the walk's position $X_n=x$ and the edge local time $l_n(\cdot)$ at time $n$, if $l_n(\{x-1, x\})<l_n(\{x, x+1\})$, then with probability $1$, $X_{n+1}=x-1$; if $l_n(\{x, x+1\})<l_n(\{x-1, x\})$, then with probability $1$, $X_{n+1}=x+1$; and if $l_n(\{x-1, x\})=l_n(\{x, x+1\})$, then $X_{n+1}=x\pm 1$ with probability $1/2$ each. The edge local time $l_n(\cdot)$ is then updated to $l_{n+1}(\cdot)$ by adding $1$ to the local time at the newly traversed edge. Such a walk is called a true self-avoding walk because it is repelled from the more visited regions, and the laws of $(X_i)_{1\leq i\leq n}$ are consistent as $n$ varies, in contrast to the {\em self-avoiding walk}.

\begin{figure}
\begin{center}
\includegraphics[width=16cm]{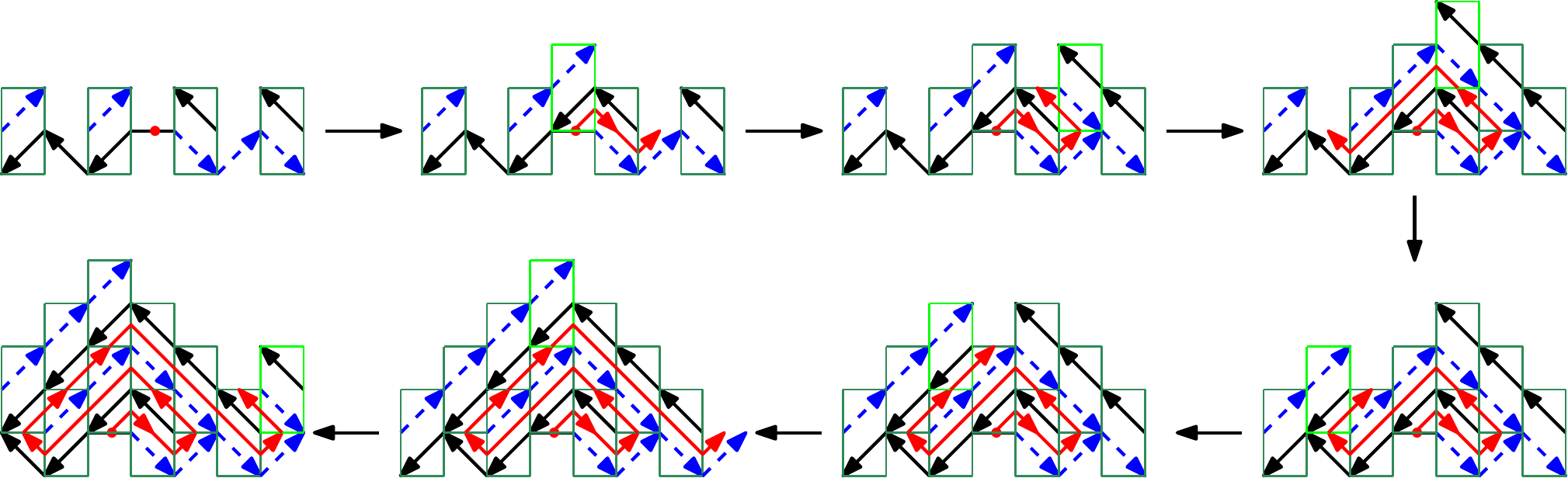}
\caption{Coupling between a true-self avoiding walk (with a special initial edge local time) with a collection of forward/backward coalescing random walks (with a special boundary condition). The horizontal coordinate of the lattice-filling curve is the walk's position and the vertical coordinate is its edge local time.}
\label{fig:TSW}
\end{center}
\end{figure}

Interestingly, if we plot the evolution of the position of the walk together with its edge local time, then there is an almost sure coupling with a
collection of forward/backward coalescing random walks. Figure~\ref{fig:TSW} illustrates such a coupling, where the initial edge local time is given by $l_0(\{2x, 2x+1\})=1$ for all $x\in\N_0$, $l_0(\{2x-1, 2x\})=0$ for all $x\in\N$, and $l_0(\{-x-1, -x\})= l_0(\{x, x+1\})$ for all $x\in \Z$. Such an $l_0$ corresponds to a special boundary condition for the coalescing random walks along the $x$-axis, as shown in Figure~\ref{fig:TSW}. The lattice-filling curve between the forward and backward coalescing random walks encodes the evolution of the position of the walk (the horizontal coordinate) and the edge local time (the vertical coordinate), and the area filled in by the curve is just the time the walk has spent. The collection of forward/backward coalescing random walks converge to the forward/backward Brownian web, albeit with the boundary condition that the path in the forward, resp.\ backward web starting at $(0,0)$ is the constant path, and time for the webs now runs in the horizontal direction in Figure~\ref{fig:TSW}. It can then be shown (see \cite{NR06}) that  the lattice-filling curve squeezed between the coalescing random walks also converge to a space-filling curve $(X_t, L_t(X_t))_{t\geq 0}$ in $\R\times [0,\infty)$, squeezed between paths in the forward and backward Brownian webs. Here, $X_t$ is the position of the so-called {\em true self-repelling motion}, and $(L_t(x))_{x\in\R}$ is its occupation time density at time $t$. Comparing with the discrete model, one sees that when a deterministic point $(x,y)\in \R\times [0,\infty)$ is first reached by the space-filling curve $(X_\tau, L_\tau(X_\tau))$ at some random time $\tau$, $(L_\tau(a))_{a\leq x}$, resp.\ $(L_\tau(a))_{a\geq x}$, must equal the path in the forward, resp.\ backward, Brownian web starting at $(x,y)$. Furthermore,
the area under $(L_\tau(a))_{a\in\R}$ must be exactly equal to $\tau$. In other words, from the realization of the forward/backward Brownian webs, one can determine when a deterministic point $(x,y)$ is reached by the space-filling curve. By considering a deterministic countable dense set of $(x,y)\in \R\times [0,\infty)$, one can then determine the entire trajectory of $(X_t)_{t\geq 0}$. This is how the true self-repelling motion was constructed by T\'oth and Werner in \cite{TW98}, which heuristically is a process with a drift given by the negative of the gradient of its occupation time density. It has the unusual scaling invariance of $(X_{at}/a^{2/3})_{t\geq 0} \stackrel{\rm dist}{=} (X_t)_{t\geq 0}$ for any $a>0$, and it has finite variation of order $3/2$. For further details, see \cite{TW98, D12}.

\subsubsection{Planar aggregation models}

\begin{figure}
\begin{center}
\includegraphics[width=15cm]{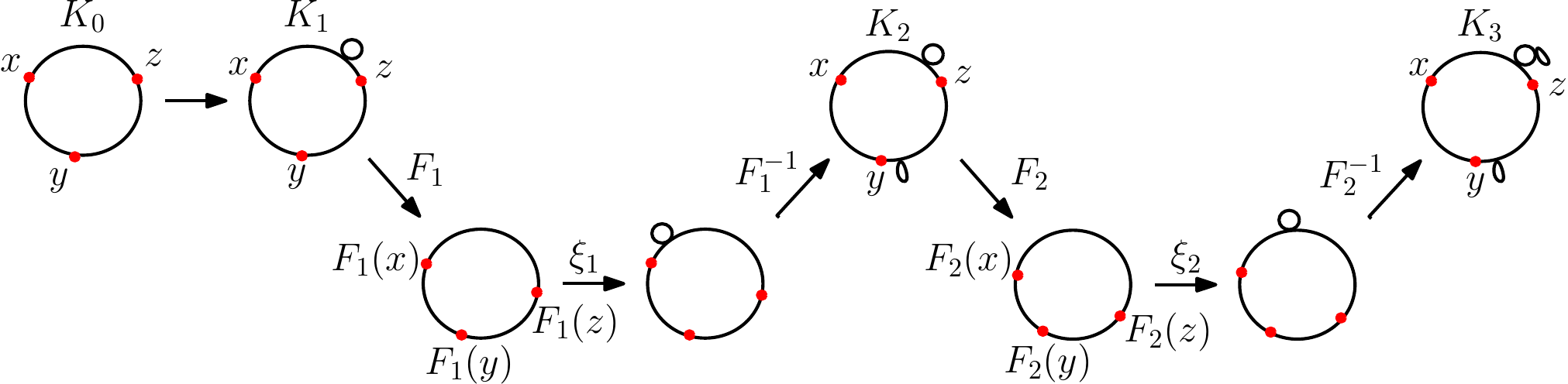}
\caption{Hastings-Levitov planar aggregation model.}
\label{fig:agg}
\end{center}
\end{figure}

In \cite{NT12}, Norris and Turner discovered an interesting connection between the Hastings-Levitov planar aggregation model and the coalescing Brownian flow. We briefly explain the model and the connection here. At time $0$, we start with a unit disk denoted by $K_0$. A small particle $P_1$, which we assume to be a disk of radius $\delta>0$ for simplicity, is attached to $K_0$ at a uniformly chosen random point on the boundary $\partial K_0$. This defines the new aggregate $K_1=K_0\cup P_1$. To define the next aggregate $K_2$, we apply a conformal map $F_1$ which maps $K_1$ back to the disk $K_0$ (more accurately, $F_1$ is a conformal map from $K_1^c\cup \{\infty\}$ to $K_0^c\cup \{\infty\}$, uniquely determined by the condition that $F_1(z)=Cz +O(1)$ for some $C>0$ as $|z|\to\infty$). A new particle $P_2$ of radius $\delta$ is then attached randomly at the boundary of $K_0=F_1(K_1)$. Reversing $F_1$ then defines the new aggregate $K_2= F_1^{-1}(F_1(K_1)\cup P_2)$. The dynamics can then be iterated as illustrated in Figure~\ref{fig:agg}. The new particles added to the aggregates $(K_n)_{n\in\N}$ at each step are no longer disks due to distortion by the conformal maps, and it was shown in~\cite{NT12} that new particles tend to pile on top of each other and form protruding fingers. Interestingly, the image of $\partial K_0$ under $F_n$ forms a coalescing flow on the circle, where for $x,y\in \partial K_0\subset \partial K_n$, the length of the arc between $F_n(x)$ and $F_n(y)$ on the unit circle is proportional to the probability that a new particle will be attached to the corresponding part of $\partial K_n$ between $x$ and $y$. In the limit that the particle radius $\delta\downarrow 0$, while time is sped up by a factor of $\delta^{-3}$, this coalescing flow can be seen as a localized disturbance flow studied in \cite{NT15}, which converges to the coalescing Brownian flow on the circle w.r.t.\ the weak flow topology described in Section~\ref{SS:flow}.

\subsubsection{Drainage networks and directed forests}

\begin{figure}
\begin{center}
\includegraphics[width=16cm]{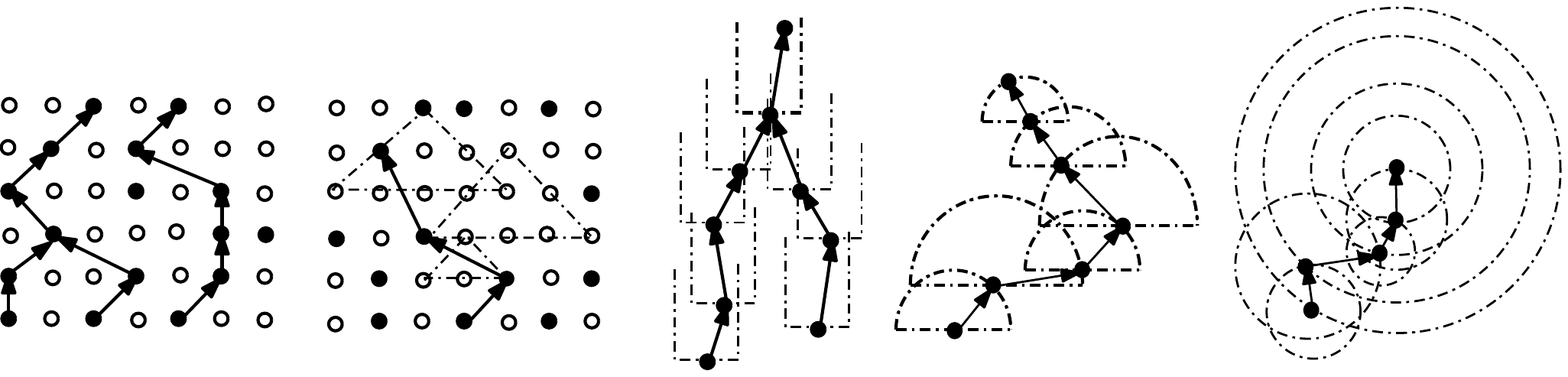}
\caption{From left to right, illustration of paths in the drainage network of \cite{GRS04},  the drainage network of \cite{RSS13}, Poisson trees \cite{FLT04}, the directed spanning forest \cite{BB07}, and the radial spanning tree \cite{BB07}.}
\label{fig:drainage}
\end{center}
\end{figure}

Drainage networks are a class of models where coalescing paths arise naturally. First, a random subset of $\Z^{d+1}$ (or $\R^{d+1}$) is determined, which represents the water sources. Next, from each source $(x,s)$, where $x\in\R^d$ and $s\in\R$, exactly one directed edge is drawn toward some other source $(y, t)$ with $t>s$, representing the flow of water from $(x,s)$ to $(y,t)$. Examples include \cite{GRS04, RSS15}, where the authors study a drainage network on $\Z^{d+1}$ with each vertex being a water source independently with probability $p\in (0,1)$. From each source $(x,s)$, a directed edge is then drawn to the closest source in the next layer $\Z^d\times \{s+1\}$, and ties are broken by choosing each closest source with equal probability. There are also other variants such as in \cite{ARS08}, where the directed edge from the source $(x,s)\in \Z^{d+1}$ connects to the closest source in the $45^o$ light cone rooted at $(x,s)$, or as in \cite{RSS13} where the directed edge from $(x,s)$ connects to the closest source in $\Z^{d}\times \{s, s+1, \ldots\}$ measured in $\ell_1$ distance (see Figure~\ref{fig:drainage}).

There are also continuum space versions where the water sources form a homogeneous Poisson point process in $\R^{d+1}$. In the {\em Poisson trees} model considered in \cite{FLT04}, a directed edge is drawn from each source $(x,s)\in\R^{d+1}$ to the source in $\{(y,t): |y-x|\leq 1, t>s\}$ with the smallest $t$-coordinate. In the {\em directed spanning forest} model considered in \cite{BB07}, a directed edge is drawn from each source $(x,s)$ to the closest source in $\R^d \times (s,\infty)$, measured in Euclidean distance.

Natural questions for such drainage networks include whether the directed edges form a single component, i.e., it is a tree rather than a forest. This has been shown to be the case for $d=1$ and $2$ for the drainage networks on $\Z^{d+1}$ described above~\cite{GRS04, ARS08, RSS13} and for the Poisson trees~\cite{FLT04}, as well as for the directed spanning forest in $\R^{1+1}$~\cite{CT13}.

The collection of paths in the drainage network, obtained by following the directed edges, can be regarded as a collection of (dependent) coalescing random walk paths. Because the dependence is in some sense local in all the models described above, it is natural to conjecture that for $d=1$, the collection of directed paths in the drainage network (after diffusive scaling) should converge to the Brownian web. This has indeed been verified for various drainage networks on $\Z^{1+1}$~\cite{CDF09, CV11, RSS13}, and for the Poisson trees on
$\R^{1+1}$~\cite{FFW05}. The main difficulty in these studies lies in the dependence among the paths. The model considered in \cite{RSS13} is a discrete analogue of the directed spanning forest~\cite{BB07}, and the dependence is handled using simultaneous regeneration times along multiple paths. Such arguments are not directly applicable to the directed spanning forest, and its convergence to the Brownian web remains open.

The directed spanning forest was introduced in \cite{BB07} as a tool to study the {\em radial spanning tree}, where given a homogeneous Poisson point process $\Xi\subset \R^d$, a directed edge is drawn from each $x\in \Xi$ to the closets point
$y\in \Xi\cup\{\rm o\}$ that lies in the ball with radius $|x|$ and centered at the origin o. If one considers directed paths in the radial spanning forest restricted to the region (in polar coordinates) $\{(r, \theta): r\in [an, n], |\theta+\pi/2|\leq n^{-b}\}$ for some $a\in (0,1)$ and $b\in (0, 1/2)$, then after proper scaling, these paths are believed to converge to the so-called {\em Brownian bridge web}~\cite{FVV14}, which consists of coalescing Brownian bridges starting from $\R\times (-\infty, 0)$ and ending at $(0,0)$ (it can also be obtained from the Brownian web via a deterministic transformation of $\R^2$). For a couple of toy models, such a convergence has been established in \cite{CV14, FVV14}.

\subsubsection{Supercritical oriented percolation}\label{SS:supperc}

\begin{figure}
\begin{center}
\includegraphics[width=4.5cm]{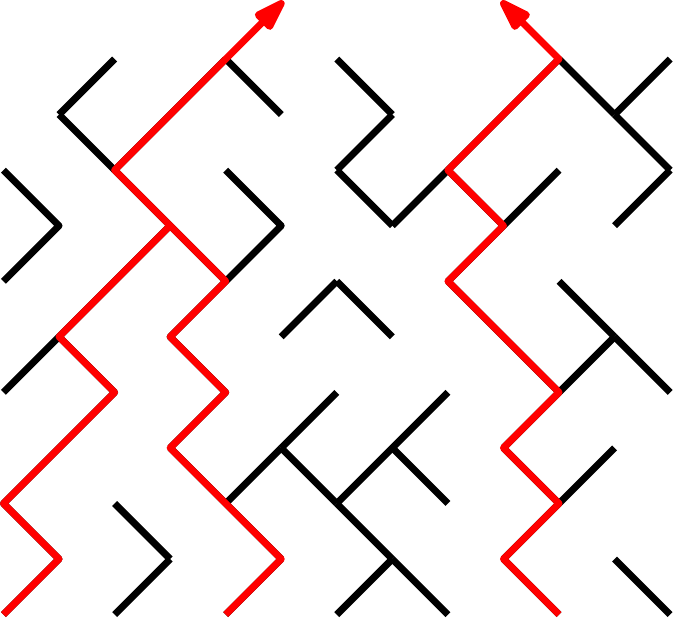}
\caption{Rightmost infinite open paths in supercritical oriented percolation on $\Z^2_{\rm even}$.}
\label{fig:OP}
\end{center}
\end{figure}

The oriented bond percolation model on $\Z^2_{\rm even}:=\{(x,t)\in \Z^2: x+t \mbox{ is even}\}$ is defined by independently setting each
oriented edge of the form $(x,t) \to (x\pm1, t+1)$ to be either open with probability $p$, or closed with probability $1-p$, with retention parameter $p\in [0,1]$. The set of vertices in $\Z^2_{\rm even}$ that can be reached from $z\in\Z^2_{\rm even}$ by following open oriented edges is called the open cluster at $z$. It is known that there is a critical $p_c\in (0,1)$ such that for $p>p_c$, the open cluster at $(0,0)$ is infinite with positive probability, while for $p\leq p_c$, it is finite with probability $1$~\cite{D84}. When the open cluster at $z\in \Z^2_{\rm even}$ is infinite, $z$ is called a percolation point.

In the supercritical regime $p>p_c$, the percolation points appear with a positive density,
and starting from each percolation point $z=(x,t)$, we can find a rightmost path $r_z$ among all the infinite open oriented paths starting from $z$ (see Figure~\ref{fig:OP}). It was shown by Durrett~\cite{D84} that each $(r_z(n))_{n\geq t}$ satisfies a law of large numbers with drift $\alpha(p)>0$, and subsequently Kuczek~\cite{K89} showed that $r_z(n)-n\alpha$ satisfies a central limit theorem with variance $n \sigma(p)^2$ for some $\sigma(p)>0$. Kuczek's CLT can be easily extended to path level convergence to Brownian motion. It is then natural to ask what is the joint scaling limit of $r_z$ for all percolation points $z$. Wu and Zhang~\cite{WZ08} conjectured that the scaling limit should be the Brownian web. This may be surprising at first, because for different percolation points $z_1$ and $z_2$, $r_{z_1}$ and $r_{z_2}$ both depend on the infinite future, and hence could be strongly dependent on each other. However, it was shown by Sarkar and Sun in \cite{SS13} that each
$r_z$ can be approximated by a thin percolation exploration cluster, such that different clusters evolve independently before they meet and quickly merge after they meet. It is then shown in~\cite{SS13} that after proper centering and scaling, the collection of rightmost infinity open paths from percolation points indeed converge to the Brownian web.

\subsection{Brownian web, critical percolation, and noise}\label{SS:webnoise}

There are close parallels between the Brownian web and the scaling limit of critical planar percolation, starting from the topology. Indeed,
as noted in Section \ref{SS:tube}, both the paths topology of \cite{FINR04} and the tube topology of \cite{BGS13} were inspired by similar topologies for planar percolation. The Brownian web is the scaling limit of discrete coalescing random walks on $\Z^2_{\rm even}$ as shown in Figure~\ref{fig:webdual}, which can be seen as a dependent planar percolation model and enjoys the same duality as bond percolation on $\Z^2_{\rm even}$. The scaling limit of critical planar percolation is invariant under conformal maps, while the Brownian web is invariant under diffusive scaling of $\R^2$, a signature of a model at criticality. The Brownian web and the scaling limit of critical planar percolation are both (and the only known examples of) two-dimensional black noise in the language of Tsirelson~\cite{T04a, T04b}, which is intimately linked to the notion of noise sensitivity and the existence of near-critical scaling limits. Indeed, the Brownian net can be regarded as a near-critical scaling limit obtained by perturbing the Brownian web. As explained in Section~\ref{S:couple}, such a perturbation can be carried out by Poisson marking the set of ``pivotal'' points (i.e., the $(1,2)$-points of the Brownian web) according to a natural local time measure, and then turning these points into branching points. For planar percolation, exactly the same procedure has been proposed in \cite{CFN06} to construct near-critical scaling limit from the critical scaling limit, and it has been rigorously carried out in
\cite{GPS13a} and \cite{GPS13b}, where defining the natural local time measure on pivotal points alone has been a very challenging task.
The fact that changing the configuration at a countable subset of $\R^2$, i.e., the Poisson marked pivotal points, is sufficient to alter the scaling limit indicates noise sensitivity and that the scaling limit is a black noise. Allowing the intensity of the Poisson marking of the pivotal points to vary continuously leads to dynamical evolutions of the critical scaling limit, where unusual behaviour may emerge at random ``dynamical times''.

In what follows, we will briefly review Tsirelson's theory of noise~\cite{T04a, T04b}, and in what sense are the Brownian web and the scaling limit of critical planar percolation both black noise. This is the key common feature that lies behind the many parallels between the two models. Such parallels will likely extend to other examples of two or higher-dimensional black noise, once they are found. We will also briefly review results on the dynamical Brownian web.

\subsubsection{Brownian web, black noise, and noise sensitivity}

Intuitively, the Brownian web and the scaling limit of critical planar percolation should satisfy the property that, the configurations are independent when restricted to disjoint domains in $\R^2$. This calls for the continuous analogue of the notion of product probability spaces (or independent random variables), which led Tsirelson to define (see \cite[Definition 3c1]{T04b})

\bd
A $($one-dimensional$)$ {\em continuous product of probability spaces}  consists of a probability space $(\Omega, \Fi, P)$ and
sub-$\sigma$-fields $\Fi_{s,t}\subset \Fi$ $($for all $s<t$$)$ such that $\Fi$ is generated by $\bigcup_{s<t}\Fi_{s,t}$ and for any $r<s<t$,
\be\label{noisefactor}
\Fi_{r,t} = \Fi_{r,s}\otimes \Fi_{s,t},
\ee
which means that $\Fi_{r,t}$ is generated by $\Fi_{r,s} \cup \Fi_{s,t}$, and $\Fi_{r,s}$ is independent of $\Fi_{s,t}$ $($i.e., $P(A\cap B)=P(A)P(B)$ for all $A\in \Fi_{r,s}$ and $B\in \Fi_{s,t}$$)$.
\ed
We can think of $\Fi_{s,t}$ as the $\sigma$-field generated by a family of observables that only depend on what happens in the interval $(s,t)$.
\br\label{R:factorize}
In dimension two or higher, we should equip $(\Omega, \Fi, P)$ with a family of sub-$\sigma$-fields $\Fi_D$, indexed by a Boolean algebra of domains $D\subset \R^2$ (take for example the Boolean algebra of sets generated by open rectangles), while the factorization property
\eqref{noisefactor} becomes
\be\label{noisefactor2}
\Fi_{D} = \Fi_{D_1} \otimes \Fi_{D_2}
\ee
whenever $D_1\cap D_2=\emptyset$ and $\overline{D} = \overline{D_1\cup D_2}$. It turns out that \eqref{noisefactor2} cannot hold without
some restrictions on the domains. The lack of a canonical choice of the family of domains in dimensions two and higher is still an issue to be resolved~\cite[Sec.~1.6]{T14}.
\er

Besides independence on disjoint domains, we also expect the Brownian web and the scaling limit of percolation to be translation invariant.
This additional assumption leads to the notion of a noise (see \cite[Definition 3d1]{T04b}).

\bd[Noise]
A noise is a continuous product of probability spaces $(\Omega, \Fi, P)$, equipped with a family of sub-$\sigma$-fields $\Fi_D$ indexed by
an algebra of domains $D\in \Di$ in $\R^d$, which is homogeneous in the following sense.  There exists a group of isomorphisms $(T_h)_{h\in\R^d}$ on $(\Omega, \Fi, P)$ such that $T_{h+h'}=T_{h}\circ T_{h'}$ and $T_h$ sends $\Fi_D$ to $\Fi_{D+h}$ for all $D\in \Di$ and $h, h'\in\R^d$.
\ed

In a tour de force, the scaling limit of critical planar percolation was shown to be a two-dimensional noise by Schramm and Smirnov in \cite{SS11}, where the sub-$\sigma$-field $\Fi_D$ (for $D$ with piecewise smooth boundary)  is generated by the indicator random variables for open crossing of quads (homeomorphic images of the unit square) contained in $D$. The Brownian web was shown to be a two-dimensional noise by Ellis and Feldheim in \cite{EF12}, where $\Fi_D$ (for open rectangles $D$) is generated by the Brownian web paths restricted to $D$. What is remarkable is that both noises are so-called {\em black noise} as defined by Tsirelson~\cite{T04a, T04b}, and they are the only known examples in dimension two or higher.
\medskip

The notion of a noise being {\em black} turns out to be equivalent to the notion of a noise being {\em sensitive}, while a noise being {\em classical} (such as white noise or Poisson noise) turns out to be equivalent to a noise being {\em stable}. A {\em non-classical} noise is a noise that contains in some sense a sensitive part. To  explain the underlying ideas, it is instructive to first discretize and then pass to the continuum, instead of giving directly the definition for the continuum noise.

Let us fix $\Lambda=(-1, 1)^2$, and consider its discretization $\Lambda_\delta := \Lambda \cap S_\delta \Z^2_{\rm even}$, where
$S_\delta :(x,t)\to (\delta x, \delta^2 t)$. To each $z\in \Lambda_\delta$, we assign an i.i.d.\ symmetric random variable $\sigma_z\in \{+1, -1\}$, which can be defined via the coordinate map on the probability space $\Omega_\delta=\{+1, -1\}^{\Lambda_\delta}$, equipped with the discrete topology and the uniform probability measure $P_\delta$. What noise we obtain in the continuum limit $\delta\downarrow 0$ depends crucially on the observables we choose, which will generate the $\sigma$-field $\Fi$ for the limiting noise.

If we choose our observables to be linear functions of $\sigma:=(\sigma_z)_{z\in \Lambda_\delta}$, i.e., $f_\delta:=\delta^{3/2}\sum_{z\in \Lambda_\delta} f(z) \sigma_z$ for continuous functions $f:\Lambda \to \R$, then their joint distributions converge to a non-trivial limit and the limiting noise is the classical white noise, where $\Fi$ is generated by the family of Gaussian random variables indexed by such continuous $f$, which we interpret as the integral of $f$ w.r.t.\ the underlying white noise.

If we choose our observables to be non-linear functions of $(\sigma_z)_{z\in \Lambda_\delta}$, then non-classical noise may appear in the limit. As illustrated in Figure~\ref{fig:webdual}, $\sigma:=(\sigma_z)_{z\in \Lambda_\delta}$ uniquely determines the collection of coalescing random walks in $\Lambda_\delta$. If we choose our observables to be the indicator random variables for whether there is a random walk
path crossing a prescribed tube $T$ in $\Lambda$ (see Section~\ref{SS:tube}), then their joint distributions converge to a non-trivial limit and the limiting noise is the Brownian web, where $\Fi$ is generated by the tube crossing events.

Roughly speaking, whether the limiting noise is classical or not depends on whether there are non-trivial observables that are noise sensitive as $\delta\downarrow 0$, defined as follows.
\bd[Noise sensitivity]
For $\eps\in (0,1)$, let $\sigma^\eps:=(\sigma_z^\eps)_{z\in \Lambda_\delta}$ be obtained from $(\sigma_z)_{z\in \Lambda_\delta}$ by independently replacing each $\sigma_z$ with an independent copy of $\sigma_z$ with probability $\eps$. A sequence of random variables $f_\delta:\Omega_\delta\to \R$ with $\sup_\delta E[f_\delta^2(\sigma)]<\infty$ is called noise sensitive if for each $\eps>0$,
\be\label{sens}
\lim_{\delta\downarrow 0}\big(E[f_\delta(\sigma) f_\delta(\sigma^\eps)] - E[f_\delta(\sigma)]^2\big) = 0,
\ee
and the sequence is called noise stable if
\be\label{stab}
\lim_{\eps\downarrow 0} \limsup_{\delta\downarrow 0}\big|\, E[f_\delta(\sigma) f_\delta(\sigma^\eps)] - E[f^2_\delta(\sigma)]\,\big| = 0.
\ee
\ed
There is a rich theory of noise sensitivity for functions on $(\Omega_\delta, P_\delta)$, and we refer to the lecture notes by Garban and Steif~\cite{GS12}  for a detailed exposition and applications to percolation. For $\{0,1\}$-valued functions $f_\delta$, noise sensitivity implies that $f_\delta(\sigma)$
and $f_\delta(\sigma^\eps)$ become asymptotically independent if an arbitrarily small, but fixed portion of $\sigma$ is resampled, while noise stability implies that as $\eps\downarrow 0$, the probability that $f_\delta(\sigma)$ and $f_\delta(\sigma^\eps)$ coincide tends to $1$
uniformly in $\delta$ close to $0$.

There is a simple criterion for noise sensitivity/stability using Fourier analysis on $(\Omega_\delta, P_\delta)$. Observe that
the set of functions $\chi_S:=\prod_{i\in S}\sigma_i$, for all $S\subset \Lambda_\delta$, is an orthonormal basis for $L^2(\Omega_\delta, P_\delta)$. Therefore any $f_\delta\in  L^2(\Omega_\delta, P_\delta)$ admits the orthogonal decomposition
\be\label{FourierWalsh}
f_\delta = \sum_{S\subset \Lambda_\delta} \hat f_\delta(S) \chi_S = E[f_\delta] + \sum_{k=1}^{\infty} \sum_{|S|=k}  \hat f_\delta(S) \prod_{i\in S} \sigma_i.
\ee
\bd[Spectral measure and energy spectrum]
The coefficients $\hat f_\delta(S)$ in \eqref{FourierWalsh} are called the {\em Fourier-Walsh} coefficients for $f_\delta$, the measure $\mu_{f_\delta}$ on $\{S: S\subset \Lambda_\delta\}$ with $\mu_{f_\delta}(S) := \hat f_\delta^2(S)$ is called the {\em spectral measure} of $f_\delta$, and the measure $E_{f_\delta}$
on $\N_0$ with $E_{f_\delta}(k) = \sum_{|S|=k} \hat f_\delta^2(S)$ is called the {\em energy spectrum} of $f_\delta$ (see e.g.~\cite{GS12}).
\ed
It is easily seen that
\be
E[f_\delta(\sigma) f_\delta(\sigma^\eps)] = E[f_\delta]^2+ \sum_{k=1}^{\infty} \eps^k \sum_{|S|=k}  \hat f_\delta^2(S) = E[f_\delta]^2+ \sum_{k=1}^{\infty} \eps^k E_{f_\delta}(k),
\ee
and hence a sequence $f_\delta$ with $\sup_\delta E[f_\delta^2]<\infty$ is noise sensitive if and only if $E_{f_\delta}(\{1,\ldots, k\})\to 0$
for every $k\in\N$, and $f_\delta$ is stable if and only if $E_{f_\delta}$ is a tight family of measures on $\N_0$.

As we take the continuum limit $\delta\downarrow 0$, what will happen to the spectral measure of the observables that generate the $\sigma$-field of the limiting noise? There are several possibilities: either all non-trivial square-integrable observables will in the limit have a spectral measure that is supported on $\Ci_{\rm fin}(\Lambda):= \{S\in \Lambda: |S|<\infty\}$, in which case the limiting noise is called {\em classical} (and noise stable); or all non-trivial observables will in the limit have a spectral measure supported on $\Ci(\Lambda)\backslash \Ci_{\rm fin}(\Lambda)$, with $\Ci(\Lambda)$ being the set of closed subsets of $\Lambda$, in which case the noise is called {\em black} (and noise sensitive); or there are non-trivial observables of both types, supported either on $\Ci_{\rm fin}(\Lambda)$ or $\Ci(\Lambda)\backslash \Ci_{\rm fin}(\Lambda)$, in which case the noise is called {\em non-classical}.

In practice, classifying a noise via discrete approximation as described above is subtle, because such approximations are not unique, and convergence of the spectral measure on $\Ci$ (equipped with Hausdorfff topology) does not imply convergence of the energy spectrum.  There could be observables $f$ for the noise whose spectral measure is supported on $\Ci_{\rm fin}$, and yet the approximating $f_\delta$ is noise sensitive (see \cite[Sec.~5c]{T04a} for the notion of {\em block sensitivity} which overcomes such issues). Instead, we can classify noise directly. The following definition extends \cite[Def.~3d2]{T04a} to general dimensions (see also \cite{T14} for an abstract formulation without associating the sub-$\sigma$-fields with domains in $\R^d$).

\bd[Non-classical noise and black noise] Let $(\Omega, \Fi, P)$ be a noise equipped with a family of sub-$\sigma$-fields $\Fi_D$ indexed by an algebra of domains $D\in \Di$ in $\R^d$. For any $f\in L^2(\Omega, \Fi, P)$, its spectral measure $\mu_f$ is defined to be the unique
positive measure on $\Ci$ $($the space of compact subsets of $\R^d$$)$ with the property that
\be
\mu_f(M\in \Ci: M\subset D) = E[ E[f \, |\, \Fi_D]^2], \qquad D\in \Di.
\ee
If for every $f\in L^2(\Omega, \Fi, P)$, its spectral measure $\mu_f$ is supported on $\Ci_{\rm fin}:=\{M\in \Ci: |M|<\infty\}$, then the noise is called {\em classical}; otherwise the noise is called {\em non-classical}. If the spectral measure of every non-constant $f\in L^2(\Omega, \Fi, P)$ is supported on $\Ci\backslash \Ci_{\rm fin}$, then the noise is called {\em black}.
\ed
White noise and Poisson noise are classical, which can be seen from their chaos expansions. The scaling limit of critical planar percolation is a two-dimensional black noise as explained in \cite{SS11}, while the Brownian web was first shown to be a one-dimensional black noise in \cite[Section 7]{T04a} and \cite[Section 7]{T04b} (with sub-$\sigma$-fields $(\Fi_{s,t})_{s<t}$ associated with the strips $\R\times (s,t)$), and then shown to be a two-dimensional black noise in \cite{EF12}. Interestingly, if a noise is black, then for all non-constant $f\in L^2(\Omega, \Fi, P)$, $\mu_f$ is in fact supported on $\Ci_{\rm perf}:=\{M\in \Ci: M \mbox{ is a perfect set}\}$, a collection of uncountable sets~\cite[Theorem 6d3]{T04a}.

The procedure of independently resampling each of a collection of $\pm1$-valued symmetric random variables with probability $\eps$
can in fact be carried out for noise. This allows one to define directly noise sensitivity vs stability for each $f\in L^2(\Omega, \Fi, P)$, where
conditions \eqref{sens}--\eqref{stab} just need to be modified to remove the $\delta$-dependence, and $\sigma$ and $\sigma^\eps$ represent the noise and its $\eps$-resampled version. A noise is then classical if and only if all $f\in L^2(\Omega, \Fi, P)$ are stable,
and black if and only if all non-constant $f\in L^2(\Omega, \Fi, P)$ are sensitive (see \cite[Sec.~4]{T04b} and \cite[Sec.~5]{T04a} for further details). In particular, $\eps$-resampling a black noise leads to an independent copy of the noise. To see non-trivial dependence, a different resampling procedure is needed, which is where dynamical evolution via Poisson marking of pivotal points comes in.

\subsubsection{Dynamical Brownian web}

Because the Brownian web is a black noise, the correct way to resample the randomness underlying the Brownian web in order to see non-trivial dependence, is to resample the ``pivotal'' points (i.e., the $(1,2)$-points of the Brownian web).
This was carried out in \cite{NRS10} and is intimately linked to the authors' marking construction of the Brownian net in the same paper
(see Section~\ref{S:couple}). Given a Poisson point process $M_\lambda\subset\R^2$ with intensity measure given by $\lambda (>0)$ times the intersection local time measure $\ell({\rm d}z)$ on the $(1,2)$-points of the Brownian web (see Prop.~\ref{P:refloc}), a resampled
web is obtained by simply flipping the sign of each $(1,2)$ point in $M_\lambda$ as defined in \eqref{sgnw}, instead of turning each such point into a branching point as in the construction of the Brownian net in Theorem~\ref{T:coupling}. Since there is a natural coupling between $M_\lambda$ for all $\lambda \geq 0$, we obtain a stochastic process in the space of compact sets of paths $(\Hi, d_\Hi)$,
which is in fact reversible w.r.t.\ the law of the Brownian web. This defines the so-called {\em dynamical Brownian web}. The dynamical Brownian web can also be constructed as the scaling limit of the {\em dynamical discrete web} introduced in \cite{HW09b}, where {\em discrete web} refers to the collection of coalescing simple symmetric random walks on $\Z^2_{\rm even}$ illustrated in Figure~\ref{fig:webdual}, and dynamical evolution is introduced by independently flipping at rate $1$ the sign of the arrow from each $z\in \Z^2_{\rm even}$. If space-time is rescaled by $(\eps, \eps^2)$, and the dynamical time is slowed down by a factor of $\eps$, it can then be shown that the dynamical discrete web converges to the dynamical Brownian web~\cite{NRS10}.

Similar questions have been studied for percolation models~\cite{G99}. There has been much study about dynamical percolation (see \cite{S09} for a survey). Furthermore, dynamical critical site percolation on the triangular lattice has been shown to converge to a continuum limit~\cite{GPS13b}, the proof of which follows the same line as in the Brownian web setting (via Poisson marking of the pivotal points), except that the implementation in the percolation setting has been much more difficult.

Key questions in the study of dynamical percolation, dynamical Brownian web, and ``dynamical'' processes in general, are the existence of exceptional times when a given property that holds almost surely at a deterministic time fails. For dynamical percolation, where the stationary law of the process is that of the percolation model at criticality, one such property is the non-existence of an infinite cluster~\cite{S09}.  For the dynamical discrete web, questions on exceptional times have also been studied \cite{FNRS09} and \cite{J11}. In particular, it was shown in~\cite{FNRS09} that there exists a.s.\ a random set of dynamical times with full Lebesgue measure, when the random walk in the discrete web starting from the origin does not satisfy the Law of the Iterated Logarithm (LIL); the authors also extended this result to the dynamical Brownian web.

This is in striking contrast with the dynamical simple symmetric random walk on $\Z$ introduced in \cite{BHPS03}, where the sign of each increment of the walk $(X_n-X_{n-1})_{n\in\N}$ is independently flipped at rate 1 in dynamical time. In contrast with the dynamical discrete web, there is no exceptional time at which the LIL fails in this setting. The difference between the two models can be accounted for by the effect of each individual switching. In \cite{BHPS03}, a switch affects only one increment of the walk, which induces an order $1$ perturbation in sup-norm. In contrast, changing the direction of a single arrow in the discrete web can change the path of the walk by a ``macroscopic'' amount. Indeed, the switching arrows encountered by the walk are comparable to the Poisson marked $(1,2)$-points encountered by a Brownian motion in the dynamical Brownian web, where each switching at a marked $(1,2)$-point leads the Brownian motion to embark on an excursion from the original path, and the size of the excursion has a heavy-tailed distribution with infinite mean.
This is why the random walk is more sensitive to the dynamics on the discrete web than to the dynamics on the random walk considered in \cite{BHPS03}.

\section{Open questions}\label{S:open}

In this section, we discuss some open questions, including how the Brownian web and net could be useful
in the study of more general one-dimensional interacting particle systems, conjectures on the geometry of the Brownian net,
and some other miscellaneous questions.

\subsection{Voter model perturbations and Brownian net with killing}\label{SS:perturb}

One further direction of research is to use the Brownian web and net to study one-dimensional
interacting particle systems with migration, birth/branching, coalescence, and death. Such a
particle system often arises as the spatial genealogies of a dual particle system. In suitable
parameter regimes, such a particle system may converge to the Brownian net, which incorporates
diffusive particle motion, branching and coalescence, and killing effect can be added by defining
the so-called {\em Brownian net with killing}, which we will explain below. Using the convergence of
the discrete particle systems to a continuum model which is amenable to analysis,  one hopes to
draw conclusions for the discrete particle systems, such as local survival vs extinction, as well
as convergence to equilibrium.

In dimensions $3$ and higher, a study in this spirit has been carried out by Cox, Durrett and Perkins in \cite{CDP13}, where they
studied a class of interacting particle systems on $\Z^d$ whose transition rates
can be written as perturbations of the transition rates of a (possibly non-nearest neighbor) voter model.
Their interest stemmed from the observation that several models in ecology and evolution, as well as in statistical mechanics --- the spatial Lotka-Volterra model \cite{NP99}, the evolution of cooperation \cite{OHLN06} ,  non-linear voter models \cite{MDDGL99} --- can be seen as perturbations of the voter model in a certain range of their parameter space. In general, these models are difficult
to study either because of a lack of monotonicity (Lotka-Volterra model)
or because of the intrinsic complexity of the model (certain models of evolution of cooperation).
However, by considering them as small perturbations of the voter model, it was shown in \cite{CDP13} that for $d\geq3$,
the properly rescaled local density of particles converges to the solution of a reaction diffusion equation, and that properties
of the underlying particle system (such as coexistence of species in the Lotka-Volterra model) can be derived from the behavior
of this PDE. Implementing an analogous program in dimensions $1$ and $2$ is still an open problem, although in dimension $2$, it is believed that the behavior of voter model perturbations should be similar to the case $d=3$.

The genealogies of the voter model is given by coalescing random walks~\cite{L04, HL75}, and the type of perturbations of the voter model in \cite{CDP13} lead to branching of the genealogies. In dimensions $d\geq 3$, because of transience of the random walk, such a collection of branching-coalescing genealogies was shown in \cite{CDP13} to converge (under suitable scaling and moment assumptions on the random walk and the branching mechanism) to a system of {\it independent} branching Brownian motions. The branching mechanism is given explicitly in terms of the ``microscopic'' descriptions of the system and also takes into account coalescence at the mesoscopic level. In dimension $1$, such a collection of branching-coalescing spatial genealogies should converge to the Brownian net under suitable scaling of parameters and space-time. One is then hopeful that in dimension $1$, one can also draw conclusions for the particle systems from such a convergence. One model that falls in this class is the biased annihilating branching process, introduced by Neuhauser and Sudbury in \cite{NS93}. Another class is a family of one-dimensional models of competition with selection and mutation (see e.g., \cite{AS13}). However, proving convergence of the particle systems to the Brownian net is a serious challenge.



\subsubsection{Brownian net with killing}\label{SS:netkill}

We briefly recall here the Brownian net with killing, which allows the modelling of the death of a particle.  A discrete space analogue has been introduced in \cite{MNR13}. Let $b,\kappa\geq0$ be such that $b+\kappa\leq1$, corresponding respectively to the branching and killing parameters of the system. From each point in $\Z^2_{\rm even}$, an arrow
is drawn from $(x,t)$ to either $(x-1,t+1)$ or $(x+1,t+1)$ with equal probability $\frac{1}{2}(1-b-\kappa)$; two arrows are drawn
with probability $b$; no arrow is drawn with probability $\kappa$. Similar to the discrete net (see Fig.~\ref{fig:couple}),
we can consider the set of all paths starting in $\Z^2_{\rm even}$ by following arrows until the path terminates.
This defines an infinite collection of coalescing random walks which branch with probability $b$ and are killed with probability $\kappa$.
This model encompasses several classical models from statistical mechanics.  For $p\in(0,1)$
and $b=p^2$, $\kappa=(1-p)^2$, one recovers the standard
one dimensional oriented percolation model. When $b=\kappa=0$,
the trajectories are distributed as coalescing random walks.

In \cite{NRS15}, it is shown that when the branching and killing parameters
depend on some small scaling parameter $\eps$
such that $\lim_{\eps\rightarrow0} b_\eps / \eps = b_\infty$ and $\lim_{\eps\rightarrow 0} \kappa_\eps/ \eps^2=\kappa_\infty$,
then the collection of branching-coalescing random walks with
killing introduced above converges (after proper rescaling)
to a continuum object called the Brownian net with killing with parameter $(b_\infty, \kappa_\infty)$.
The latter scaling of the parameters was motivated by the stochastics $q$-color Potts model on $\Z$.
Indeed, in \cite{NRS15}, it was conjectured that
this model at large inverse temperature $\beta$, is asymptotically dual
to a system of (continuous time)
branching-coalescing random walks with killing,
with branching and killing parameters
$$
b_\beta = \frac{q}{2} e^{-\beta} \quad \mbox{and} \quad \kappa_\beta = q e^{-2\beta},
$$
where $q$ is the number of colors in the system.
Thus the Brownian net could be used
to construct a natural scaling limit for the stochastic Potts model
at low temperature in dimension $1$.

The Brownian net with killing can also be directly constructed as follows. Let ${\mathcal N}_b$
be a Brownian net with branching parameter $b$. For every realization of the Brownian net,
one can introduce a natural time-length measure on the set of points of type (p,p) in $\R^2$ w.r.t.\ the
Brownian net (see Theorem~\ref{T:netpts}). More precisely, for every Borel set $E\subset\R^2$,
define the time-length measure as
\be
\cL(E) = \int_\R |\{x : (x,t)\in E \mbox{ and is of type (p,p)}\}|\ {\rm d}t \label{ex-cl}.
\ee
The Brownian net with killing
with parameter $(b,\kappa)$ is then obtained by killing the Brownian net paths
at the points of a Poisson point process on $\R^2$ with intensity measure $\kappa \cL$.

It was shown in \cite{NRS15} that the Brownian net with killing undergoes a phase transition as the killing parameter $\kappa$ varies, similar to the percolation transition in oriented percolation. There is reason to believe~\cite{H00} that the Brownian net with killing belongs to the same universality class as oriented percolation in $\Z^{1+1}$, in which case it may serve as a simpler model to study this universality class.

\subsection{Fractal structure of the Brownian net}\label{SS:geometry}

We explain here some open questions concerning the geometry of the Brownian net.  Recall from Remark~\ref{R:backbone} the set $\Ni(*,-\infty)$ of paths started at time $-\infty$, called the \emph{backbone} of the Brownian net $\Ni$. Let $N:=\{(x,t)\in\R^2:x=\pi(t)\mbox{ for some }\pi\in\Ni(*,-\infty)\}$ be the graph of the backbone. By the image set property (Prop.~\ref{P:imgset}), every continuous function $\pi:\R\to\R$ that lies entirely inside $N$ is a path in the backbone. In view of this, much information can be read off from the random closed set $N\sub\R^2$.

If we set $\xi_t:=\{x\in\R:(x,t)\in N\}$, then $(\xi_t)_{t\in\R}$ is just the branching-coalescing point set in stationarity  (see Section~\ref{S:bcptset}), which is in fact a Markov process~\cite[Sec.~1.9--1.10]{SS08}. At each deterministic time $t$, $\xi_t$ is distributed a Poisson point process with intensity 2, which is a reversible law for this Markov process. As noted in Remark~\ref{R:backbone}, by reversibility, the law of $N\subset \R^2$ (and likewise the collection of paths $\Ni(*,-\infty)$) is symmetric under time reversal.

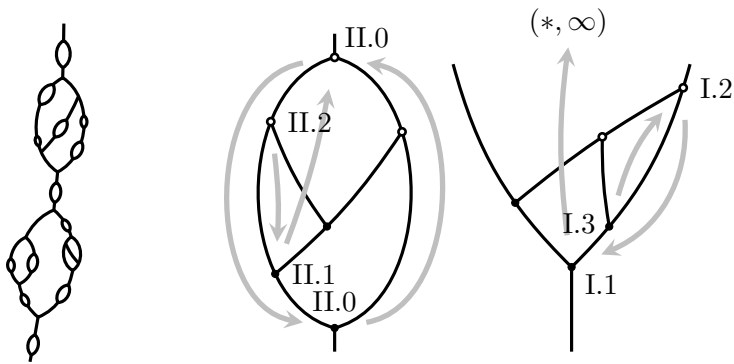
\begin{figure}[tcb]
\begin{center}
\begin{tikzpicture}
[>=stealth,mijn/.style={circle,outer sep=7pt},scale=0.45] 

\newcommand{\bubble}[4]
{\draw[green,opacity=0,name path=#2c] (#3) circle (#4);
\path[name intersections={of=#1 and #2c}];
\coordinate (#2t) at (intersection-1);
\coordinate (#2b) at (intersection-2);
\filldraw[very thick,fill=white,name path=#2] (#2t) to[bend right=60] (#2b)
 to[bend right=60] (#2t);
\node[red,opacity=0] at (#3) {#2};}

\newcommand{\cross}[4]
{\draw[green,opacity=0,name path=#2l] (#3) ++($-0.5*(#4)$)--++(#4);
\path[name intersections={of=#1 and #2l}];
\coordinate (#2t) at (intersection-1);
\coordinate (#2b) at (intersection-2);
\draw[very thick,name path=#2] (#2t) to[bend right=10] (#2b);
\node[red,opacity=0,below right] at (1,7) {#2};}

\newcommand{\lijnpunt}[4]
{\draw[green,opacity=0,name path=#2] (-5,#1)--++(10,0);
\path [name intersections={of=#3 and #2}];
\coordinate (#4) at (intersection-1);}

\begin{scope}[xshift=-9cm,yshift=-5cm,scale=1]
\draw[very thick,name path=A] (0,0) to[bend right=5] (1,10);
\bubble{A}{B}{0.5,2.9}{1.6}
\bubble{B}{C}{-0.3,3}{0.7}
\bubble{A}{D}{1,7}{1.4}
\bubble{C}{E}{0.2,2.9}{0.3}
\bubble{B}{F}{1,2}{0.3}
\bubble{B}{G}{1.2,3.2}{0.5}
\bubble{A}{H}{0.1,0.5}{0.3}
\bubble{A}{I}{0.8,5}{0.3}
\bubble{A}{J}{1,9.1}{0.4}
\cross{D}{K}{1,7.2}{2,3}
\cross{G}{L}{1.25,3.2}{2,-3}
\bubble{K}{M}{0.9,6.9}{0.3}
\bubble{D}{N}{0.5,7.9}{0.35}
\bubble{D}{O}{0.1,6.5}{0.2}
\bubble{D}{P}{1.3,6.3}{0.3}
\bubble{B}{Q}{-0.1,1.85}{0.2}
\bubble{B}{R}{1.1,4.05}{0.2}
\bubble{C}{S}{-0.6,2.85}{0.2}
\bubble{D}{T}{1.5,7.1}{0.2}
\end{scope}

\begin{scope}

\coordinate (A) at (0,4);
\coordinate (B) at (0,-4);
\draw[very thick] (A)--++(0,0.7);
\draw[very thick] (B)--++(0,-0.7);
\draw[very thick,name path=ABl] (A) to[bend right=75] (B);
\draw[very thick,name path=ABr] (B) to[bend right=75] (A);

\lijnpunt{-2.4}{hC}{ABl}{C}
\lijnpunt{1.8}{hD}{ABr}{D}

\draw[very thick,name path=CD] (C) to[bend right=8] (D);

\lijnpunt{-1}{hE}{CD}{E}
\lijnpunt{2.1}{hF}{ABl}{F}

\draw[very thick,name path=EF] (E) to[bend left=8] (F);

\node[mijn] (a) at (A) {};
\node[mijn] (b) at (B) {};
\node[mijn] (c) at (C) {};
\node[mijn] (f) at (F) {};

\draw[->,black!25,line width=2pt] (f) to[bend left=5] (c);
\draw[->,black!25,line width=2pt] (c) to[bend right=5] (a);
\draw[->,black!25,line width=2pt] (a) to[bend right=80] (b);
\draw[->,black!25,line width=2pt] (b) to[bend right=80] (a);

\node[above right] at (A) {$\rm II.0$};
\node[above] at ($(B)+(0,0.2)$) {$\rm II.0$};
\node[right] at ($(C)+(0.2,0)$) {$\rm II.1$};
\node[right] at ($(F)+(0.2,0)$) {$\rm II.2$};

\draw[line width=1pt,fill=white] (A) circle (3pt);
\filldraw (B) circle (3pt);
\filldraw (C) circle (3pt);
\draw[line width=1pt,fill=white] (D) circle (3pt);
\filldraw (E) circle (3pt);
\draw[line width=1pt,fill=white] (F) circle (3pt);
\end{scope}

\begin{scope}[xshift=7cm,yshift=-0.7cm]


\coordinate (A) at (0,-1.5);
\coordinate (B) at (-3.5,4.5);
\coordinate (C) at (3.5,4.5);

\draw[very thick] (A)--(0,-4);
\draw[very thick,name path=AB] (A) to[bend left=15] (B);
\draw[very thick,name path=AC] (A) to[bend right=15] (C);

\lijnpunt{0.4}{hD}{AB}{D}
\lijnpunt{3.8}{hE}{AC}{E}
\draw[very thick,name path=DE] (D) to[bend left=5] (E);

\lijnpunt{-0.3}{hF}{AC}{F}
\lijnpunt{2.35}{hG}{DE}{G}
\draw[very thick,name path=FG] (F) to[bend left=5] (G);

\draw[line width=1pt,fill=white] (G) circle (3pt);
\draw[line width=1pt,fill=white] (E) circle (3pt);
\filldraw (A) circle (3pt);
\filldraw (D) circle (3pt);
\filldraw (F) circle (3pt);

\node[mijn] (a) at (A) {};
\node[mijn] (b) at (B) {};
\node[mijn] (c) at (C) {};
\node[mijn] (d) at (D) {};
\node[mijn] (e) at (E) {};
\node[mijn] (f) at (F) {};
\node[mijn] (g) at (G) {};

\draw[->,black!25,line width=2pt] (f) to[bend left=12] (e);
\draw[->,black!25,line width=2pt] (e) to[bend left=35] (a);
\draw[->,black!25,line width=2pt] (a) to[bend left=7] (-0.1,5);

\node[left] at ($(F)+(-0.1,0.1)$) {$\rm I.3$};
\node[right] at ($(E)+(0.2,0)$) {$\rm I.2$};
\node[below right] at ($(A)+(0.1,0.1)$) {$\rm I.1$};
\node[above] at (-0.1,5) {$(\ast,\infty)$};

\end{scope}

\end{tikzpicture}
\end{center}
\caption{Left: schematic depiction of the fractal structure of the
  backbone. Middle and right: illustration of the classification of
  separation and meeting points on the backbone with the map $\phi$ (gray
  arrows). The leftmost picture contains two separation and two meeting points
  of type $\rm II.1$; all others are of type $\rm II.0$.}\label{fig:SMclass}
\end{figure}

It is clear that $N$ has some sort of fractal structure. The main reason for this is the ``infinite'' branching rate of the Brownian particles that $\xi_t$ consists of. Indeed, in each open time interval, a path in $\Ni(*,-\infty)$ will split infinitely often into two paths, which usually coalesce again after a very short time. But during the time that there are two paths, these too will for short moments split into two further paths, and so on, ad infinitum (see Figure~\ref{fig:SMclass}).

Natural questions are: 1.\ Is this the only way in which $N$ has a fractal structure? 2.\ How can we formulate this rigorously? And 3.\ Can we prove this? Here is an attempt to answer Question~2. Let $S$ and $M$ be the set of all {\em separation points} and {\em meeting points} (i.e., points of types $\rm (p,pp)_s$ and $\rm (pp,p)$, respectively, as defined in Section~\ref{S:webpts}) that lie on $N$. Set $S_\ast:=S\cup\{(\ast,-\infty)\}$ and $M_\ast:=M\cup\{(\ast,\infty)\}$. Define $\phi:S_\ast\to M_\ast$ by
\be
\phi(z):=\mbox{ the first meeting point of the left-most and right-most paths starting at $z$},
\ee
where we set $\phi(z):=(\ast,\infty)$ if the left-most and right-most paths starting at $z$ never meet, and for definitiveness we also define $\phi(\ast,-\infty):=(\ast,\infty)$.

By the symmetry of the backbone with respect to time reversal, we define $\phi:M_\ast\to S_\ast$ analogously, following paths in $\Ni(\ast,-\infty)$ backwards in time. (Note that this is different from our usual way of looking backwards. Instead of following dual paths, here we follow \emph{forward} paths in backbone backwards in time.)

A lot about the structure of the backbone can be understood in terms of the map $\phi$. For example, separation points $z$ such that $\phi(z)=(\ast,\infty)$ are $(-\infty,\infty)$-relevant separation points as defined in Section~\ref{S:relsep}. By Proposition~\ref{P:relev}, these form a locally finite subset of $\R^2$. If we have a pair of points $z\in S$ and $z'\in M$ such that $\phi(z)=z'$ and $\phi(z')=z$, then this corresponds to a path that splits at $z$ into two paths, which meet again at $z'$ (see Figure~\ref{fig:SMclass}). If $z=(x,t)$ and $z'=(x',t')$ form such a pair, then the left-most and right-most paths starting at $z$, up to their first meeting point $z'$, enclose a compact set that we will call a \emph{bubble}.

We conjecture that bubbles are the \emph{only} source of the fractal structure of $N$. The following two conjectures make this idea precise.

\bcon[Bubble hypothesis]\label{C:bubble}\ \\[-15pt]
\begin{itemize}
\item[{\bf(a)}] Almost surely, there does not exist an infinite sequence of points $(z_k)_{k\geq 0}$ in $S_\ast\cup M_\ast$, all different from each other, such that $\phi(z_{k-1})=z_k$ for all $k\geq 1$.
\item[{\bf(b)}] Almost surely, there does not exist an infinite sequence of points $(z_k)_{k\geq 0}$ in $S_\ast\cup M_\ast$, all different from each other, such that $\phi(z_k)=z_{k-1}$ for all $k\geq 1$.
\end{itemize}
\econ

\noindent
If Conjecture~\ref{C:bubble}~(a) is correct, then for each separation or meeting point $z$ that lies on the backbone, after some moment in the sequence $(z,\phi(z),\phi^2(z),\ldots)$, the same points start to repeat. Simple geometric considerations show that the limiting cycle must have length 2, i.e., there is an integer $n\geq 0$ such that $\phi^{n+2}(z)=\phi^n(z)$. Letting $n$ denote the smallest such integer, we say that $z$ is of \emph{class} ${\rm I.n}$ if $\phi^n(z)=(\ast,\infty)$ or $(\ast,-\infty)$ and of \emph{class} ${\rm II.n}$ if $\phi^n(z)\in\R^2$. Note that separation and meeting points of type ${\rm II.0}$ are the bottom and top points of bubbles, as we have just defined them. Points of type ${\rm II.n}$ all lie inside bubbles, with larger $n$ leading to a more complex left-right crossover pattern in the containing bubble (see Figure~\ref{fig:SMclass}). The following conjecture says that most bubbles have a very simple internal structure.

\bcon[Bubble complexity]\label{C:bubble2}
Let $\Ci_{\rm I,n}$, resp.\ $\Ci_{\rm II,n}$, denote the set of meeting and separation points on the backbone $N$ which are of class ${\rm I.n}$, resp.\ ${\rm II.n}$. Then
\begin{itemize}
\item[{\bf(a)}] The sets $\Ci_{\rm II,0}$ and $\Ci_{\rm II,1}$ are dense in $N$.
\item[{\bf(b)}] The sets $\Ci_{\rm II,n}$ with $n\geq 2$ and $\Ci_{\rm I,n}$ with $n\geq 1$ are locally finite subsets of $\R^2$.
\end{itemize}
\econ

As motivation for these conjectures, we state a somewhat more applied problem. Modulo a time reversal, the set $\Wi(\{0\}\times\R)$ of all paths in the Brownian web can be interpreted as the spatial genealogies of a population living in one-dimensional space (see~Sec.~\ref{SS:voter}). Let $M$ be the set of all meeting points of paths in $\Wi(\{0\}\times\R)$, i.e.,
\be\label{MRCA}
M:=\big\{(x,t)\in\R^2:t>0,\exists\, \pi,\pi'\in\Wi(\{0\}\times\R)
\mbox{ s.t.\ }\pi(t)=x=\pi'(t),\ \pi<\pi'\mbox{ on }(0,t)\big\}.
\ee
Biologically, we can interpret a point $z=(x,t)$ that is the first meeting point of paths $\pi_{(x,0)}$ and $\pi_{(x',0)}$ as the \emph{most recent common ancestor} (MRCA) of $x$ and $x'$, that lived a time $t$ in the past. It is easy to see that $M$ is a locally finite subset of $\R\times(0,\infty)$.

Moving away from neutral evolution, one can interpret paths in the Brownian net as \emph{potential} genealogies (cf.~the biased voter model in Sec.~\ref{SS:bvoter}), where in order to determine the true genealogy of an individual, one has to have information about what happens at selection events, which correspond to separation points. In this case, replacing $\Wi(\{0\}\times\R)$ in (\ref{MRCA}) by $\Ni(\{0\}\times\R)$, we can interpret the resulting set $M$ as the set of \emph{potential most recent common ancestors} (PMRCAs) of individuals living at time zero.

\bcon[Potential most recent common ancestors]\label{C:PMRCA}
If we replace $\Wi(\{0\}\times\R)$ in \eqref{MRCA} by $\Ni(\{0\}\times\R)$, then the set $M$ is a locally finite subset of $\R\times(0,\infty)$.
\econ

The definition of a potential most recent common ancestor (PMRCA) is somewhat reminiscent of the definition of a $(0,\infty)$-relevant separation point, except that in order to determine wether a meeting point is a PMRCA, we have to follow forward paths backwards in time. This is similar to the definition of the map $\phi:M_\ast\to S_\ast$ and indeed, Conjecture~\ref{C:PMRCA} seems to be closely related to our previous two conjectures. In particular, meeting points that lie inside a bubble (properly defined w.r.t.\ the finite time horizon $0$) can never be PMRCAs.

\subsection{Miscellaneous open questions}

We collect below some other interesting questions.

\begin{itemize}

\item[(1)] Is the Brownian net integrable in the sense that the branching-coalescing point set defined in \eqref{bcptset} admits an explicit characterization, similar to the Pfaffian point process characterization of the coalescing point set discussed in Remark~\ref{R:pfaffian}?

\item[(2)] Find effective criteria for general branching-coalescing particle systems, where paths can cross, to converge to the Brownian net.
This would extend Theorem~\ref{T:netconvD} for the non-crossing case, and pave the way for the study of voter model perturbations discussed in Section \ref{SS:perturb}.

\item[(3)] Can one formulate a well-posed martingale problem for the branching-coalescing point set $(\xi_t)_{t\geq 0}$ defined in \eqref{bcptset}? This may offer an alternative route to prove convergence to the Brownian net, which is based on
generator convergence and avoids the paths topology.

\item[(4)] For the Howitt-Warren flow $(K^\uparrow_{s,t}(x, {\rm d}y))_{s<t, x\in\R}$ introduced in Section~\ref{S:HW}, which gives the transition probability kernels of a random motion in a continuum space-time random environment constructed from the Brownian web and net, can one show that for any $x\in\R$, $-\log K^\uparrow_{0,t}(0, [tx, tx+1])$ has Tracy-Widom GUE fluctuations on the scale $t^{1/3}$
as $t\to\infty$,  similar to recent results in \cite{BC15} for random walks on $\Z$ in special i.i.d.\ Beta-distributed space-time random environments (cf.\ Remark~\ref{R:HWfluc})?

\item[(5)] In \cite{EMS13}, Evans, Morris and Sen showed that coalescing stable L\'evy processes on $\R$ with stable index $\alpha\in (1,2)$ come down from infinity, i.e., starting from everywhere on $\R$, the coalescing L\'evy processes become locally finite on $\R$ at any time $t>0$. This suggests the existence  of  a family of {\em L\'evy webs}. In fact, using the tube topology reviewed in Section~\ref{SS:tube} and the property of coming down from infinity, it should be straightforward to construct the L\'evy web w.r.t.\ the tube topology. For coalescing Brownian flow on the Sierpinski gasket, this was carried out in \cite{BGS13}. Can one construct the L\'evy web in the paths topology, and what type of special points may arise in the spirit of Theorem~\ref{T:classweb}?

By scaling invariance, the density of the coalescing stable L\'evy flow on $\R$ with index $\alpha \in (1,2)$ should be $\rho_t = \frac{C}{t^{1/\alpha}}$.  Identifying $C$ would allow one to determine the sharp asymptotic rate of decay for the density of coalescing random walks on $\Z$ in the domain of attraction of a stable L\'evy process, similar to \cite[Corollary 7.1]{NRS05} and \cite{SS12} in the Brownian case.

\item[(6)] The Brownian web has been shown by Ellis and Feldheim~\cite{EF12} to be a two-dimensional black noise, equipped with a family of sub-$\sigma$-fields indexed by finite unions of open rectangles in $\R^2$. Is there a maximal extension of the family of sub-$\sigma$-fields, indexed by a Boolean algebra of domains in $\R^2$ that include the rectangles, such that the Brownian web remains a noise (cf.\ Remark~\ref{R:factorize})? A similar question is also open for the scaling limit of critical planar percolation~\cite[Sec.~1.6]{T14}.

Recently Tsirelson~\cite{T13} showed by general arguments that there exists some continuous path $f$ starting at time $-\infty$, such that
$\Fi_{f_-}$ and $\Fi_{f_+}$, the $\sigma$-fields generated by Brownian web paths restricted respectively to the left or the right of the graph of $f$, do not jointly generate the full $\sigma$-field $\Fi$. Can we characterize the set of $f$ for which the factorization property $\Fi = \Fi_{f_-} \vee \Fi_{f_+}$ holds? (The case $f\equiv 0$ has been treated in~\cite{EF12}).

\item[(7)] We expect the Brownian net with killing (see Section~\ref{SS:netkill}) to belong to the same universality class as oriented percolation (OP) on $\Z^{1+1}$, which corresponds to the so-called Reggeon field theory~\cite{CS80, M78} (see also \cite{H00} for an extensive survey on the OP universality class). The Brownian net without killing belongs to a different universality class. But is there also a field theory corresponding to the Brownian net?  A positive indication is that the Brownian net has competing effects of instantaneous coalescence vs infinite effective rate of branching, which is very much in the spirit of renormalizations in field theories.

\item[(8)] The Brownian web appears in the scaling limit of super-critical oriented percolation, as discussed in Section \ref{SS:supperc}, and the Brownian net is expected to arise if the percolation parameter $p$ is allowed to vary in a small interval. However, the most interesting question is: What is the scaling limit of critical (and near-critical) oriented percolation (OP) on $\Z^{1+1}$? The Brownian net with killing provides a simpler model to study the OP universality class. We expect that in the near-critical scaling limit, one would obtain a family of models which interpolates between the Brownian net and the scaling limit of critical oriented percolation. However, such a goal appears far out of reach at the moment, because there are no conjectures at all on what the critical scaling limit might be, and neither have critical exponents been shown to exist.

\end{itemize}

\bigskip
\noi
{\bf\large Acknowledgement} This article is based on lectures given by R.~Sun in the trimestre program {\em Disordered Systems, Random Spatial Processes and Some Applications}, at the Institute Henri Poincar\'e in the spring of 2015. We thank the program and the organizers for the opportunity to lecture on this topic. We thank all our collaborators and colleagues with whom we have had valuable discussions. In particular, R.~Sun and E.~Schertzer would like to thank Chuck Newman and Krishnamurthi Ravishankar for introducing them to this topic. R.~Sun also thanks University of Warwick, Leiden University, and Academia Sinica for hospitality, where part of these notes were written. R.~Sun is supported by AcRF Tier 1 grant R-146-000-185-112. J.M.~Swart is sponsored by GACR grant 15-08819S.


\begin{thebibliography}{BMSV06}

\bibitem[A79]{A79}
R.~Arratia.
Coalescing Brownian motions on the line.
Ph.D. Thesis, University of Wisconsin, Madison, 1979.

\bibitem[A81]{A81}
R.~Arratia.
Coalescing Brownian motions and the voter model on $\Z$.
Unpublished partial manuscript.

\bibitem[A98]{A98}
M. Aizenman.
Scaling limit for the incipient spanning clusters.
{\em Mathematics of Multiscale Materials: Percolation and Composites (Minneapolis, Minn., 1995-1996)}, 1--24,  IMA Vol. Math. Appl. 99, Springer, New York, 1998.

\bibitem[AB99]{AB99}
M. Aizenman and A. Burchard.
H\"older regularity and dimension bounds for random curves.
{\em Duke Math. J.} 99, 419--453, 1999.

\bibitem[ARS08]{ARS08}
S. Athreya, R. Roy and A. Sarkar.
Random directed trees and forest -- drainage networks with dependence.
{\em Electron. J. Probab.} 13,  2160--2189, 2008.

\bibitem[AS11]{AS11}
S.~Athreya and R.~Sun.
One-dimensional voter model interface revisited.
{\em Electron.\ Commun.\ Probab.} 16, 792--800, 2011.

\bibitem[AS13]{AS13}
S.~Athreya and J.M.~Swart.
Systems of branching, annihilating, and coalescing particles.
{\em Electron. J. Probab.} 17, 1--32, 2012.


\bibitem[bABD90]{bABD90}
D. ben-Avraham, M.A. Burschka, and C.R. Doering.
Statics and dynamics of a diffusion-limited reaction: anomalous kinetics, nonequilibrium self-ordering, and a dynamic transition.
{\em J.  Stat. Phys.} 60, 695--728, 1990.

\bibitem[B99]{B99}
P.~Billingsley.
{\em Convergence of probability measures}, 2nd edition.
John Wiley \& Sons, 1999.


\bibitem[B09]{B09}
N.~Berestycki.
Recent progress in coalescent theory. {\em Ensaios Matemáticos} 16.
Sociedade Brasileira de Matem\'atica, Rio de Janeiro, 2009.

\bibitem[BB07]{BB07}
F. Bacelli and C. Bordenave.
The radial spanning tree of a Poisson point process.
{\em Ann. Appl. Probab.} 17, 305--359, 2007.

\bibitem[BC15]{BC15}
G. Barraquand and I. Corwin.
Random-walk in Beta-distributed random environment.
arXiv:1503.04117, 2015.


\bibitem[BGS15]{BGS13}
N. Berestycki, C. Garban and A. Sen.
Coalescing Brownian flows: a new approach.
{\em Ann. Probab.} 43, 3177--3215, 2015.

\bibitem[BHPS03]{BHPS03}
I. Benjamini, O. H\'aaggstrom, Y. Peres, J. E. Steif.
Which properties of a random sequence are dynamically sensitive?
{\em Ann. Probab.} 31, 1--34, 2003.

\bibitem[BMSV06]{BMSV06}
S.~Belhaouari, T.~Mountford, R.~Sun, and G.~Valle.
Convergence results and sharp estimates for the voter model interfaces.
{\em Electron.\ J.\ Probab.} 11, 768--801, 2006.

\bibitem[CDF09]{CDF09}
C.F.~Coletti, E.S.~Dias, and L.R.G.~Fontes.
Scaling limit for a drainage network model.
{\em J.\ Appl.\ Probab.}~46, 1184--1197, 2009.

\bibitem[CDP13]{CDP13}
J.T. ~Cox, R. ~Durrett, E. ~Perkins.
Voter model perturbations and reaction diffusion equations.
{\em Ast\'erisque} 349, 2013.


\bibitem[CFN06]{CFN06}
F.~Camia, L.R.G.~Fontes, and C.M.~Newman.
Two-dimensional scaling limits via marked nonsimple loops.
{\em Bull. Braz. Math. Soc. (N.S.)} 37, 537--559, 2006.

\bibitem[CS80]{CS80}
J.L. Cardy and R.L. Sugar.
Directed percolation and Reggeon field theory.
{\em J. Phys. A: Math. Gen.} 13, L423--L427, 1980.

\bibitem[CT13]{CT13}
D. Coupier and V.C. Tran.
The 2D-directed spanning forest is almost surely a tree.
{\em Random Structures $\&$ Algorithms} 42, 59--72, 2013.


\bibitem[CV11]{CV11}
C.F.~Coletti and G.~Valle.
Convergence to the Brownian Web for a generalization of the drainage network model.
{\em Ann. Inst. H. Poincaré Probab. Statist.} 50, 899--919, 2014.

\bibitem[CV14]{CV14}
C.F. Coletti and L.A. Valencia.
Scaling limit for a family of random paths with radial behavior.
arXiv:1310.6929, 2014.

\bibitem[D84]{D84}
R. Durrett.
Oriented percolation in two dimensions.
{\em Ann. Probab.} 12, 999--1040,1984.

\bibitem[D06]{D06}
D. Dhar.
Theoretical studies of self-organized criticality.
{\em Physica A: Statistical Mechanics and its Applications} 369, 29--70, 2006.

\bibitem[D12]{D12}
L. Dumaz.
A clever (self-repelling) burglar.
{\em Electron. J. Probab.} 17, 1--17, 2012

\bibitem[D15]{D15}
D. Dhar.
Fragmentation of a sheet by propagating, branching and merging cracks.
{\em Journal of Physics A: Mathematical and Theoretical} 48, 17, 2015.

\bibitem[DbA88]{DbA88}
C.R. Doering and D. ben-Avraham.
Interparticle distribution functions and rate equations for diffusion-limited reactions.
{\em Phys. Rev. A} 38, 3035, 1988.

\bibitem[DEF$^+$00]{DEF+00}
P. Donnelly, S.N. Evans, K. Fleischmann, T.G. Kurtz, and X. Zhou.
Continuum-sites stepping-stone models, coalescing exchangeable partitions and random trees.
{\em Ann. Probab.} 28, 1063--1110, 2000.

\bibitem[DGP11]{DGP11}
A. Depperschmidt, A. Greven, and P. Pfaelhuber.
Marked metric measure spaces.
{\em Electron. Commun. Probab.} 16, 174--188, 2011.

\bibitem[DW09]{DW09}
M. Damron and C.L. Winter.
A non-Markovian model of rill erosion.
{\em Networks and Heterogeneous Media}  4, 731--753, 2009.

\bibitem[E10]{E10}
T. Ellis.
{\em Coalescing stochastic flows driven by Poisson random measure and convergence to the Brownian web}.
PhD thesis, University of Cambridge, 2010.

\bibitem[EF16]{EF12}
T. Ellis and O.N. Feldheim.
The Brownian web is a two-dimensional black noise.
{\em Ann. Inst. H. Poincar\'e Probab. Statist.} 52,  162--172, 2016.

\bibitem[EFS15]{EFS15}
A. Etheridge, N. Freeman and D. Straulino.
The Brownian Net and Selection in the Spatial Lambda-Fleming-Viot Process.
arXiv:1506.01158, 2015.

\bibitem[EMS13]{EMS13}
S.N.~Evans, B.~Morris and A.~Sen.
Coalescing systems of non-Brownian particles.
{\em Probab.\ Theory Related Fields}~156, 307--342, 2013.

\bibitem[FFW05]{FFW05}
P.A.~Ferrari, L.R.G.~Fontes, and X.-Y.~Wu.
Two-dimensional Poisson Trees converge to the Brownian web.
{\em Ann.\ Inst.\ H.\ Poincar\'e Probab.\ Statist.} 41, 851--858, 2005.

\bibitem[FINS01]{FINS01}
L.R.G.~Fontes, M.~Isopi, C.M.~Newman, D.~Stein.
Aging in 1D discrete spin models and equivalent systems.
{\em Physical Review Letters} 87, 110201, 2001.

\bibitem[FINR02]{FINR02}
L.R.G.~Fontes, M.~Isopi, C.M.~Newman, K.~Ravishankar.
The Brownian web.
{\em Proc.~Nat.~Acad.~Sciences} 99, 15888--15893, 2002.

\bibitem[FINR03]{FINR03}
L.R.G.~Fontes, M.~Isopi, C.M.~Newman, K.~Ravishankar.
The Brownian web: characterization and convergence.
math/0304119, 2003.

\bibitem[FINR04]{FINR04}
L.R.G.~Fontes, M.~Isopi, C.M.~Newman, K.~Ravishankar.
The Brownian web: characterization and convergence.
{\em Ann.\ Probab.}~32(4), 2857--2883, 2004.

\bibitem[FINR06]{FINR06}
L.R.G.~Fontes, M.~Isopi, C.M.~Newman, K.~Ravishankar.
Coarsening, nucleation, and the marked Brownian web.
{\em Ann. Inst. H. Poincar\'e Probab. Statist.} 42, 37--60, 2006.

\bibitem[FKG70]{FKG70}
C.M.~Fortuin, P.W.~Kasteleyn, J.~Ginibre.
Correlation inequalities on some partially ordered sets.
{\em Commun.~Math.~Phys.} 22, 89--103, 1970.

\bibitem[FLT04]{FLT04}
P.A. Ferrari, C. Landim and H. Thorisson.
Poisson trees, succession lines and coalescing random walks.
{\em Ann. Inst. H. Poincar\'e Probab. Statist.} 40, 141--152, 2004.

\bibitem[FN06]{FN06}
L.R.G.~Fontes and C.M.~Newman.
The full Brownian web as scaling limit of stochastic flows.
{\em Stoch. Dyn.} 6, 213–-228, 2006.

\bibitem[FNRS09]{FNRS09} L. R. G. Fontes, C. M. Newman, K. Ravishankar and E. Schertzer.
Exceptional times for the dynamical discrete web.
{\em Stoch. Proc.  Appl.} 119, 2832-2858 (2009).

\bibitem[FVV15]{FVV14}
L.R.G. Fontes, L.A. Valencia, and G. Valle.
Scaling limit of the radial Poissonian web.
{\em Electron. J. Probab.} 20 , paper no.\ 31, 2015.

\bibitem[G99]{G99}
G. Grimmett, Percolation, Springer-Verlag, Berlin--Heidelberg--New York, 1999.

\bibitem[GPS13a]{GPS13a}
C.~Garban, G.~Pete, and O.~Schramm.
Pivotal, cluster and interface measures for critical planar percolation.
{\em J. Amer. Math. Soc.} 26, 939--1024, 2013.

\bibitem[GPS13b]{GPS13b}
C.~Garban, G.~Pete, and O.~Schramm.
The scaling limits of near-critical and dynamical percolation.
arXiv:1305.5526, 2013.


\bibitem[GRS04]{GRS04}
S. Gangopadhyay, R. Roy and A. Sarkar.
Random Oriented Trees: A Model of Drainage Networks.
{\em Ann. App. Probab.} 14, 1242--1266, 2004.

\bibitem[GS12]{GS12}
C. Garban and J. Steiff.
Noise sensitivity and percolation.
{\em Probability and statistical physics in two and more dimensions}, 49--154,
Clay Math. Proc., 15, Amer. Math. Soc., Providence, RI, 2012.

\bibitem[GSW15]{GSW15}
A.~Greven, R.~Sun, and A.~Winter.
Continuum space limit of the genealogies of interacting Fleming-Viot processes on $\Z$.
arXiv:1508.07169, 2015.

\bibitem[H00]{H00}
H. Hinrichsen.
Non-equilibrium critical phenomena and phase transitions into absorbing states.
{\em Advances in Physics} 49, 815--958, 2000.

 \bibitem[HL75]{HL75}
R.~Holley, T.~Liggett.
Ergodic theorems for weakly interacting infinite systems and the voter model.
{\em Annals of Probab.} {\bf 3}, 643-663, 1975.


\bibitem[HW09a]{HW09a}
C.~Howitt and J.~Warren.
Consistent families of Brownian motions and stochastic flows of kernels.
{\em Ann.\ Probab.}~37, 1237--1272, 2009.

\bibitem[HW09b]{HW09b}
C.~Howitt and J.~Warren.
Dynamics for the Brownian web and the erosion flow.
{\em Stochastic Processes Appl.}~119, 2028--2051, 2009.

\bibitem[J11]{J11}
D. Jenkins.
Superdiffusive and subdiffusive exceptional times in the dynamical discrete web.
{\em Stochastic Processes Appl.}~125, pp. 3373--3400, 2015.

\bibitem[K89]{K89}
T. Kuczek.
The central limit theorem for the right edge of supercritical oriented percolation.
{\em Ann. Probab.} 17, 1322--1332, 1989.

\bibitem[KS91]{KS91}
I.~Karatzas, S.E.~Shreve.
{\em Brownian Motion and Stochastic Calculus.}
2nd edition, Springer-Verlag, New York, 1991.

\bibitem[KPWH95]{KPWH95}
K. Krebs, M.P. Pfannm\"uller, B. Wehefritz, and H. Hinrichsen.
Finite-Size Scaling Studies of One-Dimensional Reaction-Diffusion Systems. Part I. Analytical Results.
{\em J. Stat. Phys.} 78, 1429--1470, 1995.

\bibitem[L04]{L04}
T.~Liggett.
{\em Interacting particle systems}. Springer, Berlin-Heidelberg-New York, 1975.

\bibitem[LL04]{LL04}
Y.~Le Jan and S.~Lemaire.
Products of beta matrices and sticky flows.
{\em Probab.\ Th.\ Relat.\ Fields}~130, 109--134, 2004.

\bibitem[LR04a]{LR04}
Y.~Le Jan and O.~Raimond.
Flows, Coalecence and Noise.
{\em Annals of Probab.}~32, 1247--1315, 2004.

\bibitem[LR04b]{LR04b}
Y.~Le Jan and O.~Raimond.
Sticky flows on the circle and their noises.
{\em Probab.\ Th.\ Relat.\ Fields}~129, 63--82, 2004.

\bibitem[M78]{M78}
M. Moshe.
Recent developments in Reggeon field theory.
{\em Physics Reports} 37, 255--345, 1978.

\bibitem[M05]{M05}
I. Molchanov.
{\em Theory of random sets}.
Springer, Berlin, 2005.


\bibitem[MDD$^+$99]{MDDGL99}
J. Molofsky, R. Durrett, J. Dushoff, D. Griffeath, S. Levin,
Local frequency dependence and global coexistence.
{\em Theor. Pop. Biol.}~{55},  270--282, 1999.


\bibitem[MNR13]{MNR13}
Y. Mohylevskyy, C. M. Newman, K. Ravishankar.
Ergodicity and percolation for variants of one-dimensional voter models.
{\em ALEA.} {\bf 10}, 485--504, 2013.

\bibitem[MRTZ06]{MRTZ06}
R. Munasinghe, R. Rajesh, R. Tribe, O. Zaboronski.
Multi-scaling of the $n$-point density function for coalescing Brownian motions.
{\em Commun. Math. Phys.} 268, 717--725, 2006.

\bibitem[NP99]{NP99}
C. Neuhauser,  S.W. Pacala,
An explicitly spatial version of the Lotka-Volterra model with interspecific competition.
{\em Ann. Appl. Probab.} {\bf 9},
1226---1259, 2000.


\bibitem[NR06]{NR06}
C.M.~Newman and K.~Ravishankar.
Convergence of the T\'oth lattice filling curve to the T\'oth-Werner plane filling curve.
{\em ALEA Lat. Am. J. Probab. Math. Stat.} 1, 333--345, 2006.

\bibitem[NRS05]{NRS05}
C.M.~Newman, K.~Ravishankar, R.~Sun.
Convergence of coalescing nonsimple random walks to the Brownian web.
{\em Electron.\ J.\ Prob.} 10, 21--60, 2005.


\bibitem[NRS10]{NRS10}
C.M.~Newman, K.~Ravishankar, and E.~Schertzer.
Marking $(1,2)$ points of the Brownian web and applications.
{\em Ann.\ Inst.\ Henri Poincar\'e Probab.\ Statist.}~46, 537--574, 2010.

\bibitem[NRS15]{NRS15}
C.M.~Newman, K.~Ravishankar, and E.~Schertzer.
The Brownian Net with killing.
{\em Stoch. Proc. and App.}~125, 1148Ð1194, 2015.

\bibitem[NS93]{NS93}
C. Neuhauser and A. Sudbury.
The biased annihilating branching process.
{\em Adv. Appl. Prob.} 25, 24--38, 1993.


\bibitem[NT12]{NT12}
J.~Norris and A.~Turner.
Hastings-Levitov aggregation in the small-particle limit.
{\em Comm. Math. Phys.} 316, 809-–841, 2012.

\bibitem[NT15]{NT15}
J.~Norris and A.~Turner.
Weak convergence of the localized disturbance flow to the coalescing Brownian flow.
{\em Ann. Probab.} 43, 935--970, 2015.

\bibitem[OHLN06]{OHLN06}
 H. Ohtsuki, H. C. Hauert, E. Lieberman, M.A. Nowak.
 A simple rule for the evolution of cooperation on graphs and social networks.
 {\it Nature}. {441}, 502---505, 2006.

\bibitem[P97]{P97}
V. Privman (editor).
{\em Nonequilibrium statistical mechanics in one dimension.}
Cambridge University Press, Cambridge, 1997.

\bibitem[RSS16a]{RSS13}
R. Roy, K. Saha, and A. Sarkar.
Random directed forest and the Brownian web.
{\em Ann. Inst. H. Poincar\'e Probab. Statist.} 52, 1106--1143, 2016.

\bibitem[RSS16b]{RSS15}
R. Roy, K. Saha, and A. Sarkar.
Hack's law in a drainage network model: a Brownian web approach.
{\em Ann. Appl. Probab.} 26, 1807--1836, 2016.


\bibitem[S76]{S76}
F.~Spitzer.
{\em Principles of Random Walk}, 2nd edition.
Springer-Verlag, 1976.

\bibitem[S77]{S77}
D. Schwartz.
Applications of duality to a class of Markov processes.
{\em Ann. Probab.} 5, 522--532, 1977.

\bibitem[S09]{S09}
J. Steif.
A survey of dynamical percolation.
{\em Fractal geometry and stochastics} IV,
Progress in Probability 61, 145--174, 2009.


\bibitem[S10]{S10}
T. Sepp\"al\"ainen.
{\em Current fluctuations for stochastic particle systems with drift in one spatial dimension}.
Ensaios Matem\'aticos [Mathematical Surveys], 18.  Sociedade Brasileira de Matem\'atica, Rio de Janeiro, 2010.

\bibitem[SS08]{SS08}
R.~Sun and J.M.~Swart.
The Brownian net.
{\em Ann.\ Probab.}~36, 1153-1208, 2008.


\bibitem[SS11]{SS11}
O. Schramm and S. Smirnov.
On the scaling limits of planar percolation.
{\em Ann. Probab.} 39, 1768--1814, 2011. With an appendix by C. Garban.

\bibitem[SS12]{SS12}
A. Sarkar and R. Sun.
Brownian web and oriented percolation: density bounds.
RIMS Kokyuroku, No. 1805, {\em Applications of Renormalization Group Methods in Mathematical Sciences}, 90--101, 2012.

\bibitem[SS13]{SS13}
A.~Sarkar and R.~Sun.
Brownian web in the scaling limit of supercritical oriented percolation in dimension $1+1$.
{\em Electron.\ J.\ Probab.} 18, paper 21, 2013.

\bibitem[SS15]{SS15}
K. Saha and A. Sarkar.
Convergence of drainage networks to the Brownian web.
arXiv:1508.06919, 2015.

\bibitem[SSS09]{SSS09}
E.~Schertzer, R.~Sun, and J.M.~Swart.
Special points of the Brownian net.
{\em Electron.\ J.\ Prob.}~14, Paper 30, 805--864, 2009.


\bibitem[SSS14]{SSS14}
E.~Schertzer, R.~Sun, and J.M.~Swart.
Stochastic flows in the Brownian web and net.
{\em Mem. Amer. Math. Soc.} 227, no.~1065, 2014.

\bibitem[STW00]{STW00}
F.~Soucaliuc, B.~T\'oth, W.~Werner.
Reflection and coalescence between one-dimensional Brownian paths.
{\em Ann.\ Inst.\ Henri Poincaré Probab.\ Statist.}~36, 509--536, 2000.

\bibitem[T95]{T95}
B.~T\'oth.
The ``true'' self-avoiding walk with bond repulsion on $\Z$ : Limit theorems.
{\em Ann. Probab.} 23, 1523--1556, 1995.

\bibitem[T04a]{T04a}
B. Tsirelson.
Scaling Limit, Noise, Stability.
{\em Lecture Notes in Mathematics} 1840, 1--106, Springer, 2004.

\bibitem[T04b]{T04b}
B. Tsirelson.
Nonclassical stochastic flows and continuous products.
{\em Probability Surveys} 1, 173--298, 2004.

\bibitem[T13]{T13}
B. Tsirelson.
Random compact set meets the graph of nonrandom continuous function.
arXiv:1308.5112, 2013.

\bibitem[T14]{T14}
B. Tsirelson.
Noise as a Boolean algebra of $\sigma$-fields.
{\em Ann. Probab.} 42, 311--353, 2014.

\bibitem[TW98]{TW98}
B.~T\'oth and W.~Werner.
The true self-repelling motion.
{\em Probab.\ Theory Related Fields}~111, 375--452, 1998.

\bibitem[TYZ12]{TYZ12}
R.~Tribe, S.K.~Yip, and O.~Zaboronski.
One dimensional annihilating and coalescing particle systems as extended Pfaffian point processes.
{\em  Electron. Commun. Probab.} 17, no.~40, 2012.

\bibitem[TZ11]{TZ11}
R.~Tribe and O.~Zaboronski.
Pfaffian formulae for one dimensional coalescing and annihilating systems.
{\em Electron. J. Probab.} 16, no.~76, 2080–-2103, 2011.

\bibitem[W02]{W02}
J.~Warren.
The noise made by a Poisson snake.
{\em Electron.\ J.\ Probab.}~7, Paper No.~21, 1--21, 2002.


\bibitem[W07]{W07}
J.~Warren.
Dyson's Brownian motions, intertwining and interlacing.
{\em  Electron. J. Probab.} 12, no. 19, 573--590, 2007.

\bibitem[WB71]{WB71}
T. Williams and R. Bjerknes.
A stochastic model for the spread of an abnormal clone through the basal layer of the epithelium.
In {\em Symp. Tobacco Research Council}, London, 1971.

\bibitem[WZ08]{WZ08}
X.-Y. Wu  and  Y. Zhang.
A  geometrical  structure  for  an  infinite  oriented  cluster  and  its uniqueness.
{\em Ann. Probab.} 36, 862--875, 2008.

\bibitem[Y16]{Y14}
Jinjiong Yu.
Edwards-Wilkinson fluctuations in the Howitt-Warren flows.
{\em Stochastic Processes Appl.} 126, 948--982, 2016.

\end{thebibliography}
\end{document}